\documentclass[12pt]{article} 
\usepackage{amssymb,latexsym}
\usepackage{amsfonts}   
\input cyracc.def
\font\tencyr=wncyr10 at 12truept
\def\cyr{\tencyr\cyracc}
 at 12truept

 at 15truept

 at 12truept
 at 12truept

 at 12truept

\font\germ=eufm10  
\def\qed{\hfill $\vrule height 2.5mm  width 2.5mm depth 0mm $}

\newtheorem{theorem}{Theorem}[section]
\newtheorem{pr}[theorem]{Proposition}
\newtheorem{cor}[theorem]{Corollary}
\newtheorem{de}[theorem]{Definition}
\newtheorem{con}[theorem]{Conjecture}
\newtheorem{rem}[theorem]{Remark}
\newtheorem{lem}[theorem]{Lemma}

\newtheorem{prb}[theorem]{Problem}
\newtheorem{ex}[theorem]{Example}
\newtheorem{exs}[theorem]{Examples}
\newtheorem{quest}[theorem]{Question}
\def\s{\hbox{\germ S}}
\def\g{\hbox{\germ g}}

\def\M{\hbox{\germ M}}

\def\R{\mathbb{R}}
\def\N{\mathbb{N}}
\def\Z{\mathbb{Z}}
\def\C{\mathbb{C}}
\def\Q{\mathbb{Q}}

\def\ds{\displaystyle}
\def\wt{\widetilde}
\def\wh{\widehat}
\def\ld{\lambda}

 \def\al{\alpha}
\def\goti#1{#1\llap{$#1$\hskip.3pt}\llap{$#1$\hskip.5pt}
   \llap{$#1$\hskip.5pt}}
\def\ggoti#1{\goti{\goti#1}}
\def\blambda{\ggoti\lambda}
\def\bmu{\ggoti\mu}
\def\bnu{\ggoti\nu}
\def\sbmu{\ggoti{\scriptstyle\mu}}
\def\sbnu{\ggoti{\scriptstyle\nu}}
\def\sblambda{\ggoti{\scriptstyle\lambda}}
\def\bdoteq{\buildrel\bullet\over{=\!\!\!=}}
\begin{document}
\begin{center}
{\textbf{AN INVITATION TO THE GENERALIZED SATURATION CONJECTURE}}
\end{center}

\begin{center}
\textsc{anatol n.\ kirillov}
\end{center}
\begin{center}
{\small {\it Research Institute of Mathematical Sciences ( RIMS )}} \\
{\small {\it Kyoto 606-8502, Japan }}\\ {\small}
\vspace{5mm} 
\begin{flushright} \small\it{\cyr ``Ya ~vam ~ne ~skazhu ~za ~vsyu ~Odessu,}

{\cyr    ~vsya ~Odessa ~ochen\cprime ~~velika.''}

From a famous russian song.
\end{flushright}
\vspace{2mm}
\begin{flushright} \small {``I'm ~not ~intending ~to ~tell ~{\it all} ~about ~Tokyo,

~~The ~whole ~Tokyo ~is ~too ~big.'' }
\end{flushright}
\vspace{4mm}
{\it Communicated by T.~Kawai. ~Received March~2, 2004. Revised May 17, 2004.}

\end{center}

\vskip 0.3cm
\begin{abstract}

 We report about some results,  interesting examples,  problems and 
conjectures revolving around the parabolic Kostant partition functions,
the parabolic Kostka polynomials and ``saturation'' properties  of several
generalizations  of the Littlewood--Richardson numbers. The Contents contains
the titles of main topics we are going to discuss in the present paper.
\end{abstract}
\vskip10cm


\vskip 0.5cm
\pagebreak
\begin{center}
\textsc{Contents}
\end{center}
\nopagebreak
\begin{tabular}{rrl}
~~~~& 1 &  {\bf Introduction}
\\
& 2 & {\bf Basic definitions and notation}
\\
&   & 2.1~~Compositions and partitions
\\
&   & 2.2~~Kostka--Foulkes polynomials 
\\
&   & 2.3~~Skew Kostka--Foulkes polynomials
\\
&   & 2.4~~Littlewood--Richardson numbers and Saturation Theorem
\\
&   & 2.5~~Internal product of Schur functions, plethysm, and polynomials  
$L_{\al,\beta}^{\mu}(q)$ 
\\
&   & 2.6~~Extended and Restricted Littlewood--Richardson numbers
\\
& 3 & {\bf Parabolic Kostant partition function and its $q$-analog}
\\
&   & 3.1~~Definitions: algebraic and combinatorial
\\
&   & 3.2~~Elementary properties, and explicit formulas for $l(\eta) \le 4$
\\
&   & 3.3~~Non-vanishing, degree and saturation theorems
\\
&   & 3.4~~Rationality and polynomiality theorems
\\
&   & 3.5~~Parabolic Kostant partition function $K_{\Phi(\eta)}(\gamma)$ as
function of $\gamma$
\\
&   & 3.6~~Reconstruction theorem

\\ 
& 4 & {\bf Parabolic Kostka polynomials: ~~Definition and basic properties}

\\
& 5 & {\bf Parabolic Kostka polynomials: ~~Examples}
\\
&   & 5.1~~Kostka--Foulkes and parabolic Kostka polynomials; ~Merris' conjecture
\\
&   & 5.2~~Parabolic Kostka polynomials and the Littlewood--Richardson numbers
\\
&   & 5.3~~McMahon polytope and rectangular Narayana numbers
\\
&   & 5.4~~Gelfand--Tsetlin  polytope $GT((2^k,1^n),(1^{2k+n}))$
\\
&   & 5.5~~One dimensional sums and parabolic Kostka polynomials
\\
& 6 & {\bf Parabolic Kostka polynomials:  ~~Conjectures}
\\  
&   & 6.1~~Non-vanishing Conjecture
\\
&   & 6.2~~Positivity Conjecture 
\\
&   & 6.3~~Generalized Saturation Conjecture for parabolic Kostka polynomials
\\
&   & 6.4~~Rationality Conjecture
\\
&   & 6.5~~Polynomiality Conjecture
\\
&   & 6.6~~Generalized Fulton Conjecture
\\
&   & 6.7~~$q$-Log concavity and $P$-positivity Conjectures
\\
&   & 6.8~~Generalized Fomin-Fulton-Li-Poon Conjectures
\\
&   & 6.9~~Miscellany
\\
& 7 &  {\bf References}
\end{tabular}

\vskip 0.5cm

\section{Introduction}

This note is based on a series of lectures given by the author during 
1998--2003 years 
concerning the interrelations between the saturation properties of the
Littlewood--Richardson numbers and their several generalizations, parabolic 
$q$-Kostant partition functions and parabolic Kostka polynomials. 

In spite of the title {\it ``An invitation to the Generalized Saturation 
Conjecture''},  we will state a big
amount of conjectures (about 30) and problems (about 15) revolving around a 
very mysterious behavior of the
coefficients, and the leading term especially, of a parabolic Kostka 
polynomial. 

Remember that,
by definition, a function $f:\Omega \subset \Z^{n} \rightarrow \Z$ satisfies 
the saturation property ( on the set $\Omega$ ), if the following 
condition holds:

 $f(N\omega) \ne 0$ for some integer $N \ge 1$ and $\omega  \in \Omega,$ 
then also $f(\omega) \ne 0.$  

For example, any homogeneous function $f$  on the set $\Omega,$  i.e. that 
satisfying the  condition $f(N\omega)=N^{\al}~f(\omega)$ for some 
$\al \in \R,$ $\forall \omega \in \Omega$ and all integers $N \ge 1,$ 
possesses the saturation property; 
a subset $\Omega \subset \Z^{n}$ is
called {\it saturated} if its characteristic function has the saturation 
property.

To be more specific, let us  introduce the numbers 
$a(\ld,\mu \Vert \eta),b(\ld,\mu \Vert \eta),
c(\ld,\mu \Vert \eta)$ and $d(\ld,\mu \Vert \eta)$ which will play an 
important role in 
our paper. Namely, let $\ld$ be a partition  and $\mu,$ and $\eta$ be 
compositions such that $|\ld|=|\mu|$ and $ll(\mu) \le |\eta|,$ see 
Section~2.1 for 
explanation of notation. Let $K_{\ld\mu\eta}(q)$ be
the corresponding parabolic Kostka polynomial.   If 
$K_{\ld\mu\eta}(q) \ne 0,$  the numbers above are 
defined from the decomposition
$$ K_{\ld\mu\eta}(q)= b(\ld,\mu \Vert \eta)~q^{a(\ld,\mu \Vert \eta)}+\cdots+
d(\ld,\mu \Vert \eta)~q^{c(\ld,\mu \Vert \eta)},
$$
where we assume that $b(\ld,\mu \Vert \eta) \ne 0$  and 
$d(\ld,\mu \Vert \eta) \ne 0,$ and $a(\ld,\mu \Vert \eta) \le c(\ld,\mu \Vert \eta).$

If $K_{\ld\mu\eta}(q)=0,$ we put by
definition $a(\ld,\mu \Vert \eta)=b(\ld,\mu \Vert \eta)=c(\ld,\mu \Vert \eta)=
d(\ld,\mu \Vert \eta)=0.$                   \\

$(\maltese)$ We {\bf expect} that $d(\ld,\mu \Vert \eta) \ge 0$  and 
$d(\ld,\mu \Vert \eta) > 0,$  or equivalently $K_{\ld\mu\eta}(q) \ne 0,$
~if and only if ~$\ld-\mu \in Y_{\eta}.$  In other words, we {\bf expect} that 
$K_{\ld\mu\eta}(q) \ne 0,$ or equivalently $d(\ld,\mu \Vert \eta) > 0,$  ~if 
and only if $K_{\Phi(\eta)}(\ld - \mu) > 0,$ 
see Section~6, Positivity and Non-vanishing conjectures. 

$(\clubsuit)$  We  regard  the numbers 
$d(\ld,\mu \Vert \eta)$ as a generalization of the Littlewood--Richardson 
coefficients, see comments after Theorem~1.4, and Section~5.2, $(1^{0})$ 
for explanations.
\begin{prb} Find combinatorial and/or algebro--geometric 
interpretations of the numbers $d(\ld,\mu \Vert \eta).$
\end{prb}
 \begin{rem} 
{ \rm We {\bf expect} that for given $\ld,$ $\mu$ and $\eta$ there
exists a rational convex polytope $\Delta_{\ld,\mu}^{\eta}$ such that the 
number of integer points inside of $\Delta_{\ld,\mu}^{\eta}$ is equal to 
$d(\ld,\mu \Vert \eta).$  }
\end{rem}

One of our main observations is that the saturation 
property of the leading coefficient $d(\ld,\mu \Vert \eta),$ i.e. that

$(\clubsuit)$  ~~$d(N\ld,N\mu \Vert \eta) \ne 0$ for some integer 
$N \ge 1$ if and only if    $d(\ld,\mu \Vert \eta) \ne 0,$  

is an easy consequence ( but not conversely ! )
of the statement that the maximal degree $c(\ld,\mu \Vert \eta)$ of $q$ in a
parabolic Kostka polynomial $K_{\ld\mu\eta}(q)$ is a {\it homogeneous
degree 1} function of $\ld$ and $\mu$.  In other words, we pose the
following conjecture:
\begin{con} ( {\bf Generalized Saturation Conjecture} )

Let $\ld$ be a partition, $\mu$ and $\eta$ be compositions such that $|\ld|=
|\mu|$ and $ll(\mu) \le |\eta|.$ Then the coefficient $c(\ld,\mu \Vert \eta)$ 
is a homogeneous piecewise linear function of $\ld$ and $\mu.$  In particular,

$$ c(N\ld,N\mu \Vert \eta)=N c(\ld,\mu \Vert \eta) $$

for any positive integer $N$.
\end{con}
Here $ll(\mu)$ denotes the {\it fake} length of a composition $\mu,$ see
Section~2.1 for the definition.

We would like to note here that, in general, the Generalized Saturation 
Conjecture ( $GSC$ for short ) is false for the numbers
$a(\ld,\mu \Vert \eta),$ see Examples~4.6. 

$(\maltese)$  However, we {\bf expect} that if 
$\mu$ is a {\it partition}, then the $GSC$  does hold for the
numbers $a(\ld,\mu \Vert \eta).$

Conjecture~1.3  is obvious for the Kostka--Foulkes polynomials, since in this 
case 
$$c(\ld,\mu \Vert (1^{|\ld|}))=n(\mu)-n(\ld)= \sum_{1 \le i < l \le l(\mu)} 
min(\mu_i,\mu_j)-\sum_{1 \le i < j \le l(\ld)}min(\ld_i,\ld_j)$$ 

is easily seen to be a homogeneous piecewise linear function of $\ld$ and 
$\mu.$  However, it seems a difficult {\bf problem} to prove the $GSC$ in 
general case, especially to find an explicit 
piecewise linear formula for the numbers $c(\ld,\mu \Vert \eta).$

We would like to add also that the saturation property of the coefficient 
$d(\ld,\mu \Vert \eta)$ is an easy consequence of Non-vanishing conjecture as
well. Indeed, 
$$d(N\ld,N\mu \Vert \eta) \not= 0  \Leftrightarrow  N(\ld - \mu) \in
Y_{\eta}  \Leftrightarrow \ld - \mu \in Y_{\eta} \Leftrightarrow 
d(\ld,\mu \Vert \eta) \not= 0. $$

Now let us explain briefly a connection between our Generalized Saturation
Conjecture and the Saturation Theorem by A.~Knutson
and T.~Tao \cite{KT1}, see also \cite{Bel}, \cite{Bu1}, \cite{DW1},
\cite{Sp}  for other proofs.
\begin{theorem} ( {\bf Saturation Theorem \cite{KT1}} )

Let $\ld,\mu$ and $\nu$ be partitions such that $|\ld|+|\mu|=|\nu|$. Then

$c_{N\ld,N\mu}^{N\nu} \ne 0$ for some integer $N \ge 1$ if and only if
$c_{\ld,\mu}^{\nu} \ne 0.$ 
\end{theorem}
Here $c_{\ld,\mu}^{\nu}$ denotes the Littlewood--Richardson number 
( $LR$-number for short )
corresponding to the partitions $\ld,\mu$ and $\nu,$ 
see Section~2.4 for details.

Now we are going to explain how the Saturation Theorem follows from the $GSC$.

First of all, we observe that $c_{\ld,\mu}^{\nu}=b( \Lambda,R)$ for some
partition $\Lambda$ and a dominant sequence of rectangular shape 
partitions $R,$ see Section~2 for the definition of a dominant sequence of 
partitions.  Namely, for given partitions $\ld=(\ld_1,\cdots,\ld_r),$ $\mu$ 
and $\nu$ such  that $|\ld|+|\mu|=|\nu|,$ define partition
$$\Lambda=(\mu_1+\ld_1,\mu_1+\ld_2,\cdots,\mu_1+\ld_r,\mu),$$
and a dominant rearrangement $R$ of the sequence of rectangular shape 
partitions $\wt R= \{ (\mu_1^{\ld_1'}),\nu \}.$  Then 

$(\clubsuit)$  ~~$a(\Lambda,R) \ge
\sum_{1 \le j \le \mu_1}\nu_j'-|\mu|,$ and $a(\Lambda,R)=
\sum_{1 \le j \le \mu_1}\nu_j'-|\mu| $ if and only if  
$c_{\ld,\mu}^{\nu} \ge 1;$ in addition, $b( \Lambda,R)=c_{\ld,\mu}^{\nu},$ 
see Section~5.2 for details.

In other words, {\bf the constant term} of the polynomial
$$ K_{\ld,\mu}^{\nu}(q):=q^{(|\mu|-\sum_{1 \le j \le \mu_1}\nu_j')}
~K_{\ld,R}(q)
$$
is {\bf equal} to the Littlewood--Richardson number $c_{\ld,\mu}^{\nu}.$  See
Sections~5.2 and 6.8 where some results and conjectures about the polynomials
$K_{\ld,\mu}^{\nu}(q)$ and their generalizations $K_{A,B,\theta}^{\nu}(q),$
and $K_{A^{(1)},\cdots,A^{(k)},\theta}^{\nu}(q)$, are presented.

The next step is to apply the Duality Theorem for parabolic Kostka polynomials
$K_{\ld,R}(q)$ corresponding to a dominant sequence of {\it rectangular} 
shape partitions $R,$  see Section~4, (4.37), Duality Theorem. As a corollary, 
we see that the coefficients $a(\Lambda,R)$ and $c(\Lambda,R)$ satisfy the 
$GSC$ simultaneously. Hence, it follows from our Theorem~1.5 that

$(\clubsuit)$  ~~$a(n\Lambda,nR)=na(\Lambda,R)$ for any integer $n \ge 1.$

Finally, let us deduce the Saturation Theorem from the above considerations. 
Indeed, assume that $c_{N\ld,N\mu}^{N\nu} \ne 0,$ then

$$ Na(\Lambda,R)=a(N\Lambda,NR)=N(\sum_{1 \le j \le \mu_1} \nu_j'-|\mu|),$$

and therefore, $a(\Lambda,R)=\sum_{1 \le j \le \mu_1} \nu_j'-|\mu|.$  The last
equality means  that $c_{\ld,\mu}^{\nu} \ne 0.$      \\

In fact, our arguments show that in the particular case under consideration, 
the Generalized Saturation Conjecture {\it is equivalent} to the Saturation
Theorem. However, our main point is that, conjecturally, the $GSC$ is still 
valid for  any partition $\ld$ and compositions $\mu,$ and $\eta.$

\begin{theorem} ( {\bf Saturation Theorem for the numbers $c(\ld,R)$ } )

Let $\ld$ be a partition and $R$ be a dominant sequence of {\bf rectangular}
shape partitions. Then

$(\clubsuit)$  ~~$c(N\ld,NR)=Nc(\ld,R)$ for any integer $N \ge 1.$
\end{theorem}
Our proof of Theorem~1.5 is based on an
explicit homogeneous piecewise linear formula for the 
Lascoux--Sch\"utzenberger statistics {\it charge}, obtained 
by A.~Berenstein and A.N.K., see \cite{KB}, \cite{Kir10}, 
and a fermionic formula for the parabolic Kostka polynomials $K_{\ld,R}(q)$
corresponding to a dominant sequence of rectangular shape partitions $R,$ see
e.g. Section~5.1, $(5^{0}).$  The proof is rather technical and long. We
assume to present it in a separate publication.   \\
 
 One of our main results, see Section~4,  in support of the  $GSC$ in general
 case is:

\begin{theorem} ( {\bf Rationality theorem for parabolic Kostka polynomials, I} )
 
The formal power series 
$$ \sum_{n \ge 0} K_{n\ld,n\mu,\eta}(q)t^{n} $$
is a rational function in  $q$ and $t$  of the form
$$ P_{\ld\mu\eta}(q,t)/ Q_{\ld\mu\eta}(q,t),$$
where $P_{\ld\mu\eta}(q,t)$ and $Q_{\ld\mu\eta}(q,t)$ are mutually prime
polynomials in $q$ and $t$ with integer coefficients, $P_{\ld\mu\eta}(0,0)=1.$

Moreover,

$(\clubsuit)$ the denominator $Q_{\ld\mu\eta}$ has the following form
$$ Q_{\ld\mu\eta}(q,t) = \prod_{j \in J} (1-q^{j}~t),$$
where $J:=J_{\ld\mu\eta}$
is a finite set of non--negative integer numbers, not necessarily distinct;

$(\clubsuit\clubsuit)$  $P_{\ld\mu\eta}(1,t)=
(1-t)^{t(\ld,\mu,\eta)}~P_{\ld\mu\eta}(t),$
where $t(\ld,\mu,\eta) \in \Z_{\ge 0},$  and
$P_{\ld\mu\eta}(t)$ is a polynomial with {\bf non--negative} integer 
coefficients. 
\end{theorem}
\begin{prb} Find combinatorial and algebro-geometric interpretations of the
set $J_{\ld\mu\eta}$ and the polynomial $J_{\ld\mu\eta}(q)=
\sum_{j \in J_{\ld\mu\eta}} q^{j}.$
\end{prb}
\begin{cor} ( {\bf Polynomiality theorem for parabolic Kostka numbers} )

Let $\ld$ be a partition and $\mu,\eta$ be compositions such that $|\ld|=|\mu|$
and $ll(\mu) \le |\eta|$. Then 
there exists a polynomial ${\cal K}_{\ld\mu\eta}(t)$ with rational
 coefficients such that for all integers $N \ge 1$
$$ {\cal K}_{\ld\mu\eta}(N)=K_{N\ld,N\mu,\eta}(1).
$$
\end{cor}
\begin{cor}( {\bf Polynomiality theorems for Kostka and $LR$-numbers} ) 

$(i)$ Let $\ld$ be a partition and $\mu$ be a composition of the same size, 
then the  Kostka number $K_{N\ld,N\mu}(1)$ is a polynomial in $N$ with 
rational  coefficients.

$(ii)$ Let $\ld,\mu$ and $\nu$ be partitions, then the Littlewood--Richardson
number $c_{N\ld,N\mu}^{N\nu}$ is a polynomial in $N$ with rational 
coefficients.
\end{cor}
See Section~4, Theorem~4.14 and Corollary~4.15. We also give a multivariable
generalization of Theorem~1.6, see Theorem~4.17.

We want to {\it emphasize} here that the polynomiality property of the 
functions $N \rightarrow K_{N\ld,N\mu,\eta}(1)$ and $N \rightarrow 
c_{N\ld,\N\mu}^{N\nu}$  is an easy consequence of our Theorem~1.6, 
but not conversely: one has to  
check  that the (irredundant) denominator $Q_{\ld\mu\eta}(q,t)$ doesn't have 
factors of the form $(1-q^k t^l)$ with $l \in \Z_{ >1}.$ 
\begin{con} If $\mu$ is a partition, then the polynomial 
${\cal K}_{\ld\mu\eta}(t)$ has {\bf non--negative} rational coefficients.
\end{con}
See Section~6, Conjecture~6.10, $(\blacklozenge)$, for more general conjectures
concerning the numbers $d(\ld,\mu \Vert \eta).$

We would like to remark that the $GSC$ does not follow immediately from 
Theorem~1.6, see Section~6, Rationality Conjecture, for details.

The polynomials $P_{\ld\mu\eta}(q,t)$ may have negative coefficients, and 
rather difficult to compute. For example, we don't know the explicit formula 
for polynomial $P_{(2^6),(1^{12}),(1^{12})}(q,t).$  We expect that the 
polynomials $P_{\ld\mu\eta}(q,t)$ should have nice algebraic and 
algebro--geometric interpretations.

Our proof of Theorem~1.6 is a pure algebraic and is based on the study of
 the parabolic $q$-Kostant partition functions, see Section~3.

Corollary 1.9,$(i)$, has been proved independently by W.~Baldoni--Silva and
M.Vergne \cite{BV}, S.~Billey, V.~Guillemin and E.~Rassart \cite{BGR}, ... .
Corollary 1.9,$(ii)$, has been proved independently by A~Knutson (unpublished),
 H.~Derksen and J.~Weyman \cite{DW2}, E.~Rassart \cite{Ras}, ... . \\

The main subject of investigation of our paper is the  study of interrelations between the saturation properties of the $LR$-numbers and their 
generalizations, and the leading coefficient of the parabolic Kostka 
polynomials.

The paper does not contain complete proofs  of the main theorems. Our
goal is different. The primary  purpose of this note  is to collect together
several results, conjectures and examples revolving around a mysterious 
behavior of the initial and the leading terms of a parabolic Kostka 
polynomial.                   \\

Let us say a few words about the content of our paper. 

In Section~2 we collect
together a few definitions and notation which will be frequently used in the 
subsequent Sections. 

In Section~2.1 we remember the definitions of partitions and
compositions and some operations over them. We would like to point out here
 some non--standard conventions about
partitions and compositions used in our paper. We will denote by 
$\ld=(\ld_1,\cdots,\ld_r)$ a (proper) partition, so that if 
$\ld \ne \emptyset,$
then $\ld_r \ne 0$. We always use $\eta$ to denote a composition without zero
components. Contrary, we will use $\mu$ to denote a composition or partition 
with zero components and zeroes at the end  allowed. A typical example is
$\mu=(0,2,0,1,3,0,0).$ Thus, according to our conventions, the compositions
$(0),(0,0),\cdots $ are different and different from the empty composition 
$\emptyset.$

In Sections~2.2 and 2.3 we recall the definitions of Kostka--Foulkes and skew
Kostka--Foulkes polynomials. For more details, see \cite{But}, \cite{DLT},
\cite{Kir1}, \cite{Kir9}, \cite{KS}, \cite{KSS}.

In Section~2.4 we remember the definition of  the Littlewood--Richardson
numbers and state the {\it Saturation Theorem }, which has been proved by
A.~Knutson and T.~Tao \cite{KT1}.

We refer the reader to interesting and clearly written papers by W.~Fulton
\cite{Fu2}, \cite{Fu3} for detailed account to the so--called {\it Horn
problem} and its connections with the Saturation Theorem. 

In Section~2.5 we study the saturation properties  of 
the {\it internal product}
structural constants $g_{\al\beta\gamma}$ and those of the {\it plethysm}
  $a_{\mu,W}^{\pi}.$  It is well--known that the $LR$-numbers 
$c_{\ld,\mu}^{\nu}$ are a special case of the internal product structural
constants $g_{\al\beta\gamma},$ and in turn, the numbers 
$g_{\al\beta\gamma}$ are a special case of the plethysm structural constants 
$a_{\mu,W}^{\pi},$ see Remark~2.13. However, based on examples we arrived at 
the conclusion that, in the general case, both the numbers 
$g_{\al\beta\gamma}$ and $a_{\mu,W}^{\pi}$ do {\bf not} satisfy the 
{\it saturation property}. 

$(\maltese)$  Nevertheless, we {\bf expect} that 

$\bullet$ the numbers $a_{\mu,W}^{\pi}$ satisfy a {\it weak form} 
of Saturation  Conjecture, i.e. for any finite dimensional 
${\g}l_n$-module $W$ there exists a polynomial $p_{W}(t)$ ( $p_{W}(t)=t$ ?? )
 such that for all partitions $\pi$ and $\mu$ one has  

 if $a_{N\mu,W}^{N\pi} \ge p_{W}(N),$  then $a_{\mu,W}^{\pi} \ne 0.$  

$\bullet$  for an interesting family of polynomials  $L_{\al,\beta}^{\mu}(q)$ 
a certain analog of the $GSC$ does hold, see  Conjecture~2.22.

It seems an interesting {\bf problem} to study whether or not the $GSC$  
is valid for polynomials $M_{\mu,W}^{\pi}(q)$ which are defined via the 
decomposition of the plethysm $W \circ s_{\mu}$: 

$$ (W \circ s_{\mu}) (X)=\sum_{\pi} M_{\mu,W}^{\pi}(q)~P_{\pi}(X,q),$$
where $X=(x_1,\cdots,x_n),$ and $P_{\pi}(X,q)$ stands for the Hall--Littlewood
polynomials.

In Section~2.5 we also state  several results about polynomials 
$L_{\al,\beta}^{\mu}(q)$ and give a few examples supporting our conjectures.

In Section~2.6 we define the extended Littlewood--Richardson numbers as well
as the level $l$ extended $LR$-numbers. The latter are a  natural 
generalization of the restricted $LR$-numbers. 

$(\maltese)$ We {\bf expect} that 
Saturation Theorem, the strong $q$-log concavity and Fomin-Fulton-Li-Poon's 
conjectures I and II are still valid for the level $l$ extended $LR$-numbers.
 \\

In Section~3 we study some algebraic properties of the parabolic $q$-Kostant
partition function $K_{\Phi(\eta)}(\gamma|~q),$ mainly in a connection with
the saturation properties of the latter. For polynomials 
$K_{\Phi(\eta)}(\gamma|~q)$ we prove an analog of the $GSC,$ Rationality and
Polynomiality theorems, and a new recurrence relation. Our proof of Rationality
theorem is based on the following simple observation:
\begin{lem} Let $R(X,q) \in \Q~[q] \lbrack\!\lbrack X^{\pm 1} \rbrack\!\rbrack$ be a  rational function in q and  
$X=(x_{1}^{\pm 1},\cdots,x_{n}^{\pm 1}).$ Let
$$ R(X,q)=\sum_{m \in \Z^{n}} A_{m}(q) X^{m} $$
be a Laurent series  expansion of $R(X,q).$

Let $a_1,\cdots,a_k \in \Z^{n}$ be fixed, then
$$ \sum_{(N_1,\cdots, N_k) \in \Z_{\ge 0}^{k} } A_{N_{1}a_{1}+\cdots+
N_{k}a_{k}}(q) ~x_1^{N_1} \cdots x_k^{N_k}$$
is a rational function in  $q$ and $x_1,\cdots,x_k.$  
\end{lem}
In Section~3 we also study the parabolic Kostant partition function 
$K_{\Phi(\eta)}(\gamma)$ as a {\it function} of $\gamma,$ see Theorems~3.23
and 3.25.

A detailed treatment of the properties of the parabolic $q$-Kostant and 
Kostant partition functions  lies at the heart of the approach to the $GSC$ 
and to the study of parabolic Kostka polynomials, presented in this paper. 
However, making an effort to keep the paper in a reasonable size, 
we do not intend to consider in Section~3, and decided to postpone for 
subsequent publications, many very interesting aspects of the theory of
parabolic Kostant partition function $K_{\Phi(\eta)}(\gamma):=
K_{\Phi(\eta)}(\gamma|~q) \vert_{q=1}$
such as 

$(i)$ The special values of parabolic Kostant partition function, 
see \cite{PS}, \cite{Kir8}, \cite{Kir9}, \cite{BV};

$(ii)$ Connections with the flow polytopes, see \cite{PS}, \cite{BV};

$(iii)$ Connections with the Orlik--Solomon and Gelfand--Varchenko algebras,
 \cite{Kir15};

$(iv)$ A $q$-analog of the generalized Kostant partition function, 
see \cite{PS}.     \\
 
In Section~4 we study, mainly, the ``saturation properties''  of parabolic 
Kostka polynomials. Many examples, results and conjectures concerning
with the parabolic Kostka polynomials, have been already considered in our
 paper \cite{Kir9}. For the reader's convenience, in the present paper we 
remember some basic properties of the 
parabolic Kostka polynomials $K_{\ld\mu\eta}(q),$ and give a sketch of proofs
of Rationality and Polynomiality theorems for the latter, see
Theorems~4.14 and 4.17, and Corollaries~4.15, 4.18 and 4.19.

In the case when $\mu$ and $\eta$ correspond to a dominant sequence of 
rectangular shape partitions $R,$ we have obtained the following result:
\begin{theorem} ( {\bf Polynomiality theorem for the numbers $b(\ld,R)$} )

Let $\ld$ be a partition and $R$ be a dominant  sequence of 
rectangular shape partitions, then

$(\clubsuit)$ ~~$b(N\ld,NR)$ is a polynomial in $N$ with rational coefficients.
\end{theorem}   
Our proof of Theorem~1.6 is a largely  algebraic, whereas that of Theorem~1.12
is based on a fermionic formula for the parabolic Kostka polynomials 
$K_{\ld,R}(q).$

$(\maltese)$ We {\bf expect} that if $\mu$ is a {\it partition}, 
then $b(N\ld,N\mu \Vert \eta)$
is a polynomial in $N$ with {\it non--negative} rational coefficients, see
Section~6, Polynomiality conjecture, for a more detailed statement. 

However, in general, $b(N\ld,N\mu \Vert \eta)$ becomes a polynomial
in $N$ only starting from big enough $N,$ see Section~6, Conjecture~6.10, 
$(\blacklozenge\blacklozenge\blacklozenge),$  and Remark~6.16.  

In Section~4 we also study some natural multivariable analogues of Theorem~1.6,
and Corollaries~1.7  and 1.8. In particular, we give a sketch of proof of a 
theorem that
for any sequences of partitions $\ld^{(1)},\cdots,\ld^{(k)}$ and compositions
$\mu^{(1)},\cdots,\mu^{(k)}$ the formal power series
$$\sum_{(N_1,\cdots,N_k) \in \Z_{\ge 0}^{k}}K_{N_{1}\ld^{(1)}
+\cdots+N_{k}\ld^{(k)},N_{1}\mu^{(1)}+\cdots+N_{k}\mu^{(k)},\eta}
(q)~x_1^{N_1} \cdots x_k^{N_k}
$$
is a rational function in $q$ and $x_1,\cdots,x_k,$  which has  the denominator
of some special form, see Section~4, Theorem~4.17.

However, in general, if $k \ge 2,$ the functions 
$$(N_1,\cdots,N_k) \rightarrow
K_{N_{1}\ld^{(1)}+\cdots+N_{k}\ld^{(k)},
N_{1}\mu^{(1)}+\cdots+N_{k}\mu^{(k)},\eta}(1),~~~and$$ 
$$(N_1,\cdots,N_k) \rightarrow c_{N_{1}\ld^{(1)}+\cdots+N_{k}\ld^{(k)},
N_{1}\mu^{(1)}+\cdots+N_{k}\mu^{(k)}}^{N_{1}\nu^{(1)}+\cdots+N_{k}\nu^{(k)}}$$
 are only piecewise polynomial functions on  the set 
$\{(N_1,\cdots,N_k) \in \Z_{\ge 0}^{k} \},$ see Example~4.23.

We want to {\it emphasize} here that the special form of the denominator of the rational function $\sum_{(N_1,\cdots,N_k) \in \Z_{\ge 0}^{k}}K_{N_{1}\ld^{(1)}+\cdots+N_{k}\ld^{(k)},N_{1}\mu^{(1)}+\cdots+N_{k}\mu^{(k)},\eta}(1),$ see 
Theorem~4.17,$(\clubsuit),$ is (in our opinion) a key fact to explain a 
piecewise polynomiality of the ``mixed'' Kostka numbers $K_{N_{1}\ld^{(1)}+
\cdots+N_{k}\ld^{(k)},N_{1}\mu^{(1)}+\cdots+N_{k}\mu^{(k)},\eta}(1)$ and
``mixed'' Littlewood--Richardson coefficients. 

$(\maltese)$ Nevertheless, we {\bf expect} that in the case of parabolic 
Kostant's partition functions, the function $(n_1,\cdots,n_k) 
\rightarrow K_{\Phi(\eta)}(n_1\gamma_1+\cdots+n_k\gamma_k)$ is a polynomial
one on the whole set $\{ (n_1,\cdots,n_k) \in \Z_{\ge 0}^k \}.$

It is well-known that the Kostka--Foulkes number $K_{\ld\mu}(1)$ counts the
number of integral points in some rational convex polytope, the so-called 
Gelfand--Tsetlin polytope $GT(\ld,\mu).$  In this connection we would like 
to pose the following question (cf with {\it mixed lattice point enumerator}
theorem for integer convex polytopes by P.McMullen \cite{McM}, 
or Example~4.23) :
\begin{quest} Let $\Delta_1, \cdots, \Delta_k \in \Q^{d}$ be rational convex 
polytopes, and $L: \Z^{d} \rightarrow \Z_{ \ge 0}$ be a continuous piecewise 
linear function.

Under what assumptions on $L$ and polytopes $\Delta_1,\cdots,\Delta_k$ the
{\bf denominator} of rational function
$$\sum_{(N_1,\cdots,N_k) \in \Z_{\ge 0}^k} \{ \sum_{a \in (N_1\Delta_1+\cdots+
N_k\Delta_k) \cap \Z^{d}}q^{L(a)} \} x_1^{N_1}\cdots x_k^{N_k}$$
has only the factors of the form $(1-q^{a_{J}^{(j)}}x_{J}),$ where 
$J \subset [1,\cdots,k],$ $x_J:=\prod_{j \in J}x_{j},$ and $a_{J}^{(j)}$ are
some non--negative integers ?
\end{quest}

In Section~4, Remark~4.24, we state some preliminary results about the 
behavior
of the parabolic Kostka number $K_{\ld\mu\eta}(1)$ considered as a function
of $\ld$ and $\mu$ on ``the space of parameters'' $Z_{\eta}=\{ (\ld,\mu) \in
 \Z_{\ge 0}^{n} \times  \Z_{\ge 0}^{n} \mid \ld_1 \ge \cdots \ge \ld_n, 
\ld-\mu \in Y_{\eta} \}.$  Based on the properties of the parabolic Kostant
partition function, see Section~3, Theorem~3.25, one can show that on the
set  $Z_{\eta}$ the parabolic Kostka number $K_{\ld\mu\eta}(1)$ is a 
continuous piecewise  polynomial function in $\ld_1,\cdots,\ld_n$ and 
$\mu_1,\cdots,\mu_n.$  The main problem about the function 
$(\ld,\mu) \rightarrow  K_{\ld\mu\eta}(1)$ we are interested in, 
is to describe ``the dominant  chamber'' for the latter function, i.e. 
to describe the maximal domain 
$Z_{\eta}^{++}$ in the set  $Z_{\eta}^{+}:=\{ (\ld,\mu) \in Z_{\eta} \mid
\ld -\mu \in Y_{\eta}^{+} \}$  such that $K_{\ld\mu\eta}(1)
 \vert_{Z_{\eta}^{++}}=K_{\Phi(\eta)}(\ld-\mu).$

In Section~4 we also introduce the parabolic Hall--Littlewood polynomials
$Q_{\mu,\eta}(X;q),$ and state the rationality theorem for the latter, see
 Remark~4.35. Details and proofs  will appear in a separate publication.
Finally, we note that for the Kostka--Macdonald polynomials $K_{\ld,\mu}(q,t),$
see \cite{Ma}, Chapter~VI, Section~8, for the definition, the generating
function
$$Z_{\ld,\mu}(q,t,x):= \sum_{n \ge 0}~K_{n\ld,n\mu}(q,t)~x^n $$
is a formal power series, which  {\it is not}, in general, a { \it rational} 
function in  $q$, $t$ and $x.$

It seems a very interesting {\bf problem} to study the properties of the 
function  $Z_{\ld,\mu}(q,t,x),$ especially in connections with the characters
of affine Lie algebras of type $A$ and the Virasoro algebra.    \\
 
In Section~5 we collect together several examples which might help to 
illuminate a mysterious nature  of the leading term of a parabolic Kostka 
polynomial. See the Contents of Section~5 for exposing  with 
the list of these examples. In particular, we show that the one dimensional 
sums (1D-sums for short) which frequently appear in Statistical Mechanics,
 see e.g. \cite{HKKOTY}, \cite{KMOTU2} and the literature quoted therein, are
 a special case of the parabolic Kostka polynomials $K_{\ld\mu\eta}(q)$ 
corresponding to a  rectangular shape partition $\ld,$ see Section~5.5 for 
details. In Section~5.1 we give, among other things, a few comments about the
Merris conjecture, and in Section~5.4 that about the $LR$-numbers 
$c_{\delta_{n},\delta_{n}}^{\ld}.$  \\

In Section~6 we state a few conjectures about the coefficients 
$a(\ld,\mu \Vert \eta),
b(\ld,\mu \Vert \eta),c(\ld,\mu \Vert \eta)$ and $d(\ld,\mu \Vert \eta).$ 
In particular, we expect, see Conjectures~6.14, 6.17 and 6.23, that 

$\bullet$  ( {\bf The generalized Fulton  conjecture} )

If ~$d(n\ld,n\mu \Vert \eta)=1$ for some integer $n \ge 1,$ then
~$d(N\ld,N\mu \Vert \eta)=1$ for all $N \in \Z_{\ge 1}.$

$\bullet$ ( {\bf Generalized $d(\ld,\mu \Vert \eta)=2$ conjecture} )

If $d(n\ld,n\mu \Vert \eta)=n+1$ for some integer $n \ge 1,$ then
$d(N\ld,N\mu \Vert \eta)=N+1$ for all $N \in \Z_{\ge 1}.$

$\bullet$ ( {\bf Generalized $d(\ld,\mu \Vert \eta)=3$ conjecture} )

$(i)$ ~If ~$d(n\ld,n\mu \Vert \eta)=2n+1$ for some integer $n \ge 2,$ then
~$d(N\ld,N\mu \Vert \eta)=2N+1$ for all $N \in \Z_{\ge 1};$

$(ii)$ ~If ~$d(n\ld,n\mu \Vert \eta)=
\left(\begin{array}{c}n+2\\2\end{array}\right)$  for some integer $n \ge 2,$
 then
~$d(N\ld,N\mu \Vert \eta)=\left(\begin{array}{c}N+2\\2\end{array}\right)$ 
for all $N \in \Z_{\ge 1}.$

These two cases exhaust  the all possibilities when $d(\ld,\mu \Vert \eta)
=3.$

$\bullet$  ( {\bf $q$-Log concavity conjecture} )
 
Let $\ld$ be a partition and $R$ be a dominant sequence of 
rectangular shape partitions, then for any integer $n \ge 1,$
$$ (K_{n\ld,nR}(q))^2 \ge K_{(n-1)\ld,(n-1)R}(q)~K_{(n+1)\ld,(n+1)R}(q).
$$
See Section~6.7, Conjecture~6.17, for a more general and detailed 
statement of the latter conjecture.

$\bullet$ ( {\bf The generalized Fomin-Fulton-Li-Poon's conjecture I, 
cf \cite{Ok1}, Conjecture~1, \cite{FFLP}, Conjecture~2.7} )
$$ K_{{\wt A}^{(1)},\cdots,{\wt A}^{(k)},\theta}^{\nu}(q) \ge 
K_{A^{(1)},\cdots,A^{(k)},\theta}^{\nu}(q). $$
$\bullet$ ( {\bf The generalized Fomin-Fulton-Li-Poon's conjecture II, cf 
\cite{FFLP}, Conjecture~5.1} )
\footnote{ As we learned from the referee,
the  extension of the original Fomin-Fulton-Li-Poon conjecture II,
 \cite{FFLP},Conjecture~5.1, to the case of skew diagrams  was also  stated
by F.~Bergeron, R.~Biagnoli and M.~Rosas, see e.g. \cite{B}, \cite{BBR};
see also \cite{MacN}. The paper \cite{BBR} contains, among other things,
 many interesting results in support of the FFLP-conjecture.}
$$K_{A^{*},B^{*},\theta}^{\nu}(q) \ge K_{A,B,\theta}^{\nu}(q).$$
See Section~6.8, Conjecture~6.24, for the  explanation of notation
we have used, further details and more conjectures.
 
In the case of the $LR$-numbers the Fulton conjecture has been proved in
\cite{KTW}. 

 Some special cases of the Fomin-Fulton-Li-Poon conjecture II
have been proved in \cite{FFLP}. 
\begin{prb} When does the number $d(\ld,\mu \Vert \eta)$ equal to $1$  ? 
\end{prb}
Finally, we would like to remark that our approach to the $GSC$ is purely
algebraic and combinatorial. It seems a very interesting {\bf problem} to find
an algebro--geometric explanation of a still experimental observation that
the coefficient $c(\ld,\mu \Vert \eta)$ is a homogeneous piecewise linear 
function 
of $\ld$ and $\mu.$ In this connection we would like to pose the following
questions:
\begin{quest} ( {\bf Parabolic Kostka polynomials and semi--invariants of
quivers} )

Let $\ld$ be a partition and $\mu,$ and $\eta$ be compositions
such that $|\ld|=|\mu|$ and $ll(\mu) \le |\eta|.$

Does there exist a quiver $Q$, dimensional vector $\beta$ and 
$GL(Q,\beta)$-weight $\sigma$ such that
$$ \dim SI(Q,\beta)_{n\sigma}=d(n\ld,n\mu \Vert \eta)$$
for all integers $n \ge 1$ ?
\end{quest}

Here $SI(Q,\beta)_{\sigma}$ stands for the weight $\sigma$ subspace of the
ring of semi--invariants 

$$SI(Q,\beta):=\Q~[Rep(Q,\beta)]^{SL(Q,\beta)}.$$

See \cite{DW1} and \cite{DW2}, and the literature quoted therein, for more 
details about the ring of semi--invariants of a quiver. It seems a very 
interesting {\bf problem} to find an interpretation of the numbers 
$c(\ld,\mu \Vert \eta)$ and $d(\ld,\mu \Vert \eta)$ in terms of quivers.

\begin{quest} (  {\bf A $q$-analog of $dim SI(Q,\beta)$} )

Does there exist a natural filtration 
$$ \{ 0={\cal F}_{0} \subset {\cal F}_{1} \subset \cdots  \}$$ 
 on the ring of semi--invariants $SI(Q,\beta)$ such that for a special
quiver $Q=T_{n,n,n}$ and a special dimensional vector $\beta,$  see \cite{DW1},
Section~3,
$$ \sum_{j \ge 1} \dim({\cal F}_{j}/{\cal F}_{j-1})~q^{j} \bdoteq 
c_{\ld,\mu}^{\nu}(q)~ ? $$  
Here $c_{\ld,\mu}^{\nu}(q)$ denotes the $q$-analog of the $LR$-numbers,
see e.g. \cite{CB}, \cite{LT}; for the meaning of the symbol ``$\bdoteq$'',
see Section~1.1.
\end{quest}
We would like to end this Introduction by the following remark. 
Throughout the paper we use the term {\bf Conjecture} to mean a statement 
for which we do not have a proof, but which we have checked on a big body of
 examples (except for Conjectures  from Section~6.9). On the
other hand, we use an expression `` We {\bf expect} that ... `` to
mean a statement which we believe is bound to be true, but for which we
don't have the extensive supporting evidence. Of course, not all plausible
conjectures and reasonable guesses prove to be true. For example, see 
Remark~4.22.

\subsection{Notation}

Throughout the paper we follow Macdonald's book \cite{Ma} as for
notation related to the theory of symmetric functions, and Stanley's
book \cite{St3} as for notation related to Combinatorics. Below we give
a list of some special notation which we will frequently use.

1) If $P(q)$ and $Q(q)$ are polynomials in $q$, the symbol $P(q)\bdoteq
Q(q)$ means that the ratio $P(q)/Q(q)$ is a power of $q$.

2) If $a,k_0,\ldots ,k_m$ are (non--negative) integers, the symbol
$q^a(k_0,\ldots ,k_m)$ stands for the polynomial $\sum_{j=0}^mk_jq^{a+j}$.

3) We use the capital Latin letters $A,B,C, \cdots$ to denote the skew
diagrams/shapes, and the small or capital Greek letters $\al,\beta,\gamma,\ld,
\mu,\Lambda,M, \cdots $ to denote either partitions or compositions.

4) Let $\eta_1=(\eta_{1,1},\eta_{1,2},\cdots,\eta_{1,p})$ and $\eta_2$ be 
compositions, we say that $\eta_2$ is a {\it subdivision} of 
$\eta_1,$ if there exists a sequence of partitions $\mu^{(j)},$  
$1 \le j \le p,$ such that $|\mu^{(j)}|=\eta_{1,j}$ and $\eta_2=
(\mu^{(1)},\cdots,\mu^{(p)}).$

5) Let $P_1(q)$ and $P_2(q)$ be polynomials with real coefficients. By 
definition, the inequality $P_1(q) \ge P_2(q)$ means that the difference
$P_1(q)-P_2(q)$ is a polynomial with non--negative real coefficients.

{ \bf Acknowledgments}

This note presents an extended version of my talks given at

$\bullet$ The International Institute for Advanced Study, Japan, January 29-30
1999; The International Workshop ``Physical Combinatorics''.

$\bullet$ The Nagoya University, Japan, August 23--27; The International 
Workshop  ``Physics and Combinatorics''.

$\bullet$ The Erwin Schroedinger International Institute of Mathematical
Physics, Vienna, Austria, July 2--15, 2000; The International Workshop 
``Representation Theory 2000''.

$\bullet$ The Isaac Newton Institute for Mathematical Sciences, Cambridge, UK,
June 28--July 4, 2001; The International Workshop ``Symmetric Functions and 
Macdonald Polynomials''.

$\bullet$ The Research Institute for Mathematical Sciences, Kyoto, Japan,
August 12--16, 2001; The International Workshop ``Integrable Models, 
Combinatorics and Representation Theory''. 

$\bullet$ The Research Institute for Mathematical Sciences, Kyoto, Japan,
May 22, 2002; Colloquium talk.

$\bullet$ The Kyushu University, Japan, November 17--21, 2003; The Second 
East Asian Conference on Algebra and Combinatorics.

I am grateful to the organizes of these Workshops and Conference for a kind 
invitation, financial support and the opportunity to talk about one of my the
 most  favorite subjects.
 
Finally, I would like to thank the referee for many valuable remarks and
comments, in particular, for sending the references  on
the preprints \cite{B}, \cite{BBR} and \cite{MacN}, which have some
overlaps with the content of our Section~6.8.

\section{Basic definitions and notation}
\subsection{Compositions and partitions}

A composition 
\begin{equation}\mu =(\mu_1,\mu_2,\cdots ,\mu_r) \label{2,1}
\end{equation}   
is a sequence of
non-negative integers. The
number $r$ in (2.1) is called the {\it fake} length of the composition $\mu$, 
and denoted by $ll(\mu)$. In the sequel, it will be convenient for us to
distinguish between two such sequences which differ only by a string of 
zeros at the end. 
Thus, for example, we regard $(2,0,1),(2,0,1,0),(2,0,1,0,0),\cdots,$   as
{\it different} compositions. The size of a composition $\mu$ is defined
to be $|\mu|= \mu_1+\cdots +\mu_r.$  

 By definition, a composition $\lambda
=(\lambda_1,\lambda_2,\cdots ,\lambda_p)$ is called  {\it partition}, if
 additionally  it satisfies the following condition:
\begin{equation} \lambda_1\ge\lambda_2\ge\cdots\ge\lambda_p\ge\ 0. \label{2.2}
\end{equation}
 
The non-zero $\lambda_i$ in (2.2) are called the {\it parts} of $\lambda.$ The
number of parts is the {\it length} of $\lambda$, denoted by $l(\lambda).$
Thus, we have $l(\ld)\le\ ll(\ld):=p.$ As in the case of compositions, we 
distinguish between two sequences (2.2) if they differ only by a string of 
zeros at the end. If $|\lambda| = n$ we say that $\lambda$ is a partition of
$n.$ Denote by ${\cal P}_n$ the set of all partitions of $n.$  \\
 A partition $\ld=(\lambda_1,\lambda_2\,\cdots,\lambda_p)$ is called 
{\it proper} if $\ld_p \ne 0$.

The {\it dominance partial ordering} "$\ge$" on the set of compositions of the
same size $n$, or that of partitions
${\cal P}_n,$ is defined as follows:
\begin{eqnarray*}
&&\lambda\ge\mu ~~{\rm if ~and ~only ~if}\\ &&\lambda_1+\cdots
 +\lambda_i\ge\mu_1+\cdots +\mu_i ~~{\rm for
~all}~~ i\ge 1.
\end{eqnarray*}

The {\it conjugate} of a partition $\ld=(\ld_1,\cdots,\ld_p)$ is the
partition $\ld'=(\ld_1',\ld_2',\cdots),$ where  \\
$\ld_i'= \# \{j | \ld_j\ge i \}.$ In particular, $\ld_1'=l(\ld)$ and
 $\ld_1=l(\ld').$

For each partition $\lambda 
=(\lambda_1,\lambda_2,\cdots ,\lambda_p)$ we define
$$ n(\lambda )=\sum_{i=1}^p(i-1)\lambda_i=\sum_{1\le i<j\le p}\min
(\lambda_i,\lambda_j).
$$

The {\it concatenation} $\mu * \nu$ of two compositions
$\mu = (\mu_1,\mu_2\, \cdots, \mu_r)$ and 
$\nu=(\nu_1,\nu_2, \cdots,\nu_s)$ is defined to be the composition
\begin{equation}
 \mu * \nu=(\mu_1,\mu_2,\cdots ,\mu_r,\nu_1,\nu_2, \cdots,\nu_s). \label{2.3}
\end{equation}
For any compositions $\mu$ and $\nu$ we define 
$\mu + \nu$ to be the sum of the
sequences $\mu$ and $\nu$ : 
\begin{equation}
(\mu+\nu)_i=\mu_i+\nu_i .
\end{equation}
Thus, for example, $n \mu = (n\mu_1,n\mu_2,\cdots ,n\mu_r).$
\begin{de}
We say that a sequence of partitions $\bmu = (\mu^{(1)},\mu^{(2)},\cdots,
\mu^{(r)})$ is a dominant one, if the concatenation 
$\mu^{(1)}*\mu^{(2)}*\cdots*\mu^{(r)}$ is a {\it partition}.
\end{de}
\begin{de}
Let  $\mu=(\mu_1,\mu_2,\cdots,\mu_r)$ and $\eta=(\eta_1,\eta_2,\cdots,\eta_p)$
be compositions, we say that the composition $\mu$
is compatible with $\eta$ if the all compositions 
\begin{equation}
\mu^{(i)}=(\mu_{\eta_1+\cdots +\eta_{i-1}+1},\cdots ,\mu_{\eta_1+\cdots
+\eta_i}), ~~~1\le i\le p  \label{2.5}
\end{equation}
appear to be  partitions (possibly with zeros at the end), where by
definition we put $\eta_{0}:=0.$
\end{de}
 In other words, the composition $\mu$ is
 the  concatenation of partitions $\mu^{(i)}, ~~1\le i\le p.$  Conversely, if a
composition $\mu$ is the concatenation of partitions $\mu^{(i)},~~1\le i\le p,$
then the composition $\eta$ can be reconstructed from that $\mu$ as follows:

$$\eta=(ll(\mu^{(1)}),ll(\mu^{(2)}),\cdots,ll(\mu^{(p)})).$$

\subsection{Kostka--Foulkes polynomials }

In Sections 2.2 till that 2.6 we will assume that all partitions which will
appear, are {\it proper}.
\begin{de}
The Kostka--Foulkes polynomials are defined as the matrix elements of
the transition matrix
$$K(q)=M(s,P)
$$ from the Schur functions $s_{\lambda}(x)$ to the Hall--Littlewood
functions $P_{\mu}(x;q)$:
\begin{equation}
s_{\lambda}(x)=\sum_{\mu}K_{\lambda\mu}(q)P_{\mu}(x;q). \label{2.6}
\end{equation}
\end{de}
It is well known, see e.g.  \cite{Ma}, Chapter~I, that if $\lambda$ and $\mu$
are  partitions, then

$\bullet$ $K_{\lambda\mu}(q)\ne 0$ if and only if $\lambda\ge\mu$ with
respect to the dominance partial ordering "$\ge$" on the set of partitions.

$\bullet$ If $\lambda\ge\mu$, $K_{\lambda\mu}(q)$ is a monic of degree
$n(\mu)-n(\lambda)$ polynomial with {\it non--negative} integer 
coefficients. This result is due to A.~Lascoux and 
M.-P.~Sch\"utzenberger \cite{LS}.

$\bullet$ If $l(\mu)=n,$  then  
\begin{equation}
K_{\lambda\mu}(q):=\sum_{w\in \Sigma_n}(-1)^{l(w)}
K_{n}(w(\lambda+\delta)-\mu-\delta |~q),
\end{equation}
where $l(w)$ denotes the {\it length} of a permutation $w \in \Sigma_n,$ 
 $\delta:=\delta_n=(n-1,n-2,\cdots,1,0)$, and for any $\gamma \in \Z^n,
|\gamma|=0,$ $K_n(\gamma |~q)$ stands for a  $q$-analog of the Kostant 
partition function $K_n(\gamma)$, see e.g. \cite{Ma},~Chapter~III, Section~6,
Example~4, or Section~3 of the present paper.

\begin{theorem} Let $\ld$ and $\mu$ be partitions of the same size. 
There exists a polynomial ${\cal E}_{\ld,\mu }(t)$  with rational coefficients
such that for any integer $N \ge 1$ one has
$$  {\cal E}_{\ld,\mu }(N) = K_{N\ld,N\mu}(1).
$$
\end{theorem}
\begin{cor} The Ehrhart polynomial ${\cal E}_{\ld,\mu }(t)$ of 
the Gelfand--Tsetlin polytope 
$GT(\ld,\mu)$ is a {\bf polynomial}, even though the polytope 
$GT(\ld,\mu)$ itself does not necessary appear to be an integral one.
\end{cor}
For a definition of the Gelfand-Tsetlin polytope see, e.g. \cite{Kir10},
\cite{BGR} or \cite{LM}. For a definition and basic properties of the Ehrhart
polynomial of a convex integral polytope see, e.g. \cite{St3} or \cite{Hi}.
  
 Theorem~2.4 and Corollary~2.5  are  a particular case of a more general 
result, see Section~4, Corollary~4.15.

We refer the reader to a paper \cite{LM} which contains a rich information 
about vertices of Gelfand--Tsetlin's  polytopes. In particular, one can find 
in \cite{LM} several examples of Gelfand--Tsetlin's polytopes with some 
non-integral vertices. 

\begin{con} Let $\ld$ and $\mu$ be (proper) partitions of the same size, then
the Ehrhart polynomial ${\cal E}_{\ld,\mu }(t)$ has {\bf non--negative}
 rational coefficients.
\end{con}
We remark that Conjecture~2.6 is a special case of Polynomiality 
Conjecture from Section~6. 

Polynomiality of the function $N \longrightarrow K_{N\ld,N\mu}(1)$ has
been proved independently by several authors:  W.~Baldoni-Silva and 
M.~Vergne \cite{BV}, S.~Billey, V.~Guillemin and E.~Rassart \cite{BGR}, ... .
\begin{prb}
 Find a fermionic, i.e. a positive linear combination of products of powers 
of $t$ and $t$-binomial coefficients, formula for the polynomials 
${\cal E}_{\ld\mu}(t)$.
\end{prb}

 This problem should be a very difficult one, however, since, for example, the
polynomial 
$${\cal E}_{(n^n), ((n-1)^n,1^n)}(t)$$
coincides with the Ehrhart polynomial of the Birkhoff polytope ${\cal B}_n$ of
doubly stochastic matrices,
see \cite{Kir9}, Section~7.5. We refer the reader to a paper by M.~Beck and
D.~Pixton \cite{BP} and the literature quoted therein, for a further 
information about the Ehrhart polynomials          \\
( for $n \le 9$ ) and the volumes ( for $n \le 10$ ) of the Birkhoff 
polytope ${\cal B}_n.$ 
 
The (normalized) leading coefficient of  Ehrhart's 
polynomial ${\cal E}_{\ld\mu}(t)$ is equal to the (normalized)
volume of Gelfand--Tsetlin's polytope $GT(\ld ,\mu)$, and is known
in the literature, see e.g.  \cite{Hek}, \cite{Ok2}, as a {\it continuous}
 analog of the {\it weight multiplicity} $\dim V_{\ld}(\mu)$.

Finally, we would like to note that in general, the Ehrhart polynomial of a
convex integral polytope may have negative coefficients. The famous example
is {\it the Reeve tetrahedron}, see e.g. \cite{Kir9}, Example~7.34, {\bf 6}, 
and the literature quoted therein.

\subsection{Skew Kostka--Foulkes polynomials}

Let $\ld ,\mu$ and $\nu$ be partitions, $\ld\supset\mu$, and
$|\ld|=|\mu|+|\nu|$.

\begin{de}\label{d1.1} The skew Kostka--Foulkes polynomials
$K_{\ld\setminus\mu ,\nu}(q)$ are defined as the transition coefficients
from the skew Schur functions $s_{\ld\setminus\mu}(x)$ to
the Hall--Littlewood functions $P_{\nu}(x;q)$:
\begin{equation}
s_{\ld\setminus\mu}(x)=\sum_{\nu}K_{\ld\setminus\mu ,\nu}(q)P_{\nu}(x;q).
\label{2.3}
\end{equation}
\end{de}

It is clear that
$$K_{\ld\setminus\mu ,\nu}(q)=\sum_{\pi}c_{\mu\pi}^{\ld}K_{\pi\nu}(q),$$
where the coefficients $c_{\mu\pi}^{\nu}={\rm
Mult}[V_{\nu}:V_{\mu}\otimes V_{\pi}]$ stand for the
Littlewood--Richardson numbers.

Let us remark that 
\begin{equation}
K_{\ld\setminus\mu ,\nu}(q)=\sum_Tq^{c(T)}
\end{equation}
summed over all semistandard skew tableaux $T$ of shape $\ld\setminus\mu$ and
weight $\nu$, where $c(T)$ denotes the {\it charge} of a skew tableau $T$. \\
 In the case $\mu = \emptyset,$  the formula (2.9) is due to A.~Lascoux and
  M.-P.~Sch\"utzenberger  \cite{LS}. See also \cite{But},
Chapter II, for an 
extended exposition of \cite{LS}.  We refer the reader to \cite{Ma},
 Chapter III, Section 6, for the definition of the Lascoux--Sch\"utzenberger 
statistics {\it charge}
on the set of semistandard Young tableaux.

We will use  also the {\it cocharge} version of the skew Kostka--Foulkes
polynomials:
\begin{equation}
{\overline K}_{\ld\setminus\mu ,\nu}(q)=\sum_{\pi}c_{\mu\pi}^{\ld}
{\overline K}_{\pi\mu}(q), \label{2.5}
\end{equation}
where ${\overline K}_{\ld\mu}(q)=q^{n(\mu)}K_{\ld\mu}(q^{-1})$.

$(\spadesuit)$  We will see in Section 5.1, example $3^0$, that the skew 
Kostka-Foulkes polynomials are some special cases  of the parabolic Kostka 
polynomials.
 
\subsection{Littlewood--Richardson numbers and Saturation Theorem}
 
 The Littlewood--Richardson numbers ${c_{\ld,\mu}^{\nu}},$ $LR$-numbers for 
short, are defined as the
structural constants of the multiplication of Schur functions. More 
specifically,
let $\ld$ and $\mu$ be partitions, then
\begin{equation}
 s_{\ld}s_{\mu}=\sum_{\nu}c_{\ld,\mu}^{\nu}s_{\nu},
\end{equation}
or equivalently,
$$s_{\nu\setminus\mu}=\sum_{\ld}c_{\ld,\mu}^{\nu}s_{\ld}.$$
We have $c_{\ld,\mu}^{\nu}=0$ unless $|\nu|=|\ld|+|\mu|$ and 
$\nu\supset\ld,\mu.$  A pure combinatorial way to compute the $LR$-numbers
is given by the celebrated Littlewood--Richardson rule, see e.g. \cite{Ma},
Chapter~I, Section~9. 

\hskip -0.6cm{\bf Saturation Theorem} ( A.~Knutson and T.~Tao \cite{KT1} ) \\

 ~~~~~$c_{N\ld,N\mu}^{N\nu}\ne 0$ for some integer $N \ge 1$ if and only if
$c_{\ld,\mu}^{\nu} \ne 0.$
\vskip 0.3cm
We refer the reader to interesting and nice written papers by W.~Fulton
\cite{Fu2}, \cite{Fu3}  and A.~Zelevinsky \cite{Z2} for detailed account
to an origin of {\it Saturation Conjecture} ( now a theorem by A.~Knutson 
and T.~Tao ) and its connections with the so-called {\it Horn Problem}.

\subsection{Internal product of Schur functions, and polynomials
$L_{\al,\beta}^{\mu}(q)$}

 The irreducible characters $\chi^{\ld}$ of the symmetric group $\Sigma_n$ are
indexed in a natural way by partitions $\ld$ of $n$. If $w\in \Sigma_n$, then
define $\rho (w)$ to be the partition of $n$ whose parts are the cycle
lengths of $w$. For any partition $\ld$ of $m$ of length $l$, define the
power--sum symmetric function
$$p_{\ld}=p_{\ld_1}\ldots p_{\ld_l},$$
where $p_n(x)=\sum x_i^n$. For brevity write $p_w:=p_{\rho (w)}$.
The Schur functions $s_{\ld}$ and power--sums $p_{\mu}$ are
related by a famous result of Frobenius
\begin{equation}
s_{\ld}=\frac{1}{n!}\sum_{w\in \Sigma_n}\chi^{\ld}(w)p_w. \label{2.12}
\end{equation}
For a pair of partitions $\al$ and $\beta$, $|\al|=|\beta|=n$, let us
define the internal product  $s_{\al}*s_{\beta}$ of Schur functions
$s_{\al}$ and $s_{\beta}$:
\begin{equation}
s_{\al}*s_{\beta}=\frac{1}{n!}
\sum_{w\in \Sigma_n}\chi^{\al}(w)\chi^{\beta}(w)p_w.
\label{2.13}
\end{equation}
It is well--known, see e.g. \cite{Ma},Chapter I, Section 7, that
$$s_{\al}*s_{(n)}=s_{\al},~~ s_{\al}*s_{(1^n)}=s_{\al'},
$$
where $\al'$ denotes the conjugate partition to $\al$.

Let $\al ,\beta ,\gamma$ be partitions of a natural number $n\ge 1$,
consider the following numbers
\begin{equation}
g_{\al\beta\gamma}=\frac{1}{n!}\sum_{w\in \Sigma_n} \chi^{\al}(w)
\chi^{\beta}(w)\chi^{\gamma}(w). \label{2.14}
\end{equation}
The numbers $g_{\al\beta\gamma}$ coincide with the structural constants
for multiplication of the characters $\chi^{\al}$ of the symmetric group
$\Sigma_n$:
\begin{equation}
\chi^{\al}\chi^{\beta}=\sum_{\gamma}g_{\al\beta\gamma}\chi^{\gamma}.
\label{2.15}
\end{equation}
Hence, $g_{\al\beta\gamma}$ are {\it non--negative} integers. It is clear that
\begin{equation}
s_{\al}*s_{\beta}=\sum_{\gamma}g_{\al\beta\gamma}s_{\gamma}. \label{2.16}
\end{equation}
\begin{rem} {\rm  More generally, let $A$ and $B$ be two skew diagrams and 
$\gamma$ be a partition  all of the same cardinality $n$. Define the 
coefficients $g_{A,B,\gamma}$ and the {\it internal product} $s_A \ast s_B$ of 
skew Schur functions $s_A$ and $s_B$ as follows. Let $\chi^{A}$ and $\chi^{B}$ 
be the characters of representations of the symmetric group $\Sigma_n$ which  
correspond to the skew diagrams $A$ and $B$. The numbers $g_{A,B,\gamma}$ are 
defined via the decomposition
$$ \chi^{A}~\chi^{B}= \sum_{\gamma} g_{A,B,\gamma}~\chi^{\gamma}.
$$ 
The {\it internal product} of the skew Schur functions $s_A$ and
$s_B$ is defined as follows
$$ s_{A} \ast s_{B}= \sum_{\gamma} g_{A,B,\gamma}~s_{\gamma}.
$$
Finally, let $C$ be one more skew diagram, define the number $g_{A,B,C}$ to
be equal to $\langle s_{A}*s_{B},s_{C} \rangle,$  where $\langle ~,\rangle$
denotes the Redfield--Hall scalar product on the ring of symmetric functions,
 see \cite{Ma}, Chapter~I, Section~4.   } 
\end{rem}
\begin{rem} { \rm It is one of the most fundamental  open  problems  in 
Combinatorics and Representation Theory of the symmetric group that to find a  
combinatorial  rule for description  of the numbers $g_{\al\beta\gamma}.$ }
\end{rem}
\begin{theorem} Let $\al,\beta$ and $\gamma$ be partitions of the same size
$n.$ 

$(\clubsuit)$  If $g_{\al\beta\gamma} \ne 0,$ then $g_{N\al,N\beta,N\gamma}
 \ne 0$ for any integer $N \ge 1.$
\end{theorem}

\begin{rem} { \rm The converse statement, i.e.
 
if $g_{N\al,N\beta,N\gamma} \ne 0$ for some integer $N \ge 2,$ then 
$g_{\al\beta\gamma} \ne 0,$ 

the so-called {\bf saturation property} of the 
structural constants $g_{\al\beta\gamma}$, is {\bf not true} in general if 
$n \ge 7,$ even
under the  additional assumption that partitions $\al,\beta,\gamma$ and 
their conjugate ones  $\al',\beta',\gamma'$, all have at least two different 
parts. For example,

$g_{(6,1),(4,1^3),(3,3,1)}=0,$ but $g_{(12,2),(8,2^3),(6,6,2)} \ge 1,$
$g_{(5,2),(4,3),(4,1^3)}=0,$ but $g_{(10,4),(8,6),(8,2^3)} \ge 1,$

$g_{(6,1^2),(6,1^2),(4,3,1)}=0,$ but $g_{(12,2^2),(12,2^2),(8,6,2)} \ge 1,$
$g_{(6,2),(6,1^2),(4,2^2)}=0,$ 

but $g_{(12,4),(12,2^2),(8,4^2)} \ge 1.$
 

On the other hand,

$g_{(3,1,1),(3,2),(2,1^3)}= 1$ and $g_{(6,2,2),(6,4),(4,2^3)}=2,$ 
$g_{(2,1),(2,1),(1^3)}=1$ and $g_{(4,2),(4,2),(2^3)}=1,$

$g_{(2,2),(2,2),(2,2)}=1$ and $g_{(4,4),(4,4),(4,4)}=1,$
$g_{(2,2),(2,2),(1^4)}=1$ and $g_{(4,4),(4,4),(2^4)}=1.$
  
$(\maltese)$ However, we {\bf expect} that the formal power series
$$ \sum_{N \ge 1} g_{N\al,N\beta,N\gamma} ~t^N $$
is a {\it rational} function of $t$  
~(with the only possible pole at $t=1$ ??).   }
\end{rem}
\begin{rem} { \rm  ( {\bf Plethysm structural constants} )

 Fix integer numbers $k$ and $n \ge 2$, and a finite dimensional 
representation 
$W$ of the Lie algebra $\g l_n$.  The $k$-th tensor power $W^{\otimes k}$ of
the $\g l_n$-module $W$ has a natural structure of 
$\Sigma_k \times {\g}l_n$-module,
where $\Sigma_k$ denotes the symmetric group of order $k!$. Let
\begin{equation}
W^{\otimes k}=\sum_{\mu,\pi} a_{\mu,W}^{\pi} ~S^{\mu} \otimes V_{\pi}
\end{equation}
be the decomposition of the module $W^{\otimes k}$ into irreducible
$\Sigma_k \times {\g}l_n$-submodules. Here $\mu$ is a partition of size $k,$ 
and $S^{\mu}$ stands for the irreducible representation of the symmetric group
$\Sigma_k$ which corresponds to the partition $\mu$; $\pi$ is a partition of
length at most $n$ and $V_{\pi}$ denotes the irreducible $\g l_n$-module
with the highest weight $\pi$.

 If $W=V_{\ld}$ is the irreducible $\g l_n$-module with the highest weight
$\ld$, then the numbers  $a_{\ld,\mu}^{\pi}:=a_{\mu,V_{\ld}}^{\pi}$ coincide
with the structural constants of yet another multiplication, called {\it 
plethysm}, in the ring of symmetric functions $\Lambda$:
$$ s_{\ld} \circ s_{\mu}=\sum_{\pi}a_{\ld,\mu}^{\pi}~s_{\pi}. $$
 Note, that the plethysm is an associative, but not commutative operation.

It is well-known, see e.g. \cite{Stm}, that if $\al$ and $\beta$ are 
partitions of the same size $k$ such that $l(\al)=r$, $l(\beta)=s$ 
and $n \ge r+s$, and
furthermore, $W={\g}l_n$ is the adjoint representation, and 
$$ \pi= (k+\al_1,\cdots,k+\al_r,\underbrace{k \ldots,k}_{n-r-s},k-\beta_s,
\cdots,k-\beta_1), $$
then 
$$ a_{\mu,\g l_n}^{\pi}:=[ S^{\mu} \otimes V_{\pi}: {\g}l_n^{\otimes k}] =
g_{\al\beta\mu}. $$
Hence, the inner product structure constants $g_{\al\beta\gamma}$, and
therefore the $LR$-numbers, see Section~2.6, are certain special cases of the 
plethysm  structural constants  $a_{\mu,W}^{\pi}.$  
\begin{con}
Let $\mu$ and $\pi,$ $l(\pi) \le n,$ be partitions such that $\mu$ has at
least two different parts. Let $W$ be a finite dimensional $\g l_n$-module. 

If $a_{\mu,W}^{\pi} \ne 0,$ then $a_{N\mu,W}^{N\pi} \ne 0,$ 
for any  integer $N \ge 1.$
\end{con}

$(\maltese)$ Moreover, we {\bf expect} that if $N_1$ and $N_2$ are integers 
such that $N_1 \ge N_2,$ then 

$a_{N_1 \mu,W}^{N_1 \pi} \ge a_{N_2 \mu,W}^{N_2 \pi},$ 
 and the formal power series

$$\sum_{N \ge 1} a_{N\mu,W}^{N\pi}~t^N$$ 

is a {\it rational} function of $t$ 
~(with the only possible pole at $t=1$ ??). 

$(\clubsuit)$  We want to emphasize that the plethysm structural constants 
$a_{\mu,W}^{\pi}$ do not satisfy  the so-called {\bf saturation property},
 i.e. it's not true, in general, that if $a_{N\mu,W}^{N\pi} \ne 0$ 
for some integer $N \ge 2,$ then $a_{\mu,W}^{\pi} \ne 0.$

Using the tables of plethysms from \cite{Ag}, we have checked that

$a_{(2,2),(4,2)}^{(6,4^2,2^5)}=1,$ but $a_{(2,2),(2,1)}^{(3,2^2,1^5)}=0,$
$a_{(2,2),(4,2)}^{(4^5,2^2)}=1,$ but $a_{(2,2),(2,1)}^{(2^5,1^2)}=0.$

$(\maltese)$ Based on several examples, we {\bf expect} that if 
$a_{2\mu,W}^{2\pi} \ge 2,$  then $a_{\mu,W}^{\pi} \ne 0.$

On the other hand, Conjecture~2.14 is not true if a partition $\mu$ has
a form $(1^k)$. For example,

$a_{(2,1,1),(1,1,1)}^{(4,4,2,1,1)}=1,$ but $a_{(2,1,1),(2,2,2)}^{(8,8,4,2,2)}=0,$  $a_{(2,1,1),(1,1,1)}^{(4,3,3,1,1)}=0,$ but $a_{(2,1,1),(2,2,2)}^{(8,6,6,2,2)}=1.$   } 
\end{rem} 
\begin{quest}  Could it be true that for any finite dimensional 
${\g}l_n$-module $W$ there exists a polynomial $p_{W}(t)$ ( $p_{W}(t)=t$ ?? )
 such that for all partitions $\pi$ and $\mu$ one has 

 if $a_{N\mu,W}^{N\pi} \ge p_{W}(N),$  then $a_{\mu,W}^{\pi} \ne 0.$  
\end{quest}   

$(\spadesuit)$  { \it It is one of the most fundamental problems of Algebraic 
Combinatorics, Representation Theory, Theory of Invariants, ... that
to find a combinatorial rule for description of the numbers 
$a_{\mu,W}^{\pi}.$ }  

\qed

\begin{de} The  polynomials $L_{\al\beta}^{\mu}(q)$ are defined via the
decomposition of the internal product of Schur functions
$s_{\al}*s_{\beta}(x)$ in terms of the Hall--Littlewood functions:
\begin{equation}
s_{\al}*s_{\beta}(x)=\sum_{\mu}L_{\al\beta}^{\mu}(q)P_{\mu}(x;q).
\label{2.17}
\end{equation}
\end{de}
In a similar fashion one can define the polynomials $L_{A,B}^{\mu}(q)$,
where $A$ and $B$ are skew diagrams and $\mu$ is a partition:
$$ s_{A} \ast s_{B}(x) = \sum_{\mu} L_{A,B}^{\mu}(q)~P_{\mu}(x;q).
$$

\begin{exs}
 $(i)$ Take $n=4$, $\al=(3,1)$ and $\beta=(2,2)$.    

Then the all non--zero polynomials  $L_{(3,1),(2,2)}^{\mu}(q)$ are:

$$L_{(3,1),(2,2)}^{(3,1)}(q)=1,~L_{(3,1),(2,2)}^{(2,2)}(q)=q,
~L_{(3,1),(2,2)}^{(2,1,1)}(q)=1+q+q^2,$$
$$L_{(3,1),(2,2)}^{(1^4)}(q)=q(1,1,2,1,1).$$
 $(ii)$  Take $n=6$ and $\al = \beta = (3,2,1),$ then
$$ L_{\al,\beta}^{(6)}(q)=1,~L_{\al,\beta}^{(5,1)}(q)=2+q,
~L_{\al,\beta}^{(4,2)}(q)=(3,2,1),~L_{\al,\beta}^{(4,1,1)}(q)= (4,5,2,1)
=(1+q)(4,1,1),$$
$$ L_{\al,\beta}^{(3,1^3)}(q)= (4,9,12,11,5,2,1)=(1+q)(4,5,7,4,1,1),
~L_{\al,\beta}^{(2,1^4)}(q)= (1+q)^2~(1+q^2)^2~(2,3,0,1).$$ 
$(iii)$ Take $n=6,$ ~$\al = (4,2)$ and $ \beta = (3,2,1),$ then
$$ L_{\al,\beta}^{(5,1)}(q)=1,~L_{\al,\beta}^{(4,2)}(q)=2+q,
~L_{\al,\beta}^{(4,1,1)}(q)=(2,3,1),~L_{\al,\beta}^{(3,3)}(q)=(1,2,1),$$
$$L_{\al,\beta}^{(3,1^3)}(q)=(1+q)(1+q+q^2)(2,1,1), 
~L_{\al,\beta}^{(2,2,1,1)}(q)=(1+q)^2(1+q+q^2)(2,0,1),$$
$$ L_{\al,\beta}^{(2,1^4)}(q)= (1+q)^2(1+q+q^2)(1+q+q^2)(1,1,0,1).$$
 $(iv)$ Take $n=6$, $\al=(4,2)$ and $\beta=(2^3)$. Then
$$ L_{\al,\beta}^{(4,1,1)}(q)=q, L_{\al,\beta}^{(3,3)}(q)=q, 
L_{\al,\beta}^{(3,2,1)}(q)=1+q+q^2,
$$
$$ L_{\al,\beta}^{(2,2,1,1)}(q)=q(3,2,3,1,1), L_{\al,\beta}^{(2,1^4)}(q)=(1,1,1)(1,0,2,1,2,0,1),
$$
$$ L_{\al,\beta}^{(1^6)}(q)= (1,0,1,1,0,1)~{\wt K}_{\al,(1^6)}(q).
$$
\end{exs}   
Hereafter we shell use the notation ${\wt K}_{\al,\mu}(q)$ to denote the
polynomial $q^{n(\mu)-n(\al)}~K_{\al,\mu}(q^{-1}).$ 
\begin{rem} {\rm It is not true in general that if $\al,\beta,\mu$ are 
partition and $\al \ge \mu,$  then  the ratio  
$L_{\al,\beta}^{\mu}(1)/K_{\al,\mu}(1) \in \Z$.

For example, take $\al=\beta=(6,2,1)$ and $\mu=(3,3,2,1)$. Then
$$L_{\al,\beta}^{\mu}(q)=(2,17,44,63,64,48,29,15,6,2,1), ~{\wt K_{\al,\mu}(q)}
=(1,2,2,1)$$
and $L_{\al,\beta}^{\mu}(1)=291,$  $L_{\al,\beta}^{\mu}(-1)=1.$  

We see that ${\wt K_{\al,\mu}(q)}$ is not
a divisor of $L_{\al,\beta}^{\mu}(q)$, and  the ratio
$L_{\al,\beta}^{\mu}(1)/K_{\al,\mu}(1) \notin \Z$. Note that 
$L_{\al,\beta}^{\mu}(0)=c_{(2,1),(2,1)}^{(3,2,1)}=2$ and 
$deg L_{\al,\beta}^{\mu}(q)=10=n(\mu)$ in a good agreement with 
Conjecture~2.23. }
\end{rem}

It follows from (\ref{2.6}) and (\ref{2.16}) that
\begin{equation}
L_{\al\beta}^{\mu}(q)=\sum_{\gamma}g_{\al\beta\gamma}K_{\gamma\mu}(q).
\label{2.18}
\end{equation}
Thus, the polynomials $L_{\al\beta}^{\mu}(q)$ have non--negative integer
coefficients, and
$$L_{\al\beta}^{\mu}(0)=g_{\al\beta\mu}.$$

It follows from (\ref{2.17}) that the number $L_{\al,\beta}^{\mu}(1)$ is equal 
to $\langle s_{\al}*s_{\beta},h_{\mu} \rangle,$ where 
$\langle ~, \rangle$ denotes 
the Redfield--Hall scalar product on the ring of symmetric functions, 
see \cite{Ma}. In other words,
$$ s_{\al}(x) \ast s_{\beta}(x) = \sum_{\mu} L_{\al,\beta}^{\mu}(1) m_{\mu}(x),
$$
where $m_{\mu}(x)$ denotes the {\it monomial} symmetric function corresponding
to partition $\mu.$ Therefore, the numbers $L_{\al,\beta}^{\mu}(1)$ and 
$L_{A,B}^{\mu}(1)$ can be defined for any composition $\mu.$

\begin{rem} {\rm There is a well-known connection between the structural 
constants $g_{\al\beta\gamma}$ and the numbers $L_{\al,\beta}^{\mu}(1).$ 
Namely, let $A,$ $B$ and $C= \Gamma \setminus \gamma$ be skew diagrams such 
that the partition $\Gamma$ has 
the length at most $n,$  and $|A|=|B|=|C|.$  Then
$$ g_{A,B,C}= \sum_{w \in \Sigma_n} (-1)^{l(w)}~L_{A,B}^{w \circ ~C}(1),
$$ where $w \circ C$ stands for the composition 
$w(\Gamma+\delta_n)- \gamma - \delta_n,$ and 
$\delta_n=(n-1,n-2,\cdots,1,0).$   }
\end{rem}
The polynomials $L_{\al\beta}^{\mu}(q)$ can be considered as a
generalization of the Kostka--Foulkes polynomials. Indeed, if
partition $\beta$ consists  of one part, $\beta =(n)$, then
$$L_{\al \beta}^{\mu}(q)=K_{\al ,\mu}(q),~~
L_{\al \beta'}^{\mu}(q)=K_{\al' ,\mu}(q).
$$
\begin{pr} Let $\al,$ $\beta$ and $\mu=(\mu_1 \ge \cdots \ge \mu_r)$ 
be partitions of the same size $n.$ Then
\begin{equation}
 L_{\al,\beta}^{\mu}(1)= \ds\sum_{\sbmu}K_{\al,~\sbmu}(1)K_{\beta,~\sbmu}(1),
\label{2.19}
\end{equation}
where the sum runs over sequences of partitions $\bmu=(\mu^{(1)},\cdots,
\mu^{(r)})$ such that $|\mu^{(a)}|=\mu_a$, $1 \le a \le r$.
\end{pr}

\begin{cor} If $\mu=(r,1^s)$ is a hook partition, then
\begin{equation}
 L_{\al,\beta}^{\mu}(1)= \sum_{|\ld|=r} K_{\al\setminus\ld,~(1^s)}(1)
~K_{\beta\setminus\ld,~(1^s)}(1), \label{2.20}
\end{equation}
where the sum runs over all partitions $\ld$ of size $r,$ 
$\ld \subset \al \cap \beta.$  
\end{cor}
In particular, $L_{\al \beta}^{(1^n)}(1)= f^{\al}f^{\beta},$
where $f^{\al}$ denotes the number of standard Young tableaux of shape $\al.$
More generally \cite{Kir9},
\begin{equation}
L_{\al\beta}^{(1^n)}(q)=K_{\beta'\al}(q,q){\wt K}_{\al ,(1^n)}(q)=
K_{\al'\beta}(q,q){\wt K}_{\beta ,(1^n)}(q),
\end{equation}
where
$$ {\wt K}_{\al ,\beta}(q):=q^{n(\beta)-n(\al)}~K_{\al\beta}(q^{-1}), 
~~K_{\al\beta}(q,q):=K_{\al\beta}(q,t)\vert_{t=q},$$
and $K_{\al\beta}(q,t)$ stands for the double Kostka polynomial
introduced by I.~Macdonald \cite{Ma}, Chapter~VI, (8.11).

\begin{prb} Find a $q$--analog of the equality (\ref{2.20}) .
\end{prb}

\begin{con} ( {\bf Saturation conjecture for polynomials 
$L_{\al,\beta}^{\mu}(q)$ } )
 
 Let $\al$,~$\beta$ and $\mu$ be partitions of the same size such that 
$L_{\al,\beta}^{\mu}(q) \ne 0.$  Then

$(\blacklozenge)$  For any integer $N \ge 1,$

$\bullet$ $\max deg L_{N\al,~N\beta}^{N\mu}(q)= 
N~\max deg L_{\al,\beta}^{\mu}(q);$
 
$\bullet$ If partition $\mu$ either has at least two different parts, or
$\mu$ has a rectangular shape, but $\mu$ is different from the both partitions
$\al$ and $\beta,$ and their conjugate ones $\al'$ and $\beta',$ then

$\min deg L_{N\al,~N\beta}^{N\mu}(q)= N~\min deg L_{\al,\beta}^{\mu}(q).$

$(\blacklozenge \blacklozenge)$  $\max deg L_{\al,\beta}^{\mu}(q)=n(\mu)- A(\al,\beta),$
where $A(\al,\beta)$ stands for the $\min deg K_{\al,\beta}(q,q),$ i.e. 

$K_{\al,\beta}(q,q)= B(\al,\beta)q^{A(\al,\beta)}$+ higher degree terms.

$(\blacklozenge \blacklozenge \blacklozenge)$ ( {\bf  Saturation conjecture 
for polynomials $K_{\al,\beta}(q,q)$} )

For any integer $N \ge 1,$   ~~~$A(N\al,~N\beta)=N~A(\al,\beta).$
\end{con}

\begin{exs} 
 $(i)$ Take $n=3$, 
$$ L_{(2,1),(2,1)}^{(3)}(q) =1,~~L_{(4,2),(4,2)}^{(6)}(q) =1,$$
$$ L_{(2,1),(2,1)}^{(1^3)}(q)=1+q+q^2+q^3,
~L_{(4,2),(4,2)}^{(2^3)}(q)=1+2~q+4~q^2+3~q^3+3~q^4+q^5+q^6,$$
$$ L_{(2,1),(2,1)}^{(2,1)}(q)=1+q,~~L_{(4,2),(4,2)}^{(4,2)}(q)=2+q+q^2.$$
$$ L_{(2,1),(1^3)}^{(2,1)}(q)=1,~~L_{(4,2),(2^3)}^{(4,2)}(q)=1,
$$
$$ L_{(2,1),(1^3)}^{(1^3)}(q)=q+q^2,~~L_{(4,2),(2^3)}^{(2^3)}(q)=(1,1,2,1,1).
$$
 $(ii)$ Take $n=4$,
$$ L_{(3,1),(2,2)}^{(2,2)}(q)=q,~~L_{(6,2),(4,4)}^{(4,4)}(q)=1+q^2, $$
$$ L_{(3,1),(2,2)}^{(2,1,1)}(q)=1+q+q^2, ~~L_{(6,2),(4,4)}^{(4,2,2)}(q)= 
(1,2,3,1,1), $$
$$L_{(3,1),(2,2)}^{(1^4)}(q)=q(1,1,2,1,1),~~L_{(6,2),(4,4)}^{(2^4)}(q)=
q^2(2,2,6,5,7,4,4,1,1). $$
\end{exs}
The latter example shows that for the numbers $g_{\al\beta\gamma}$ an obvious
generalization of the Fulton conjecture, see Section~6, is false.
\begin{con} ( {\bf Rationality conjecture for polynomials 
$L_{\al,\beta}^{\mu}(q)$} )

Let $\al,\beta$ and $\mu$ be partitions of the same size. The generating
function
$$ \sum_{N \ge 0} L_{N\al,N\beta}^{N\mu}(q)~t^{N}$$
is a rational function of $q$ and $t.$ 
\end{con}
\begin{prb}
Give a combinatorial interpretation of the integer numbers 
$L_{\al,\beta}^{\mu}(-1)$.
\end{prb} 
 
\begin{prb} Find a fermionic type formula for the polynomials 
$L_{\al \beta}^{(\mu)}(q)$  which extends that for the Kostka--Foulkes
polynomials, see Section~5.1, Theorem~5.3.
\end{prb}

\subsection{Extended and Restricted Littlewood--Richardson numbers}

$(1^0)$ ( {\bf Extended Littlewood--Richardson numbers} )

Let $\ld$,$\mu$ and $\nu$ be partitions such that $|\ld|+|\mu|\ge\ |\nu|.$
Choose an integer number $N$ such that $N \ge N_0:=
\max(|\ld|+\ld_1,|\mu|+\mu_1,
|\nu|+\nu_1)$, and consider partitions $$\al_N:=(N-|\ld|,\ld),
~\beta_N:=(N-|\mu|,\mu),~\gamma_N:=(N-|\nu|,\nu).$$ It is clear that
$|\al_N|=|\beta_N|=|\gamma_N|=N$.

According to results by F.~Murnaghan \cite{Mur}, Y.~Dvir \cite{Dvi}  and E.
Vallejo \cite{Val},
if $N \ge N_0$, then the number $g_{\ld_N,\beta_N,\gamma_N}$ does not depend
 on $N$.
\begin{de} The extended Littlewood--Richardson number ~$C_{\ld,\mu}^{\nu}$ is
defined to be equal to the stable value of the numbers 
$g_{\ld_N,\beta_N,\gamma_N}.$
\end{de}
More generally, the following statement is true:
\begin{pr}
 The sequence of polynomials $\{L_{\al_{N},\beta_{N}}^
{\gamma_N}(q) \}_{N \ge 1}$ is stabilized to the polynomial  
${\cal L}_{\ld,\mu}^{\nu}(q),$ i.e. if integer $N$ is big enough, then the
polynomial  $L_{\al_{N},\beta_{N}}^{\gamma_N}(q)$  does not
depend on $N$ and equal to ${\cal L}_{\ld,\mu}^{\nu}(q).$
The latter is a polynomial with non--negative integer coefficients, and 
${\cal L}_{\ld,\mu}^{\nu}(0)=C_{\ld,\mu}^{\nu}.$ 
\end{pr}
According to another result by Y.~Dvir \cite{Dvi}, the numbers 
$C_{\ld,\mu}^{\nu}$ can be considered as  a 
generalization of the $LR$-numbers $c_{\ld,\mu}^{\nu}.$
\begin{pr} ( Y.~Dvir \cite{Dvi} )  If $|\ld|+|\mu|=|\nu|,$ then
the number $C_{\ld,\mu}^{\nu}$ coincides with the Littlewood--Richardson 
number $c_{\ld,\mu}^{\nu}.$
\end{pr}

\begin{exs}  $(i)$ Take $\ld=\mu=(2,1)$, then
$$ C_{\ld,\mu}^{(3,2,1)}=c_{\ld,\mu}^{(3,2,1)}=2,$$
$$C_{\ld,\mu}^{(3,1,1)}=6,~C_{\ld,\mu}^{(2,2,1)}=5,~ C_{\ld,\mu}^{(2,1,1,1)}=4,
~C_{\ld,\mu}^{(3,2)}=5,$$
$$ C_{\ld,\mu}^{(2,2)}=6,~C_{\ld,\mu}^{(3,1)}= ~C_{\ld,\mu}^{(2,1,1)}=9, 
~C_{\ld,\mu}^{(2,1)}=9.$$
  $(ii)$ Take $\ld=(2,1)$ and $\mu=(3,1),$ then
$$ C_{\ld,\mu}^{(3,1)}=13, ~C_{\ld,\mu}^{(2,1)}= 9.$$ 
\end{exs}

\begin{prb} Find a combinatorial rule for calculating the extended 
$LR$-numbers $C_{\ld,\mu}^{\nu}$ which extends the Littlewood--Richardson rule.
\end{prb}

\begin{con} ( {\bf  Saturation conjecture for extended $LR$-numbers } ) 

$C_{N\ld,N\mu}^{N\nu} \ne 0$ for some integer $N \ge 1$ if and only if
~~$C_{\ld,\mu}^{\nu} \ne 0.$
\end{con}

$(2^0)$ ( {\bf Restricted Littlewood--Richardson numbers, cf. 
\cite{Kac}, Exercise~13.35} )

Fix positive integers $l$ and $n \ge 2$. Denote by $\Sigma_{n,l}$ the affine 
reflection group on $\R^{n}$ generated by the reflection
$$ s_{0}=(x_n+l,x_2,\cdots,x_{n-1},x_1-l)
$$
and the symmetric group $\Sigma_n$.
\begin{de} ( {\bf Restricted Littlewood--Richardson numbers}  )

Let $\ld,\mu$ and $\nu$ be partitions such that  $|\ld|+|\mu|=|\nu|$. 
Define the {\it level $l$ } restricted
Littlewood--Richardson number $c_{\ld,\mu}^{\nu}[l]$ as follows
$$ c_{\ld,\mu}^{\nu}[l]=\sum_{w \in \Sigma_{n,l}} (-1)^{l(w)} c_{\ld,\mu}^{w
\circ \nu},
$$
where $w \circ \nu$ denotes the composition  $ w(\nu+\delta_n) - \delta_n,$ and
$\delta_n=(n-1,\cdots,1,0).$
\end{de}
It is well--known that
$$ 0 \le c_{\ld,\mu}^{\nu}[1] \le c_{\ld,\mu}^{\nu}[2] \le \cdots = 
c_{\ld,\mu}^{\nu}.
$$
In a similar fashion one can define the {\it  level $l$} extended 
Littlewood--Richardson numbers $C_{\ld,\mu}^{\nu}[l].$
\begin{con} ( {\bf Saturation conjecture for the level $l$ extended 
$LR$-numbers} )

Let $\ld,\mu$ and $\nu$ be partitions such that $|\ld|+|\mu| \ge |\nu|.$ Then

$C_{N\ld,N\mu}^{N\nu}[l] \ne 0$ for some integer $N \ge 1$ if and only if
$C_{\ld,\mu}^{\nu}[l] \ne 0.$
\end{con}
\begin{con} ( {\bf Polynomiality conjecture for level $l$ extended 
$LR$-numbers} )

Let $\ld,\mu$ and $\nu$ be partitions such that $|\ld|+|\mu| \ge |\nu|.$ Then

$C_{N\ld,N\mu}^{N\nu}[l]$ is a polynomial in $N$ with {\bf non--negative}
 rational coefficients.

$(\clubsuit)$ Moreover , the formal power series 
$$ \sum_{N \ge 0} C_{N\ld,N\mu}^{N\nu}[l]~t^N $$
is a rational function in $t$ of the form
$$P_{\ld,\mu}^{\nu,l}(t)/(1-t)^{r(\ld,\mu,\nu,l)+1}, 
~~P_{\ld,\mu}^{\nu,l}(0)=1, ~~P_{\ld,\mu}^{\nu,l}(1) \ne 0,$$ 
where $r(\ld,\mu,\nu,l) \in \Z_{ \ge 0}$ and $P_{\ld,\mu}^{\nu,l}(t)$ 
is a polynomial with non--negative integer  coefficients.
\end{con}

\section{Parabolic Kostant partition function and its $q$-analog}

\subsection{Definitions: algebraic and combinatorial}
 Let $\eta = (\eta_1,\eta_2,\cdots,\eta_p)$ be a  composition, $\eta_p > 0$, 
 $|\eta| = n.$ Denote by
$\Phi(\eta)$ the set of ordered pairs $(i,j)\in\Z^2$ such that
\begin{equation}
1\le i \le\eta_1+\cdots +\eta_r<j\le n   \label{2.2}
\end{equation}
for some $r$, $1\le r\le p$. For example, if $\eta =(1^n)$, then
$$\Phi(\eta)=\{(i,j)\in\Z^2|1\le i<j\le n\}.
$$
\begin{de} 
Let $\gamma =(\gamma_1,\gamma_2,\cdots,\gamma_n)\in\Z^n$ be a sequence of 
integers such that $|\gamma|=0$, define a  {\it parabolic 
$q$-Kostant partition
function}  $K_{\Phi(\eta)}(\gamma|~q)$ via  the decomposition
\begin{equation}
\prod_{(i,j)\in\Phi(\eta)}(1-qx_i/x_j)^{-1}
= \sum_{\gamma }K_{\Phi(\eta)}(\gamma |~q)x^{\gamma}, \label{2.3}
\end{equation}
where the sum runs over the all sequences $\gamma
=(\gamma_1,\gamma_2,\cdots,\gamma_n) \in \Z^n$ such that $|\gamma|=0.$ 
\end{de}

\begin{de} Let  $K_{\Phi(\eta)}(\gamma)$ denote the parabolic Kostant 
partition function, that is to say, the value of the polynomial 
$K_{\Phi(\eta)}(\gamma |~q)$ at $ q:=1.$
\end{de}

\begin{rem} ( {\bf Combinatorial definition of $q$-Kostant partition 
function} )

{\rm One can give an equivalent pure combinatorial definition
of the parabolic $q$-Kostant partition function $K_{\Phi(\eta)}(\gamma|~q)$
 as follows.

Let  $\eta$ be a composition, $|\eta|=n$. Denote by 
$SM_{\eta}(\gamma)$ the set of all skew--symmetric integer matrices
 $M=(m_{i,j})_{1 \le i,j \le n}$ such that      \\
$(i)$  $m_{i,j} \ge 0,$ ~if $1 \le i \le j \le n $ ;   \\
$(ii)$  $m_{i,j}=0,$ ~if $r_{k-1}< i \le j \le r_k$  for some $k,$
 $1 \le k \le p,$
where $r_k:=\sum_{j \le k} \eta_j,$ and $r_0:=0 ;$  \\
$(iii)$  $\sum_{j=1}^n m_{i,j} = \gamma_i,$ for all $i, 1 \le i \le n.$ \\

For each $M \in SM_{\eta}(\gamma)$ we define the {\it magnitude} of $M,$
denoted by $|| M||,$ to be the  sum 
$\sum_{1 \le i \le j \le n}m_{i,j}.$  Then
\begin{equation}
 K_{\Phi(\eta)}(\gamma |~q)=\sum q^{||M||},
\end{equation}
where the sum runs over  all  matrices $M \in SM_{\eta}(\gamma).$  }

\end{rem}
Therefore, $K_{\Phi(\eta)}(\gamma)= Card~| SM_{\eta}(\gamma)|.$ 

\begin{rem} ( { \bf Generalized $q$-Kostant partition function \cite{PS}} )

{ \rm Let $\mit\Sigma \subseteq \Phi(1^n)$ be a subset, 
following  \cite{PS}
one can define the {\it generalized} Kostant partition function 
$K_{\mit\Sigma}(\gamma)$ and its $q$-analog $K_{\mit\Sigma}(\gamma |~q),
$ from the decomposition
$$ \prod_{(i,j) \in \mit\Sigma}(1-q~x_i/x_j)^{-1} =\sum_{\gamma} 
K_{\mit\Sigma}(\gamma |~q)x^{\gamma},
$$
where the sum runs over all sequences $\gamma \in \Z^n$ such
 that $|\gamma|=0.$  Moreover, by definition,
$K_{\mit\Sigma}(\gamma)=K_{\mit\Sigma}(\gamma |~q) \vert_{q=1}.$ 

Equivalently, 
$$K_{\mit\Sigma}(\gamma |~q)= \sum_{M} q^{||M||},$$
where the sum runs over the set of $n$ by $n$ skew--symmetric matrices 
$M=(m_{i,j})$  such that

$(i)$ $m_{i,j} \ge 0$ if $ 1 < i \le j \le n,$

$(ii)$ $m_{i,j}=0$ if $ (i,j) \notin \mit\Sigma,$
 
$(iii)$ $\sum_{j}m_{i,j}=\gamma_{i}$ for all i, $ 1 \le i \le n.$

$(\spadesuit)$ Most of our results about the parabolic $q$-Kostant partition 
function  $K_{\Phi(\eta)}(\gamma |~q)$, including, for example, 
Theorems~3.17, 3.20, 3.23, 3.25, 3.30 and 3.31, with  a small modifications, 
are still 
valid  for the function  $K_{\mit\Sigma}(\gamma |~q).$  Since we don't use the 
generalized Kostant partition function in  the present paper, we leave this 
interesting subject for subsequent publications.}  
\end{rem}

\subsection{Elementary properties, and explicit formulas for $l(\eta) \le 4$}

Using the above combinatorial definition of the function 
$K_{\Phi(\eta)}(\gamma |~q),$ one can describe some elementary, but useful,
properties of the latter.
\begin{pr}
$(i)$ Let $\eta_i$, $i=1,~2,$ be two compositions
 and $\gamma_i \in Y_{\eta_i},$  $i=1,~2,$ then
$$K_{\Phi(\eta_1*\eta_2)}(\gamma_1*\gamma_2 |~q)=K_{\Phi(\eta_1)}(\gamma_1 |~q)~K_{\Phi(\eta_2)}(\gamma_2 |~q). $$
$(ii)$ Let $\eta$ be a composition and $\gamma \in Y_{\eta},$ then
\begin{equation}
K_{\Phi(\eta)}(\gamma | ~q)=K_{\Phi\left(\overleftarrow{\eta}\right)}
\left(- \overleftarrow{\gamma} |~q\right),
\end{equation}
where for any composition $\beta=(\beta_1,\cdots,\beta_{r-1},\beta_{r})$  the
symbol $\overleftarrow{\beta}$ stands for the composition 
$(\beta_{r},\beta_{r-1},\cdots,\beta_1).$

$(iii)$ Let $\eta_1$ and $\eta_2$ be compositions such that $\eta_2$ is a 
subdivision of $\eta_1,$ so that $\eta_1 \ge \eta_2$.  Then
$$ K_{\Phi(\eta_1)}(q) \le K_{\Phi(\eta_2)}(q).$$
\end{pr} 
See Section~1, Notation, for the definition when a 
composition $\eta_2$ is a subdivision of that $\eta_1.$
We remark that the last statement is false if one assumes only that 
$\eta_1 \ge \eta_2$ with respect to the dominance partial ordering on the
 set of compositions, see example below.
 
\begin{ex} Take $\gamma=(3,0,-1,-1,0,-1),$ then $(2,3,1) \ge (2,2,2),$ but
$$K_{\Phi(2,2,2)}(\gamma |~q)=q^3+2~q^4 \le K_{\Phi(2,3,1)}(\gamma |~q)=
q^3+3~q^4.$$ On the other hand, 
$K_{\Phi(2,1,1,2)}(\gamma |~q)=q^3(1,3,2,1) \ge K_{\Phi(2,2,2)}(\gamma |~q).$
\end{ex}
\begin{pr} ( {\bf Recurrence relation for parabolic $q$-Kostant partition 
function} )

Let $\eta =(\eta_1,\cdots,\eta_p)$ be a composition of size $n$,  $\gamma \in 
Y_{\eta}$. Define ${\wt\eta}=(\eta_1,\cdots,\eta_{r_{p-1}})$.  Then
\begin{equation}
K_{\Phi(\eta)}(\gamma |~q)=q^{-\gamma_n} \sum_{\beta}K_{\Phi({\wt\eta})}
(\gamma_1 - \beta_1,\cdots,\gamma_{r_{p-1}}-\beta_{r_{p-1}} |~q),
\end{equation}
where the sum runs over $\beta \in \Z_{\ge 0}^{r_{p-1}}$ such that
$|\beta|=-\gamma_n$. 
\end{pr}
The next proposition describes  several particular cases of Theorem~3.31
below, namely, the cases  when a  parabolic $q$-Kostant partition function 
admits an explicit simple expression.
\begin{pr} ( {\bf Explicit formulas for $l(\eta) \le 4$ } )  \\ 
$(i)$  Let $\eta=(\eta_1,\eta_2)$ be a two component composition and 
$\gamma \in Y_{\eta}.$  Let us introduce integer vectors $\ld=(\gamma_1,\cdots,
\gamma_{\eta_{1}})$ and $\mu = (- \gamma_{\eta_{1}+1},\cdots,- \gamma_{\eta_1+\eta_2}).$
Then $\ld$ and $\mu$ are compositions of the same size, and
\begin{equation}
 K_{\Phi(\eta)}(\gamma |~q)=|{\cal P}_{\ld\mu}|~ q^{ |\ld|},
\end{equation}
where ${\cal P}_{\ld\mu}$ denotes the set of transportation matrices of type $(\ld ; \mu),$
i.e. the set of $l(\ld)$ by $l(\mu)$ matrices with non--negative integer 
entries, and the row sums $\ld_{i},$  and the column sums $\mu_{j}.$ 

$(ii)$ Let $\eta=(1^3)$ and $\gamma \in Y_{(1^3)},$ i.e. $\gamma_1 \ge 0$ and
$\gamma_1+\gamma_2 \ge 0.$ Then
$$ K_{(1^3)}(\gamma |~q)=q^{\max(\gamma_1,\gamma_1+\gamma_2)}~
\left[\begin{array}{c}\ \min(\gamma_1,\gamma_1+\gamma_2)\\ \ 1\end{array}\right]_{q}. $$
$(iii)$ Let $\eta=(\eta_1,\eta_2,\eta_3)$ be a three component composition of
size $n$, and $\gamma \in Y_{\eta}^{+}$ belongs to the dominant chamber.
Then
\begin{equation}
K_{\Phi(\eta)}(\gamma |~q)=q^{-\gamma_n} \prod_{j=1}^{\eta_1}B_{q}(
\gamma_j + \eta_2 ; \eta_2),
\end{equation}
where for  $n \ge k$
\begin{equation}
B_{q}(n;k)=\sum_{j=0}^{n-k}\left(\begin{array}{c}\j+k-1\\
 \ \j\end{array}\right)q^j = 1/(k-1)! ( \partial / \partial q )^{k-1}
\lbrack (q^{k-1}-q^{n})/(1-q) \rbrack.
\end{equation}

$(iv)$ Let $\eta=(1,\eta_2,\eta_3,\eta_4)$ be a four component composition of
size $n$, $\eta_1=1$, and $\gamma \in Y_{\eta}^{+}$ belongs to the dominant 
chamber. Then
$$ K_{\Phi(\eta)}(\gamma |~q)=q^{- \gamma_n}\sum_{\beta}B_{q}(\beta_1+\eta_3;
\eta_3)\prod_{j=2}^{\eta_2+1}B_{q}(\beta_j+\gamma_j+\eta_3;\eta_3)q^{\beta_j},
$$
where the sum runs over all vectors  $\beta \in \Z^{\eta_2+1}_{\ge 0}$ such 
that $|\beta|=\gamma_1.$

In particular, if $(\gamma_1,\gamma_2,\gamma_3,\gamma_4) \in Y_{(1^4)}^{+}$,
i.e. $\gamma_1 \ge 0,\gamma_2 \ge 0$ and $\gamma_3 \ge 0,$ 
then

$$ K_{\Phi(1^4)}(\gamma_1,\gamma_2,\gamma_3,\gamma_4)=q^{-\gamma_4} \lbrace ~q
\left[\begin{array}{c}\gamma_1+2\\ \ 2\end{array}\right]_{q}
 \left[\begin{array}{c}\gamma_2\\ \ 1\end{array}\right]_{q}
+\sum_{2j \le \gamma_1} q^{2j} 
 (\left[\begin{array}{c}\gamma_1+1-2j\\ \ 1\end{array}\right]_{q})^2 
~ \rbrace 
$$
$$ = q^{-\gamma_n}~
 \sum_{j=1}^{\gamma_1}
q^{j} \left[\begin{array}{c}\gamma_1+1-j\\ \ 1\end{array}\right]_{q} \left[\begin{array}{c}\gamma_2+1-j\\ \ 1\end{array}\right]_{q}.
$$
\end{pr}
Therefore, if $(\gamma_1,\gamma_2,\gamma_3,\gamma_4) \in Y_{(1^4)}^{+},$  then
~~$K_{\Phi(1^4)}(\gamma_1,\gamma_2,\gamma_3,\gamma_4)=
\left(\begin{array}{c}\gamma_1+3\\ \ 3\end{array}\right)+\gamma_2 
\left(\begin{array}{c}\gamma_1+2\\ \ 2\end{array}\right)$.

We remark that $B_{q}(n;l)\vert_{q=1}=
\left(\begin{array}{c}n\\ \l\end{array}\right)$.
\begin{rem} {\rm It is well-known, see e.g. \cite{St3}, \cite{LSt} and 
the literature quoted therein, that on the set of transportation matrices 
of size  $n$ by $m$, the function $|{\cal P}_{\ld\mu}|$ is a continuous 
piecewise polynomial function in $\ld_1,\cdots,\ld_n,\mu_1,\cdots,\mu_m$ of
degree $(n-1)(m-1)$. }
\end{rem}
\begin{quest} It follows from the above Proposition and the formula (5.41)
from Section~5.1, that if $N$ is big enough integer such that $\nu_N:= 
\gamma+N\delta_{\eta_1+\eta_2}$ is a partition, and if we put by definition
 $\ld_N:=N\delta_{\eta_2}$ and $\mu_N:=
N~(\delta_{\eta_1}+(\underbrace{\eta_2, \ldots,\eta_2}_{\eta_1})),$   then
$$ |{\cal P}_{\ld\mu}|=c_{\ld_N,\mu_N}^{\nu_N}.$$
$(\clubsuit)$ Is it true that if N is a big enough integer, then
$${\cal P}_{\ld\mu}(q) \bdoteq c_{\ld_N,~\mu_N}^{\nu_N}(q),$$
where $c_{\ld,\mu}^{\nu}(q)$ denotes the $q$-analog of the $LR$-numbers, 
introduced C.~Carre and B.~Leclerc, and A.~Lascoux, B.~Leclerc and 
J.-Y.~Thibon, see e.g. \cite{LLT} ?
\end{quest}
For the definition of polynomials ${\cal P}_{\ld,\mu}(q)$ see Section~5.4, 
(5.48).

\subsection{Non--vanishing, Degree and Saturation theorems}

It is clear from the very  definition  that $K_{\Phi(\eta)}(\gamma |~q)$ is
a polynomial in $q$ with non-negative integer coefficients.    
For example, if $\eta=(1^n)$, the function $K_{\Phi(1^n)}(\gamma |~q)$
coincides with the  $q$-analog $K_n(\gamma |~q)$ of the Kostant partition 
function $K_n(\gamma),$ see e.g \cite{BV}. 
It is not difficult to see \cite{Kir9} that 

~~~~~~~~~~$K_n(\gamma\ |~q)\ne 0$ if and only if $\gamma \in Y_n,$
where 
$$  Y_n:=\{(\gamma_1,\dots, \gamma_n) \in \Z^n | \sum_{i=1}^{k}\gamma_i \ge 0,
~ 1 \le k\le n,~ \sum_{i=1}^{n} \gamma_i = 0 \}.
$$
Our next goal is to generalize this result to the case of the parabolic 
 $q$-Kostant partition function $K_{\Phi(\eta)}(\gamma |~q)$ corresponding 
to an arbitrary composition $\eta.$

\begin{de} Let $\eta=(\eta_1,\dots,\eta_p)$ be a composition of size
$n$, denote by $Y_{\eta}$ the set of sequences $(\gamma_1,\dots,
\gamma_n) \in \Z^n,$  $|\gamma|=0$,  such that for each integer $k, 
0 \le k \le p-1$, the
 following inequalities are valid:
 $$\sum_{j=1}^{r_k}\gamma_{j} +\sum_{a \in \Omega_k}\gamma_a \ge 0 \quad 
\mbox{for all subsets} \quad \Omega_k
 \subseteq [\eta_k +1,\dots,\eta_{k}+\eta_{k+1}],$$
where $r_k:=\sum_{j \le k} \eta_j$, if~$k \ge 1,$ and~$r_0:=0;$  
 by definition, we put $\eta_0:=0$.
\end{de}  

In particular, we have  
$\gamma_1 \ge 0,\cdots,\gamma_{\eta_1} \ge 0,$ 
  and $\gamma_{r_{p-1}+1} \le 0,\cdots,\gamma_n \le 0.$
\begin{de}
 Denote by $Y_{\eta}^{+}$ the {\it dominant chamber} in the set
$Y_{\eta}$, i.e. the subset of $Y_{\eta}$ consisting of all vectors 
$\gamma=(\gamma_1,\dots,\gamma_n)$ such that $ \gamma_1\ge \dots \ge 
\gamma_{n-1}\ge 0.$
\end{de}
$(\spadesuit)$  We want to stress that if $\gamma \in Y_{\eta}^{+}$, then
$\gamma_{r_{p-1}+1}= \cdots=\gamma_{n-1}=0,$ and $\gamma_n \le 0.$
 
\begin{theorem} ( {\bf Non-vanishing and Degree Theorem for parabolic 
$q$-Kostant partition function} )

 Let $\eta=(\eta_1,\dots,\eta_p)$ be a composition of size 
$n,$ and $\gamma \in \Z^n $ such that $|\gamma|=0$.
 Then

~~~~~$K_{\Phi(\eta)}(\gamma|~q) \ne 0,$ if and only if ~$\gamma \in Y_{\eta}.$ 
 Moreover,
\begin{equation}
 deg K_{\Phi(\eta)}(\gamma|~q)= \sum_{k=1}^{p-1}~(p-k) \ds
(\sum_{j = r_{k-1}+1}^{r_{k}} \gamma_{j}).
\end{equation}
\end{theorem}
Remember that $r_k=\sum_{j \le k}~\eta_j$ if $k \ge 1$, and $r_0=0$.
\begin{ex} Take $\gamma=(2,1,0,-1,0,-1,-1)$ and $\eta=(1,2,2,1,1).$  Using
formula (3.31), let us
compute the degree of the parabolic $q$-Kostant partition function 
$K_{\Phi(\eta)}(\gamma |~q).$   Namely,
 
$\deg K_{\Phi(\eta)}(\gamma |~q)=2+(2+1)+(2+1-1)+(2+1-1-1)=8.$ In fact,

$K_{\Phi(\eta)}(\gamma |~q)=q^3(3,21,52,65,42,13).$
\end{ex}
 If $\gamma \in Y_{\eta}$, so that $K_{\Phi(\eta)}(\gamma|~q) \ne 0,$
we denote by $r(\gamma,\eta)q^{s(\gamma,\eta)}$ its leading term. For example,

$$ r(\gamma,(1^n)) =1, ~~ s(\gamma, (1^n))=\sum_{i=1}^{n-1}(n-i)\gamma_i,
$$
$$ r((3,0,-1,-1,0,-1),(2,3,1))=3,~s((3,0,-1,-1,0,-1),(2,3,1))=4.$$

In general, the number $r(\gamma,\eta)$ can be equal to any positive integer.
As for the number $s(\gamma,\eta),$ it follows from Theorem~3.15  that
$s(\gamma,\eta)= (\gamma,\delta_{\Phi(\eta)}),$
where $\delta_{\Phi(\eta)}$ denotes the vector with components 
$(\delta_{\Phi(\eta)})_{i}= p-k $ if $ r_{k-1} < i \le r_k, ~k=1,\cdots,p$.

  Moreover, the numbers $s(\gamma,\eta)$  satisfy the  so-called 
{\it saturation} property.
\begin{cor} ( {\bf Saturation theorem for parabolic Kostant partition functions} )

For any positive number $N$ we have
$$ s(N\gamma,\eta)=Ns(\gamma,\eta).$$
\end{cor}
\begin{con} ( {\bf Unimodality conjecture for parabolic Kostant partition
functions} )

Let $\eta$ be a composition of size $n,$ and 
$\gamma \in \Z^n $ such that $|\gamma|=0$. Then, $K_{\Phi(\eta)}(\gamma|~q)$
is a unimodal polynomial in the variable $q$.
\end{con}

\subsection{Rationality and polynomiality theorems}

\begin{theorem} ( {\bf Rationality theorem for parabolic Kostant partition 
function, I} )

Let $\eta$ be a composition and $\gamma \in Y_{\eta}.$  Then 
$$ \sum_{n \ge 0} K_{\Phi(\eta)}(n\gamma |~q)t^n = 
P_{\eta\gamma}(q,t)/Q_{\eta\gamma}(q,t),$$
where $P_{\eta\gamma}(q,t)$ and $Q_{\eta\gamma}(q,t)$ are mutually prime 
polynomials in $q$ and $t$ with integer coefficients, $P_{\eta\gamma}(0,0)=1.$

Moreover,

$(\clubsuit)$ the denominator $Q_{\eta\gamma}$ has the following 
form:
$$ Q_{\eta\gamma}(q,t) = \prod_{j \in J} (1-q^{j}~t),$$
where  $J:=J_{\eta\gamma}$ is a finite set of non--negative integer numbers, 
 not necessarily distinct;

$(\clubsuit\clubsuit)$  $P_{\eta\gamma}(1,t)=(1-t)^{t(\eta,\gamma)}
~P_{\eta\gamma}(t),
~~P_{\eta\gamma}(1) \ne 0,$ ~where $t(\eta,\gamma) \in \Z_{\ge 0},$ and 
$P_{\eta\gamma}(t)$ is a polynomial with non--negative integer coefficients.
\end{theorem}
$(\maltese)$ We {\bf expect} that if $\gamma_1$ and $\gamma_2$ belong to the
set $Y_{\eta},$ and $\gamma_1 \ge \gamma_2,$ i.e. 
$ \sum_{ j \le k}\gamma_{1,j} \ge \sum_{j \le k}\gamma_{2,j},$ $\forall 
k \ge 1,$  then
$$P_{\eta\gamma_2}(t) - P_{\eta\gamma_1}(t) \ge 0.$$
In other words, the latter difference is a polynomial with non--negative
coefficients.  

\begin{cor} ( {\bf Polynomiality theorem for parabolic Kostant partition
 function} )    

Let $\eta$ be a composition and $\gamma \in Y_{\eta}$. 
 There exists a polynomial ${\cal K}_{\eta\gamma}(t)$ with rational
 coefficients such that for any integer number $N \ge 1$,
${\cal K}_{\eta\gamma}(N)=K_{\Phi(\eta)}(N \gamma).$
\end{cor}
\begin{con} 
The polynomials ${\cal K}_{\eta\gamma}(t)$ have  {\bf non--negative} rational 
coefficients.  
\end{con} 
\begin{theorem} ( {\bf Rationality theorem for parabolic Kostant partition
function, II} )

Let $\eta$ be a composition and $\gamma_1,\cdots,\gamma_k  \in Y_{\eta}$. 
Then the generating function
$$ \sum_{(N_1,\cdots,N_k) \in \Z_{\ge 0}^{k}} K_{\Phi(\eta)}(N_1 \gamma_1+
\cdots+N_k \gamma_k |~q)~x_{1}^{N_1} \cdots x_{k}^{N_k}$$
is a rational function in $q$ and the variables $X_{k}=(x_1, \cdots, x_k)$ of
the form $P(q,X_k)/Q(q,X_k),$ where 
$P:=P_{\gamma_1,\cdots,\gamma_k,\eta}(q,X_k)$ and $Q(q,X_k)
:=Q_{\gamma_1,\cdots,\gamma_k,\eta}(q,X_k)$ are mutually prime polynomials in
$q$ and $X_k$ with integer coefficients, $P(0,0)=1.$ 

$(\clubsuit)$ Moreover, the denominator $Q(q,X_k)$ has the following structure:
$$ Q(q,X_k)=\prod_{ \emptyset \ne W \subset \{ 1,\cdots, k \}}~  
\prod_{a_W \in J_W} (1-q^{a_W}~x_{W}),$$
where $x_W:= \prod_{i \in W}x_i,$ and for each non--empty subset $W \subset  
\{ 1,\cdots, k \},$ $ J_W$ denotes a certain set, depending on $W$ and
$\gamma_1,\cdots,\gamma_k,$ of non--negative integers, not necessarily
distinct.  
\end{theorem}
$(\maltese)$ We {\bf expect} that if $W= \{b \}, 1 \le b \le k,$ then $J_W=
J_{ \gamma_{b},\eta}.$
\begin{cor} ( {\bf Piecewise polynomiality theorem for parabolic Kostant 
partition function} )

Let $\eta$ be a composition and $\gamma_1,\cdots,\gamma_k \in Y_{\eta}.$
There exists a piecewise polynomial function ${\cal K}(t_1,\cdots,t_k):=
{\cal K}_{\gamma_1,\cdots,\gamma_k}(t_1,\cdots,t_k)$ with rational
 coefficients such that for any non--negative integer
numbers $N_1,\cdots,N_k$, ~~${\cal K}(N_1,\cdots,N_k)=K_{\Phi(\eta)}
(N_1 \gamma_1+\cdots+N_k \gamma_k).$
\end{cor}
$(\maltese)$  We {\bf expect} that the restriction of the function
${\cal K}(t_1,\cdots,t_k)$ on ``the dominant chamber'' ${\cal N}_k:=
\{ (N_1 \ge N_2 \ge \cdots \ge N_k) \in \Z_{\ge 0}^{k} \}$ is a polynomial 
with {\bf non--negative} rational coefficients. 

\begin{ex} {\rm Take $\gamma_1= (2,1,0,-1,-1,-1),$ $\gamma_2=(1,1,-1,-1)$ and
$\eta=(1^5)$.  Then $Q_{\gamma_1,\gamma_2}(1,x,y)=(1-x)^7(1-y)^4,$  and 
$$P_{\gamma_1,\gamma_2}(1,x,y)=(1,26,71,26)+(1,-57,-223,-93)y+(0,33,224,115)y^2-(0,8,66,50)y^3.$$
Therefore, in our example the function $(n,m) \rightarrow K_{\Phi(1^5)}(n(2,1,-1,-1,-1)+m(1,1,-1,-1))$ is a polynomial one on the whole set 
$\{(n,m) \in \Z_{\ge 0}^2 \}$.  }
\end{ex}
$(\maltese\maltese)$ We {\bf expect} that in fact the function
$(n_1,\cdots,n_k) 
\rightarrow K_{\Phi(\eta)}(n_1\gamma_1+\cdots+n_k\gamma_k)$ is a 
{\bf polynomial} one on the whole set 
$\{(n_1,\cdots,n_k) \in \Z_{\ge 0}^k \},$ cf {\it mixed lattice point 
enumerator theorem} by P.McMullen \cite{McM}.

\subsection{Parabolic Kostant partition function $K_{\Phi(\eta)}(\gamma)$
as function of $\gamma$}

 In this Section we state a few theorems, problems  and one conjecture 
about behavior of the parabolic Kostant partition function 
 $K_{\Phi(\eta)}(\gamma),$ considered as a {\it function} of $\gamma,$ on 
the set $Y_{\eta}$.
\begin{theorem} ( {\bf Polynomial expression for  the restriction of the 
parabolic Kostant partition function $K_{\Phi(\eta)}(\gamma )$ on the 
dominant chamber $Y_{\eta}^{+}$} )

Let $\eta=(\eta_1,\cdots,\eta_p)$, $p \ge 3$, $\eta_p \ne 0$, be a
composition, consider 
vector $l=(l_1,l_2,\cdots,l_{r_{p-2}})$, where
$l_i=\sum_{j=k+1}^{p-1} \eta_j $ if $r_{k-1} < i \le r_{k}, 1 \le k \le p-2.$ 
Let $\widehat\eta=(\eta_1,\cdots,\eta_{p-2})$.  If $\gamma \in Y_{\eta}^{+},$
 then
\begin{equation}
K_{\Phi(\eta)}(\gamma )= \sum_{\beta}K_{\Phi(\widehat\eta)}
(\beta_1-l_1,\cdots,\beta_{r_{p-2}}-l_{r_{p-2}}) \prod_{j=1}^{r_{p-2}} 
\left(\begin{array}{c}\gamma_j+l_j\\ \ \beta_j\end{array}\right),
\end{equation}
where the sum runs over $\beta \in \Z_{ \ge 0}^{r_{p-2}}$ such that $|\beta|=
|l|=\sum_{1 \le i < j \le p-1} \eta_i \eta_j$.
\end{theorem}

\begin{cor} Being restricted on the dominant chamber $Y_{\eta}^{+},$ the 
function $F_{\eta}(\gamma):=K_{\Phi(\eta)}(\gamma)$
is a {\bf polynomial} in $\gamma_1,\cdots,\gamma_{r_{p-2}}$  of  
degree $|l|= \sum_{1 \le i <j \le p} \eta_i \eta_j - 
\eta_p(n-\eta_p)$  with rational 
coefficients. 
\end{cor}
\begin{theorem} ( {\bf Piecewise polynomiality theorem for function  
$\gamma \longrightarrow K_{\Phi(\eta)}(\gamma)$} )

 On the set $Y_{\eta}$ the function 
$\gamma \longrightarrow F_{\eta}(\gamma):=
K_{\Phi(\eta)}(\gamma)$  is a continuous piecewise polynomial 
function of degree $\sum_{1 \le i < j \le p}\eta_i \eta_j -n +1$.
\end{theorem}
We see that if $\eta_p > 1$, then the dominant chamber $Y_{\eta}^{+}$ is
strictly contained in some maximal polynomiality domains of the function 
$F_{\eta}$.
\begin{prb} Count the number and describe a structure of the polynomiality
domains of the function $F_{\eta}$.
\end{prb}
\begin{con} Restriction of the function $|l|!~F_{\eta}$ on the dominant chamber
$Y_{\eta}^{+}$, denoted by $F_{\eta}^{+}$,  is a polynomial in 
$\gamma_1,\cdots,\gamma_{r_{p-2}}$  with 
{\bf non--negative} integer coefficients.
\end{con}
\begin{prb} Find a combinatorial interpretation of the coefficients of the
polynomial $F_{\eta}^{+}.$
\end{prb}

\subsection{Reconstruction theorem}

The leading term $|l|!~G_{\eta}(\gamma)$  of the polynomial 
$F_{\eta}^{+}(\gamma),$
 i.e. the degree $|l|$ homogeneous part of $F_{\eta}^{+}(\gamma),$ admits the
following description.
\begin{de} For any composition $\eta=(\eta_1,\cdots,\eta_p),$ such that 
$\eta_p > 0$ and $p \ge 3,$ define the operator
$$  {\cal D}_{\eta}=\prod_{1 \le i \le \eta_1 < j \le r_{p-2}} 
(\partial/\partial\gamma_i -\partial/\partial\gamma_j), $$
acting on the quotient ring of the ring of polynomials 
${\Q}~[\gamma_1,\cdots,\gamma_n]$ by the ideal generated by the sum
$\gamma_1+\cdots+\gamma_n.$
\end{de}
Let $\gamma=(\gamma_1,\cdots,\gamma_n) \in \Z^n,$ $|\gamma|=0.$
\begin{theorem} ( {\bf Characterization of polynomials $G_{\eta}(\gamma)$} )

The polynomials $G_{\eta}(\gamma)$ are uniquely determined by
the following properties

$(i)$ $G_{\eta}(\gamma)$ is a homogeneous polynomial of degree 
 $|l|=\sum_{1 \le i < j \le p-1} \eta_i \eta_j,$

$(ii)$ ${\cal D}_{\eta}G_{\eta}(\gamma)=
\prod_{j=1}^{\eta_1}( \gamma_j^{\eta_{p-1}}/ \eta_{p-1}! ) ~~G_{(\eta_2,\cdots,\eta_p)}(\gamma_{\eta_1+1},\cdots,\gamma_n),$

$(iii)$ $G_{\eta_1,\eta_2}(\gamma)=1.$
\end{theorem}
\begin{theorem} ( {\bf Reconstruction Theorem } )

Let $G_{\eta}(\gamma)=\sum_{\beta} b_{\eta}(\beta) \prod_{j=1}^{r_{p-2}} 
\gamma_j^{\beta_j}/\beta_j !$,
summed over $\beta \in \Z_{\ge 0}^{r_{p-2}}$ such that $|\beta|=|l|.$ Then
$$ F_{\eta}(\gamma)=\sum_{\beta} b_{\eta}(\beta)~\prod_{j=1}^{r_{p-2}} 
\left(\begin{array}{c}\gamma_j+l_j\\ \ \beta_j\end{array}\right).
$$
\end{theorem}
\begin{cor} Let $l$ be the vector defined in Theorem~3.23, then
$$b_{\eta}(\beta)=K_{\Phi(\eta)}(\beta-l).$$ 
In particular,
$G_{\eta}(\gamma)$ is a polynomial with {\bf non--negative} rational 
coefficients.
\end{cor}

Finally, we state  a result which is a refinement of 
Proposition~3.7, and gives partly a $q$-analog of the recurrence relation 
(3.32). 
\begin{theorem} ( {\bf A $q$-analog of Theorem~3.23} )   

Let $\eta=(\eta_1,\cdots,\eta_p)$, $p \ge 3$, $\eta_p \ne 0$,
 be a composition. Define  $\widehat\gamma=(\gamma_1,\cdots,\gamma_{r_{p-3}}
\underbrace{0,\ldots,0}_{\eta_{p-2}})$  and $\widehat\eta=(\eta_1,\cdots,
\eta_{r_{p-2}})$. If $\gamma \in Y_{\eta}^{+}$, then  \\

$K_{\Phi(\eta)}(\gamma |~q)=q^{-\gamma_n}\sum_{\beta}K_{\Phi(\widehat\eta)}(
\widehat\gamma-\beta  |~q)\prod_{j=1}^{r_{p-3}}B_{q}(\beta_j+\eta_{p-1};
\eta_{p-1}) \prod_{j=r_{p-3}+1}^{r_{p-2}}B_{q}(\gamma_j+\beta_j + \eta_{p-1};
\eta_{p-1}),$

where the sum runs over vectors $\beta \in \Z_{\ge 0}^{r_{p-2}}$ such that
$|\beta|=\sum_{j=1}^{r_{p-3}}\gamma_j$, and polynomials $B_{q}(n;k)$ are 
defined in Proposition~3.8, formula (3.30).
\end{theorem}
\begin{rem} {\rm The ``classical'' case $\eta=(1^n)$ and $q=1,$ which 
corresponds
to the Kostant partition function $K_{n}(\gamma)$, has been studied by
F.~Berezin and I.M~Gelfand \cite{BeG}, B.~Kostant \cite{K}, 
B.V.~Lidskii \cite{Lid1},\cite{Lid2}, D.~Peterson, A.N.~K. \cite{Kir8},
 \cite{Kir9}, A.~Postnikov and R.~Stanley \cite{PS}, W.~Baldoni-Silva 
and M.~Vergne \cite{BV}, S.~Billey, V.~Guillemin and 
E.~Rassart \cite{BGR}, J.~De Loera and B.~Sturmfels \cite{LSt}, ... . 
In particular, if $\eta=(1^n)$ and $q=1$, Theorem~3.19 has been proved by B.V.
Lidskii \cite{Lid1} in 1984, and by D.~Peterson (unpublished). The case of 
arbitrary $\eta$ and $q$ has been studied by the author (unpublished,
but see \cite{Kir9}). The case of generalized Kostant partition functions and 
$q=1$ has been studied by A.~Postnikov and R.~Stanley (unpublished, but see
\cite{PS}).       }
\end{rem} 

\section{Parabolic Kostka polynomials:  \\
 Definition and basic properties  }

\begin{de} (\cite{Kir9},\cite{KS}) Let $\lambda$ be a partition and $\mu$ 
and $\eta$ be compositions
such that $|\ld|=|\mu|,$ ~$|\eta| = n$  and $ll(\mu) \le n.$ 
 Define the {\it parabolic Kostka
polynomial}
 $K_{\lambda\mu\eta}(q)$ as follows:
\begin{equation}
K_{\lambda\mu\eta}(q):=\sum_{w\in \Sigma_n}(-1)^{l(w)}
K_{\Phi(\eta)}(w(\lambda+\delta)-\mu-\delta |~q),   \label{2.4}
\end{equation}
where $\delta
:=\delta_n=(\hbox{$n-1,n-2,$}\ldots ,1,0)$.
\end{de}

If a composition $\mu$ is compatible with $\eta$ and
corresponds to
the sequence of partitions (possibly with zeros at the end) 
${\bmu}=(\mu^{(1)},\mu^{(2)},\cdots,\mu^{(r)})$,
we will denote
the parabolic Kostka polynomial $K_{\lambda\mu\eta}(q)$ by
$K_{\lambda,~{\sbmu}}(q)$
or $K_{\lambda,(\mu^{(1)},\mu^{(2)},\cdots,\mu^{(r)})}(q).$
If a sequence of partitions $\bmu=(\mu^{(1)},\cdots,\mu^{(r)})$ consists of 
only {\it rectangular} shape partitions
 $\mu^{(a)} = (\mu_a^{\eta_a}):=R_a,  ~~1\le a\le r$, we will write 
$R=(R_1,R_2,\dots,R_r)$ instead of $\bmu$, and $K_{\lambda,R}(q)$  instead of
$K_{\lambda,~{\sbmu}}(q).$

 Let us elucidate Definition~4.1 by a simple, but interesting example.
\begin{ex} Take $\ld=(6,2,2,2)$, $\mu=(2^6)$ and $\eta=(2^3)$. There are 4
contributions to the RHS(4,33), namely,
$$ K_{\ld\mu\eta}(q)=K_{\Phi(\eta)}(\gamma_1 |~q)-K_{\Phi(\eta)}(\gamma_2 |~q)-
K_{\Phi(\eta)}(\gamma_3 |~q)+K_{\Phi(\eta)}(\gamma_4 |~q),$$
where $\gamma_1=\ld-\mu=(4,0,0,0,-2,-2)$, $\gamma_2=(4,0,0,0,-3,-1)$,
$\gamma_3=(4,0,-1,1,-2,-2)$ and $\gamma_4=(4,0,-1,1,-3,-1)$. It is not
difficult to see that $K_{\Phi(\eta)}(\gamma_1 |~q)=q^4(1,4,10,12,9)$,
$K_{\Phi(\eta)}(\gamma_2 |~q)=q^4(1,4,7,10,8)$, $K_{\Phi(\eta)}(\gamma_3 |~q)=
q^5(2,7,10,7)$ and $K_{\Phi(\eta)}(\gamma_4 |~q)=q^5(2,5,8,6)$. Hence,
$K_{\ld\mu\eta}(q)=q^6,$ and $degK_{\ld\mu\eta}(q)=6 < 
degK_{\Phi(\eta)}(\ld-\mu |~q)=8$.
\end{ex}
\begin{rem} {\rm Using in Definition~4.1 the $q$-analog 
$K_{\mit\Sigma}(\gamma |~q)$ of the 
generalized Kostant partition function, see Section~3.1, Remark~3.4, one 
can define the ``generalized'' Kostka polynomials  $K_{\ld\mu\mit\Sigma}(\gamma |~q)$. They
form  an interesting family of polynomials to study.} 
\end{rem}
\begin{theorem}\label{t3.3} {\rm (\cite{SW})} Let $\lambda$ be a partition,
and $\mu$ be a composition compatible with $\eta.$  Then
\begin{equation}
K_{\lambda,~{\sbmu}}(1):=K_{\lambda,(\mu^{(1)},\mu^{(2)},\cdots,\mu^{(r)})}(1)=
{\rm Mult}[V_{\lambda}:\otimes_{i=1}^rV_{\mu^{(i)}}], \label{3.6}
\end{equation}
i.e. $K_{\lambda,(\mu^{(1)},\mu^{(2)},\cdots,\mu^{(r)})}(1)$ is equal to the
multiplicity of the irreducible
highest weight $\lambda$ $\g l(n)$--module $V_{\lambda}$ in the tensor
product of irreducible highest  weight $\mu^{(i)}$
representations $V_{\mu^{(i)}}$,  \\
$1\le i\le r$, of the Lie algebra $\g l(n)$.
\end{theorem}
In the case when all partitions $\mu^{(i)}$ have rectangular shapes,
Theorem~\ref{t3.3} has been proved in \cite{Kir}.
\vskip 0.3cm
\begin{rem} { \rm  We {\bf expect} that $K_{\ld\mu\eta}(1) \ge 0$ for any 
partition $\ld$
and compositions $\mu$ and $\eta.$ It seems a challenge problem 
to find a combinatorial and/or representation-theoretic interpretations of the
numbers $K_{\ld\mu\eta}(1)$ and $K_{\ld\mu\eta}(-1)$ for  general $\ld,\mu$
and $\eta$. In particular,  \\

$(\clubsuit)$  When does the number $K_{\ld\mu\eta}(1)$ equal to $1$  ? }
\end{rem}
\begin{exs} In these examples we will use notation $P_{\ld\mu\eta}(q,t),$ 
$Q_{\ld\mu\eta}(q,t)$ and $J_{\ld\mu\eta}(q),$ which will be explained  in 
Theorem~4.14.

$(i)$  Take $\ld=(3,2,1)$,  $\mu=(2,2,2)$ and $\eta=(1^3)$. Then 
$K_{\ld\mu\eta}(q)=K_{\ld\mu}(q)=q+q^2,$ and
$$ \sum_{n \ge 0} K_{n\ld,n\mu,\eta}(q)~t^n = (1-qt)^{-1}(1-q^2t)^{-1}.
$$
$(ii)$ Take the same $\ld$, but $\mu=(0,2,2,2)$ and $\eta=(1^4)$. Then

$K_{\ld\mu\eta}(q)=q^3(-1,-1,0,1,2,1)$, $K_{2\ld,2\mu,\eta}(q)=
q^5(1,0,-2,-4,-4,-1,0,3,3,4,2,1)$. 

Moreover,
$$ P_{\ld\mu\eta}(q,t)=1-q^2(1,3,2,1)t+\cdots+q^{33}(-1,1,1,0,-1,-1)~t^7,$$
$$ Q_{\ld\mu\eta}(q,t)=(1-q^3t)(1-q^7t) \prod_{j=2}^{8}(1-q^j~t),
~~J_{\ld\mu\eta}(q)=q^2(1,2,1,1,1,2,1),
$$
see Theorem~4.14, $(\spadesuit),$ for the definition of polynomials 
$J_{\ld\mu\eta}(q).$  \\
 
$(iii)$ Take again $\ld=(3,2,1)$, but $\mu=(0,2,0,2,2)$ and $\eta=(1,2,1,1)$.
Then

$K_{\ld\mu\eta}(q)=q^3(1,0,-4-3,2,4,2)$, $K_{2\ld,2\mu,\eta}(q)=
q^7(3,5,6,-3,-13,-17,-11,3,9,12,6,3)$. Moreover,
$$ P_{\ld\mu\eta}(q,t)=1-q^3(1,2,6,5,0,-2)t+\cdots+q^{70}(1,-1,-2,1,2)~t^{12},$$
$$ Q_{\ld\mu\eta}(q,t)= \prod_{j=3}^{9} (1-q^j~t)^2,
~~J_{\ld\mu\eta}(q)=q^3(2,2,2,2,2,2,2).
$$
$(iv)$ Take the same $\ld=(3,2,1),$ but $\mu=(0,2,0,2,0,2)$ and 
$\eta=(1,2,2,1).$  Then

$K_{\ld\mu\eta}(q)=q^4(1,2,-8,-6,8,5),$  $K_{2\ld,2\mu,\eta}(q)=-q^7-2q^8
+\cdots+22q^{17}+12q^{18}.$  Moreover,

$P_{\ld\mu\eta}(q,t)=1-q^3(2,3,2,11,10,-4,-2)~t+\cdots+
q^{132}(1,-1,-2,1,2)~t^{22},$ and 

$J_{\ld\mu\eta}(q)=q^3(2,4,4,3,4,4,3).$  In other words,
$$ Q_{\ld\mu\eta}(q,t)=(1-q^3t)^2(1-q^4t)^4(1-q^5t)^4(1-q^6t)^3(1-q^7t)^4(1-q^8t)^4(1-q^9t)^3. $$ 
\end{exs}
$(\clubsuit)$  We would like to remark that the reasons for the equality 
below are  elusive.
$$q^{62}~P_{(3,2,1),(0,2,0,2,2),(1,2,1,1)}(q,t) \vert_{t^{12}} = 
P_{(3,2,1),(0,2,0,2,0,2),(1,2,2,1)}(q,t) \vert_{t^{22}}.$$

$(\spadesuit)$ These examples show that for general $\ld,$ $\mu$ and $\eta$, 
the polynomials 
$K_{\ld\mu\eta}(q)$ may have negative coefficients, the numbers 
$a(\ld,\mu \Vert \eta)$ may be negative  and may not be a homogeneous function 
in $n,$ and those
$b(\ld,\mu \Vert \eta)$ may not satisfy the (generalized) 
Fulton conjecture.     \\

Our nearest goal is to describe several cases when the polynomials  
$K_{\ld\mu\eta}(q)$ have
only {\it non--negative} coefficients. However, we want to point out that 
there are many other cases when
the all coefficients of a parabolic Kostka polynomial  are non--negative.
 
\begin{ex} Take $\ld=(6,3,2,1)$, $\mu=(2,1,2,1,2,1,2,1)$ and $\eta=(2^4)$.
Then $K_{\ld\mu\eta}(q)=q^{11}(4,18,24,14,4)$. It is interesting to compare
the polynomial $K_{\ld\mu\eta}(q)$ with the $q$-analog of the $LR$-numbers 
$c_{\mu^{(1)},\cdots,\mu^{(r)}}^{\ld}(q)$ introduced by C.~Carre, A.~Lascoux, 
B.~Leclerc and J.-Y.~Thibon, see e.g. \cite{LLT}. Namely, one can show that 
$c_{(2,1),(2,1),(2,1),(2,1)}^{(6,3,2,1)}(q) = q^8(2,7,12,15,14,9,4,1)$.
\end{ex}
\begin{pr} Let $\ld$ be a partition and $\bmu=(\mu^{(1)},\cdots,\mu^{(r)})$
be a sequence of partitions. \\
 If inequalities $ll(\mu^{(i)})\ge\ l(\ld)$  holds for all i, then
\begin{equation}
K_{\ld,~{\sbmu}}(q)\bdoteq Mult[V_{\lambda}:\otimes_{i=1}^rV_{\mu^{(i)}}].
 \label{4.22}
\end{equation}
\end{pr}
\begin{pr} Let $\ld$ be a partition and $\bmu=(\mu^{(1)},\mu^{(2)})$ be
a {\bf dominant} sequence of partitions. Then
\begin{equation}
K_{\ld,~{\sbmu}}(q)\bdoteq c_{\mu^{(1)},\mu^{(2)}}^{(\ld)}.  \label{4.23}
\end{equation}
\end{pr}
See Introduction, Section~1.1, for the explanation of the meaning of the
symbol ``$\bdoteq$''. \\

\hskip -0.6cm{\bf Positivity Theorem} {\rm (\cite{Kir1},\cite{KSS})}
{\it Let $\ld$ be a partition, and
$\bmu =(\mu^{(1)},R_2,\cdots,R_r)$ be a sequence of (proper) partitions
such that

$(a)$ $(R_2,\cdots,R_r)$ is a dominant sequence of rectangular shape 
partitions,

$(b)$ {\bf either} $ll(\mu^{(1)}) \ge l(\ld)$, 

{\bf or} $\ld\supset\mu^{(1)}$ and the complement $\ld\setminus\mu^{(1)}$ is a
disjoint union of partitions $\ld^{(1)},\ld^{(2)},\cdots,\ld^{(p)}$.

Then the parabolic Kostka polynomial $K_{\ld,(\mu^{(1)},R_2,\cdots,R_r)}(q)$
has non-negative integer coefficients.}
\vskip 0.3cm
\begin{con} ( {\bf Positivity conjecture for parabolic Kostka polynomials,
 cf  \cite{Kir5},\cite{KS} } )

Let $\ld$ be a partition and $\bmu = (\mu^{(1)},\mu^{(2)},\cdots,\mu^{(r)})$ 
be a sequence of  (proper) partitions such that 
$(\mu^{(2)},\cdots,\mu^{(r)})$ is a {\bf dominant} sequence of  partitions. 
Assume that 

{\bf either}  $\ld\supset\mu^{(1)}$  and the complement 
$\ld\setminus\mu^{(1)}$ is a
disjoint union of partitions $\ld^{(1)},\ld^{(2)},\cdots,\ld^{(p)},$
~~~~~{\bf or} $ll(\mu^{(1)}) \ge l(\ld)$.

Then
$$K_{\ld,~{\sbmu}}(q) \in \N~[q].
$$
\end{con}
$(\maltese)$ In particular, we {\bf expect}  \cite{Kir5}, \cite{KS}  that 
if $\ld$ and $\mu$ are  partitions and $\eta$ is a composition, then
$$ K_{\ld\mu\eta}(q) \in \N~[q].
$$
\begin{rem}  {\rm According to (\ref{3.6}) and Conjecture 4.10, if $\mu$
is a (proper) partition, then the
parabolic Kostka polynomials $K_{\ld\mu\eta}(q)$ may be
considered as a $q$--analog of the tensor product multiplicities.
Another $q$--analog of the tensor product multiplicities has been
introduced by C.~Carre and B.~Leclerc \cite{CB}, and A.~Lascoux, B.~Leclerc 
and J.-Y.~Thibon \cite{LLT}. Formulas
(\ref{4.22}) and (\ref{4.23}) show that in general these two $q$--analogs are
different. However, it was conjectured in  \cite{Kir5}, Conjecture~6.5 and in 
 \cite{KS}, Conjecture~5, that, in fact,
these two $q$--analogs {\it coincide} in the case when a partition $\mu$
and a composition $\eta$ correspond to a {\it dominant} sequence of
rectangular shape partitions.}
\end{rem}

\hskip -0.6cm{\bf Duality Theorem} (\cite{Kir6,KS}) {\it Let $\ld$
be a partition, and $R$ be a dominant sequence of rectangular
shape partitions, $R=((\mu_a^{\eta_a}))_{a=1}^r$. Denote by $R'$
a dominant rearrangement of the sequence of rectangular shape
partitions $((\eta_a^{\mu_a}))_{a=1}^r$ obtained by transposing
each of the rectangular in $R$. Then
\begin{equation}
K_{\ld'R'}(q)=q^{n(R)}K_{\ld R}(q^{-1}), \label{3.12}
\end{equation}
where $n(R)=\ds\sum_{1\le a<b\le p}\min (\mu_a,\mu_b)\min (\eta_a,\eta_b)$.}
\vskip 0.3cm
Note that the left hand side of (\ref{3.12}) is computed in $\g
l(m)$, where $m=\sum\mu_a$ is the total number of columns in the
rectangles of $R$, whereas the right hand side of (\ref{3.12}) is
computed in $\g l(n)$, where $n=\sum\eta_a$ is the total number of
rows in the rectangles of $R$. 
\begin{cor} We have

$(i)$  $a(\ld,R )= n(R)-c(\ld',R'),$ 

$(ii)$ $b(\ld,R) = d(\ld',R')$ 
\end{cor}
\begin{con} Let $\ld$ and $\mu$ be partitions, and $\eta_1$ and $\eta_2$
be compositions such that $\eta_2$ is a subdivision of $\eta_1.$ Then
$$ K_{\ld\mu\eta_1}(q) \le K_{\ld\mu\eta_2}(q).$$
\end{con}
We remark that Conjecture~4.13 is false if one assumes only that $\eta_1 \ge
\eta_2$ with respect to the dominance partial ordering on the set of 
compositions, see Example~3.6. 
\begin{theorem} ( {\bf Rationality theorem for parabolic Kostka polynomials, I} )
 
The formal power series 
$$ \sum_{n \ge 0} K_{n\ld,n\mu,\eta}(q)t^{n} $$
is a rational function in  $q$ and $t$  of the form
$$ P_{\ld\mu\eta}(q,t)/ Q_{\ld\mu\eta}(q,t),$$
where $P_{\ld\mu\eta}(q,t)$ and $Q_{\ld\mu\eta}(q,t)$ are mutually prime
polynomials in $q$ and $t$ with integer coefficients and $P_{\ld\mu\eta}(0,0)
=1$.

Moreover,

$(\clubsuit)$ the denominator $Q_{\ld\mu\eta}$ has the following 
form:
$$ Q_{\ld\mu\eta}(q,t) = \prod_{j \in J} (1-q^{j}~t),$$
where $J:=J_{\ld\mu\eta}$
is a finite set of non--negative integer numbers, not necessarily distinct;

$(\clubsuit\clubsuit)$  $P_{\ld\mu\eta}(1,t)=(1-t)^{t(\ld,\mu,\eta)}
~P_{\ld\mu\eta}(t),$
where $t(\ld,\mu,\eta) \in \Z_{\ge 0},$ ~$P_{\ld\mu\eta}(1) \ne 0,$ and
$P_{\ld\mu\eta}(t)$ is a polynomial with non--negative integer coefficients. 
\end{theorem}
 
$(\spadesuit)$  It is convenient to depict the set $J_{\ld\mu\eta}$ in the 
polynomial  $J_{\ld\mu\eta}(q)=\sum_{j \in J_{\ld\mu\eta}} q^j.$

$(\maltese)$ We {\bf expect} that if $\mu_1$ and $\mu_2$ are partitions such
that $\mu_1 \ge \mu_2$ with respect to the dominance partial ordering, see e.g.
Section~2.1, then
$$P_{\ld,\mu_{2},\eta}(t)-P_{\ld,\mu_{1},\eta}(t) \ge 0,$$
i.e. the latter difference is a polynomial with non--negative coefficients.

\begin{cor} ( {\bf Polynomiality theorem for parabolic Kostka numbers} )

Let $\ld$ be a partition, and $\mu$ and $\eta$ be compositions
such that $\ld-\mu \in Y_{\eta}.$ 
There exists a polynomial ${\cal K}_{\eta\mu\eta}(t)$ with rational 
coefficients such that 

$(\clubsuit)$  for any integer number $N \ge 1$, ${\cal K}_{\eta\mu\eta}(N)=
K_{N\ld,N\mu,\eta}(1).$
\end{cor}
\begin{con} 
If $\mu$ is a partition and $\eta$ is a composition, then the polynomial
${\cal K}_{\eta\mu\eta}(t)$ has {\bf non--negative} rational coefficients.
\end{con}
 Theorem~4.14 is a corollary of the corresponding theorem for parabolic
Kostant's partition function ( Theorem~3.17 ) . In  Section~6,  Rationality
 Conjecture, we state  a few conjectures  about the structure of the numerator
 $P_{\ld\mu\eta}(q,t).$  
\begin{theorem} ( {\bf Rationality theorem for parabolic Kostka polynomials, II})

Let $ \blambda=(\ld^{(1)},\cdots,\ld^{(k)})$ be a sequence of partitions, 
$\bmu=(\mu^{(1)},\cdots,\mu^{(k)})$ be a sequence of compositions 
and $\eta$ be a composition such that $|\ld^{(j)}|=|\mu^{(j)}|$ and 
$ll(\mu^{(j)}) \le |\eta|$ for all $ 1 \le j \le k.$  Then the generating
function
$$ \sum_{(n_1,\cdots,n_k) \in \Z_{ \ge 0}^{k}} K_{n_1\ld^{(1)}+\cdots+n_k\ld^{(k)},n_1\mu^{(1)}+\cdots+n_k\mu^{(k)},\eta}(q)~x_1^{n_1} \cdots x_k^{n_k}
$$
is a rational function in $q$ and the variables $X_k:=(x_1,\cdots,x_k)$ of the
form $P(q,X_k)/Q(q,X_k),$
where $P(q,X_k):=P_{~{\sblambda},~{\sbmu},\eta}(q,X_k)$ and $Q(q,X_k):
=Q_{~{\sblambda},~{\sbmu},\eta}(q,X_k)$ are mutually prime polynomials in $q$ 
and $X_k$ with integer coefficients, $P(0,0)=1.$

$(\clubsuit)$ Moreover, the denominator $Q(q,X_k)$ has the following structure:
$$ Q(q,X_k)=\prod_{ \emptyset \ne W \subset \{ 1,\cdots, k \}}~  
\prod_{a_W \in J_W} (1-q^{a_W}~x_{W}),$$
where $x_W:= \prod_{i \in W}x_i,$ and for each non--empty subset $W \subset  
\{ 1,\cdots, k \},$  $J_W$ denotes a certain set, depending on $W$ and
$\blambda,\bmu,\eta,$ of non--negative integer numbers, not necessarily
distinct. 
\end{theorem}
$(\maltese)$ We {\bf expect} that in general, all the sets $J_{W},$  
$\emptyset \ne W \subset \{1,\cdots,k \},$ are non trivial, i.e. each contain
 at least one positive element. 
\begin{cor} ( {\bf Piecewise polynomiality theorem for parabolic Kostka 
numbers} )

Let $ \blambda=(\ld^{(1)},\cdots,\ld^{(k)})$ be a sequence of partitions, 
$\bmu=(\mu^{(1)},\cdots,\mu^{(k)})$ be a sequence of compositions 
and $\eta$ be a compositions such that $|\ld^{(j)}|=|\mu^{(j)}|$ and 
$ll(\mu^{(j)}) \le |\eta|$ for all $ 1 \le j \le k.$
There exists a piecewise polynomial function ${\cal K}(t_1,\cdots,t_k):=
{\cal K}_{\sblambda,~\sbmu,\eta}(t_1,\cdots,\gamma_k)$ with rational
 coefficients such that for any non--negative integer
numbers $N_1,\cdots,N_k$, 
$${\cal K}(N_1,\cdots,N_k)=K_{N_{1}\ld^{(1)}+\cdots+N_{k}\ld^{(k)},N_{1}\mu^{(1)}+\cdots+N_{k}\mu^{(k)},\eta}(1).$$
\end{cor}
$(\maltese)$ We {\bf expect} that if all $\mu^{(1)},\cdots,\mu^{(k)}$ are 
partitions, then the restriction of the function

${\cal K}_{~{\sblambda},~{\sbmu},\eta}(t_1,\cdots,t_k)$ on ``the dominant chamber'' 
${\cal N}_k:=\{ (N_1 \ge
N_2 \ge \cdots \ge N_k) \in \Z_{\ge 0}^{k} \}$ is a polynomial 
with {\bf non--negative} rational coefficients. 
\begin{cor} ( {\bf Piecewise polinomiality theorem for  $LR$-numbers} )

Let $\blambda=(\ld^{(1)},\cdots,\ld^{(k)}),\bmu=(\mu^{(1)},\cdots,\mu^{(k)})$  
and $\bnu=(\nu^{(1)},\cdots, \nu^{(k)})$ be three sequences of partitions. 
There  exists a piecewise
polynomial function $LR_{~{\sblambda},~{\sbmu}}^{~{\sbnu}}(t_1,\cdots,t_k)$ 
such that for any non--negative integers $N_1,\cdots,N_k,$
$$LR_{~{\sblambda},~{\sbmu}}^{~{\sbnu}}(N_{1},\cdots,N_{k})=c_{N_{1}\ld^{(1)}+
\cdots+N_{k} \ld^{(k)},N_{1}\mu^{(1)}+\cdots+N_{k} \mu^{(k)}}^{N_{1}\nu^{(1)}
+\cdots+ N_{k}\nu^{(k)}}.$$ 
\end{cor}
$(\maltese)$  We {\bf expect} that the restriction of the function
$LR_{~{\sblambda},{\sbmu}}^{{~\sbnu}}(t_1,\cdots,t_k)$ on ``the dominant 
chamber'' ${\cal N}_k:=\{ (N_1 \ge
N_2 \ge \cdots \ge N_k) \in \Z_{\ge 0}^{k} \}$ is a polynomial 
with {\bf non--negative} rational coefficients.

\begin{prb} Describe the polynomiality domains of the function
$$ (N_1,\cdots,N_k) \rightarrow c_{N_{1}\ld^{(1)}+\cdots+N_{k} \ld^{(k)},N_{1}\mu^{(1)}+\cdots+N_{k} \mu^{(k)}}^{N_{1}\nu^{(1)}+\cdots+
N_{k}\nu^{(k)}}.$$ 
\end{prb}
 
\begin{exs} $(i)$ Take $\ld= (5,3,3,2),$  $\mu=(3,3,3,2,1,1)$ and  
$\eta=(1^{6}).$

 One can check that  

$K_{\ld\mu\eta}(q)=K_{\ld\mu}(q)=q^3(3,5,8,6,5,2,1),$  
$P_{\ld\mu\eta}(q,t)=1+q^4(3,5,4,3,1)~t -                   \\
q^7(1,3,2,1,0,1,3,3,2,1)~t^2- q^{12}(2,9,14,18,18,20,17,14,8,4,1,1)~t^3  \\  
+q^{16}(3,6,10,17,28,35,39,36,30,24,19,11,5,1)t^4     
-q^{21}(-1,0,4,3,6,6,13,16,16,10,5,1,1)t^5      \\
-q^{26}(1,2,7,10,16,19,22,23,23,20,17,10,6,4,1)t^6       
+q^{33}(1,4,8,14,17,20,23,23,19,16,6,1)t^7    \\
-q^{40}(-1,-1,1,4,6,7,3)~t^8 - q^{48}~(1+q+q^2)^2~t^9,$ 

$J_{\ld\mu\eta}(q)=q^3(3,2,3,2,2,1,1).$ In other words,
$$ Q_{\ld\mu\eta}(q,t)= (1-q^3~t)^3(1-q^4~t)^2(1-q^5~t)^3(1-q^6~t)^2(1-q^7~t)^2
(1-q^8~t)(1-q^9~t).
$$
Therefore,  the dimension of the Gelfand--Tsetlin polytope $GT(\ld,\mu)$ is
 equal to 9, and
$$ \sum_{n \ge 0} K_{n\ld,n\mu}(1)~t^n = (1+21~t+78~t^2+64~t^3+9~t^4)/(1-t)^
{10},$$
$$ \sum_{n \ge 0} K_{n\ld,n\mu}(-1)~t^n = (1-3t+6t^2-4t^3+t^4)/(1-t^2)^4(1+t).
$$ 
$(ii)$ Take $\ld=(3,2,1)$ and $\mu=\eta=(1^6)$. Then $K_{\ld\mu}(q)=q^4(1,2,2,3,3,2,2,1),$

$ P_{\ld\mu\eta}(q,t)=1+q^6(1,2,2,1,1)t+q^{12}(1,2,5,4,6,4,3,1,1)t^2  \\
+q^{20}(1,1,1,0,1,-2,-1,-2,-1,-1)t^3 
-q^{29}(2,2,4,4,4,3,3,1)t^4-q^{37}(1,1,2,1,2,1,1)t^5,$
$$ J_{\ld\mu\eta}(q)=q^4(1,2,1,1,1,1,1,1).$$
 Therefore, the dimension of the Gelfand--Tsetlin polytope $GT(\ld,\mu)$ is
 equal to 7, and
$$ \sum_{n \ge 0} K_{n\ld,n\mu}(1)~t^n = (1+8~t+35~t^2+32~t^3+9~t^4)/(1-t)^
{8},$$
$$\sum_{n \ge 0} K_{n\ld,n\mu}(-1)~t^n = (1+5~t^2+3~t^4)/(1-t^2)^4.
$$
\end{exs}
\begin{rem} {\rm We see that in both examples $J_{\ld\mu\eta}(q) \le 
K_{\ld\mu\eta}(q),$ and the initial and the leading terms of the polynomials
$J_{\ld\mu\eta}(q)$ and $K_{\ld\mu\eta}(q)$ are the same.  These 
observations may be { \bf not} true if $\mu$ is an arbitrary composition, 
e.g. if  $\ld=(3,2,1), \mu=(0,2,0,2,0,2)$ and $\eta=(1,2,2,1),$ then

$K_{\ld\mu\eta}(q)=q^4(1,2,-8,-6,8,5),$ but $J_{\ld\mu\eta}(q)=
q^3(2,4,4,3,4,4,3),$
see Examples~4.6. 

$(\spadesuit)$ It was the surprising and unexpected thing for the author to
find that even though $\mu$ and $\eta$ are partitions, the above inequality
$$ J_{\ld\mu\eta}(q) \le K_{\ld\mu\eta}(q)$$
may be wrong. For example, take $\ld=(2,2,2,1,1)$ and $\mu=\eta=(1^8)$. Then
$$ K_{\ld\mu\eta}(q)=q^3(1,1,2,{\bf 2},3,3,4,3,3,2,2,1,1),
~~but ~J_{\ld\mu\eta}(q)= q^3(1,1,2,{\bf 3},2,2,2,1,1,1,1,1,1).$$
Furthermore, one can show that $P_{(2,2,2,1,1),(1^8),(1^8)}(q,t)=$

$1+q^6(-1,1,2,2,2,1,1)t+\cdots+q^{117}(1,1,2,2,2,1,1,-1)t^{13}+q^{130}t^{14},$

see Section~5.4 for more details about the polynomials 
$P_{(2^k,1^n),(1^{2k+n}),(1^{2k+n})}(q,t).$
\\

$(\maltese)$  However, we {\bf expect} that  if ~$\mu$ is a  partition, then
the initial and the leading 
terms of the polynomials $J_{\ld\mu\eta}(q)$ and $K_{\ld\mu\eta}(q)$ are 
the same. 

$(\maltese)$  Moreover, we {\bf expect} that if $\mu$ is an arbitrary 
composition, then

$j_{\max}:=\max \{j \mid j \in J_{\ld\mu\eta} \}=c(\ld,\mu \Vert \eta)$ and 
$\# \{j \in J_{\ld\mu\eta} \mid j=j_{\max} \} \le d(\ld,\mu \Vert \eta),$
see Section~6.4, Rationality conjecture, for more detailed statements.}
\end{rem}

\begin{exs}  $(1)$ Take $\ld^{(1)}=(3,2,1),$ $\ld^{(2)}=(2,2),$ $\mu^{(1)}=
(1^6),$
$\mu^{(2)}=(1^4)$ and $\eta=(1^6).$  Then one can check that 
$$Q(q,x,y)=Q_{(3,2,1),(1^6),(1^6)}(q,x)~Q_{(2,2),(1^4),(1^4)}(q,y)(1-q^7xy)(1-q^8xy),$$
where $Q_{(3,2,1),(1^6),(1^6)}(q,x)=(1-q^{5}x) \prod_{j=4}^{11}(1-q^{j}x),$
see Example~4.20, $(ii),$ and 

$Q_{(2,2),(1^4),(1^4)}(q,y)=(1-q^2y)(1-q^4y),$ $P_{(2,2),(1^4),(1^4)}(q,y)=1.$

The expression for $P(q,x,y)$ is rather long, so we give here only the 
formula for its value at $q=1.$  Namely,

$P(1,x,y)=[1+8x+35x^2+32x^3+9x^4+(6x-44x^2-118x^3-81x^4-18x^5)y$

$+(-3x^2+40x^3+143x^4+66x^5+9x^6)y^2-(16x^4+48x^5+21x^6)y^3](1-x).$ 

Let us remark that in our case $Q(1,x,y)=(1-x)^9(1-y)^2(1-xy)^2,$ and because
of the well-known identity
$$ (1-x_1 \cdots x_k)^{-1}\prod_{j=1}^{k}(1-x_j)^{-1}=\sum_{(n_1,\cdots,n_k)
\in \Z_{\ge 0}^{k}} \min(n_1,\cdots,n_k)x_1^{n_1} \cdots x_k^{n_k},$$

this example shows that the Kostka number $K_{n(3,2,1)+m(2,2),(n^6)+(m^4)}(1)$
considered as a function of $n$ and $m$ on the
set $\{(n,m) \in \Z_{\ge 0}^{2} \},$ has at least {\bf two} different 
polynomiality region, namely,
``the dominant chamber'' ${\cal N}_2=\{ (n,m) \mid n \ge m \}$ and that
$\{ (n,m) \mid n \le m \}.$ Moreover, since
$$ K_{n(3,2,1)+m(2,2),(n^6)+(m^4)}(1)=
c_{n(5,4,3,2,1)+m(3,2,1),n(3,2,1)+m(2,2)}^{n(6,5,4,3,2,1)+m(4,3,2,1)},$$
we see that if 

$\ld^{(1)}=(3,2,1),\ld^{(2)}=(2,2),\mu^{(1)}=(5,4,3,2,1),
\mu^{(2)}=(3,2,1),\nu^{(1)}=(6,5,4,3,2,1)$ and $\nu^{(2)}=(4,3,2,1),$ then

$(\clubsuit)$  the Littlewood--Richardson number  $c_{n\ld^{(1)}+
m\ld^{(2)},n\mu^{(1)}+m\mu^{(2)}}^{n\nu^{(1)}+m\nu^{(2)}}$ considered as a 
function of $n$ and $m$  on the set $\{(n,m) \in \Z_{\ge 0}^{2} \},$ has 
the same ( at least) two different polynomiality regions.

$(2)$ Now take $\ld^{(1)}=(3,2,1),$ $\ld^{(2)}=(2,2,1),$ $\mu^{(1)}=
(1^6),$ $\mu^{(2)}=(1^5)$ and $\eta=(1^6).$  Then one can check that 
$$Q(q,x,y)=Q_{(3,2,1),(1^6),(1^6)}(q,x)~Q_{(2,2,1),(1^5),(1^5)}(q,y),$$
where $Q_{(2,2,1),(1^5),(1^5)}(q,y)=
(1-q^{2}y)(1-q^3y)(1-q^4y)(1-q^5y)(1-q^6y).$

Therefore, in this case the function $ (n,m) \rightarrow K_{n(3,2,1)+m(2,2,1),
n(1^6)+m(1^5),\eta}(1)$ is a polynomial function in $n$ and $m$ on the whole 
set $ \{(n,m) \in \Z_{\ge 0}^{2} \}.$
\end{exs}
It seems interesting to compare the above-described examples with the following
result by P.~McMullen \cite{McM}:

Let $\Delta_1,\cdots,\Delta_k \subset \R^d$ be {\bf integer} convex polytopes, 
and $t_1,\cdots,t_k \in \N^k.$ Given any integer polytope 
$\Gamma \subset \R^d,$  denote by $N(\Gamma):= \# (\Gamma \cap \Z^d).$ \\

$(\spadesuit)${\bf Mixed lattice point enumerator theorem (P.~McMullen, 
\cite{McM})} \\

$N(t_1\Delta_1+\cdots+t_k\Delta_k)$ is a 
{\bf polynomial} in $t_1,\cdots,t_k$ with rational coefficients of total degree
 at most $d$. Moreover, the terms of degree $d$ are given by 
$Vol(t_1\Delta_1+\cdots+t_k\Delta_k),$ the so-called {\it mixed volume} of
the polytopes $\Delta_1,\cdots,\Delta_k.$

In other words, the generating function $\sum_{(n_1,\cdots,n_k) \in 
\Z^k_{ \ge 0}} N(n_1\Delta_1+\cdots+n_k\Delta_k)x^{n_1} \cdots x_k^{n_k}$ is
a {\it rational} function in $x_1,\cdots,x_k$ with the (irredundant) 
denominator of the form $\prod_{j=1}^{k}(1-x_j)^{a_j}$ for some non--negative
 integers $a_{1},\cdots,a_{k}.$

\begin{rem} ( {\bf Parabolic Kostka number $K_{\ld\mu\eta}(1)$ as a function of
$\ld$ and $\mu$} )

{\rm Let $\eta$ be a composition, $l(\eta)=p.$ It follows 
from Theorem~3.25 that on the set
$$ Z_{\eta} = \{(\ld,\mu) \in \Z_{\ge 0}^{n} \times \Z_{\ge 0}^{n}
 \mid \ld_1 \ge \ld_2 \ge \cdots \ge \ld_n, ~~~ \ld - \mu \in Y_{\eta} \} $$
the function $(\ld,\mu) \longrightarrow K_{\ld\mu\eta}(1)$ is a continuous
piecewise polynomial function ${\cal K}_{\eta}(\ld,\mu)$  in 
$\ld_1,\cdots,\ld_n,\mu_1,\cdots,\mu_n$ of degree 
$\sum_{1 \le i < j \le p}\eta_i \eta_j -n+1.$ 

It is a challenge  problem to describe the polynomiality domains
of the function $(\ld,\mu) \longrightarrow K_{\ld\mu\eta}(1),$ and find the
corresponding polynomials ${\cal K}_{\eta}(\ld,\mu).$ In the case $\eta=(1^n)$
a  partial solution to this problem has been done by B.V.~Lidskii \cite{Lid2}.
To the best of our knowledge, if $n \ge 4,$ an explicit description of the
polynomiality domains of the function $(\ld,\mu) \longrightarrow 
K_{\ld\mu\eta}(1)$ is not known.
\begin{exs} $(i)$ Take $n=3,$ so that $\ld =(\ld_1 \ge \ld_2 \ge \ld_3 \ge 0)$
 and $\mu=(\mu_1,\mu_2,\mu_3).$ If $\mu$  is a partition, then
 $$ K_{\ld,\mu}(q)= q^{a(\ld,\mu)}~
\left[\begin{array}{c}\ N_{\ld,\mu}+1\\ \ 1\end{array}\right]_{q},$$
where 

$a(\ld,\mu)= \max \{\ld_1-\mu_1, \ld_1+\ld_2-\mu_1-\mu_2, \ld_1+2\ld_2-
2\mu_1-\mu_2, 2\ld_1+\ld_3-2\mu_1-\mu_2 \},$

$N_{\ld\mu}=\min \{\ld_1-\ld_2, \ld_2-\ld_3, \ld_1-\mu_1,
\ld_1+\ld_2-\mu_1-\mu_2 \}.$

$(\clubsuit)$  In particular, we see that $a(\ld,\mu)$ is a homogeneous 
piecewise linear function in $\ld_1,\ld_2,\ld_3$ and $\mu_1,\mu_2.$

Now let us define ``the dominant chamber''
$$ Z_{(1^3)}^{++}= \{(\ld,\mu) \in Z_{(1^3)} \mid 
\ld_3 \le \mu_2 \le \ld_2 \le \mu_1 \le \ld_1, \ld_1-\ld_2+\ld_3 \le \mu_1 \}.
$$
If $\mu$ is a partition, then

$K_{\ld\mu}(q) \vert_{Z_{(1^3)}^{++}}= K_{\Phi(1^3)}(\ld-\mu |~q)=
q^{\ld_1+\ld_2-\mu_1-\mu_2}
\left[\begin{array}{c}\ \ld_1-\mu_1+1\\ \ 1\end{array}\right]_{q}.$

$(\clubsuit)$ One can check that the domain $Z_{(1^3)}^{++}$ is the maximal 
one among domains $D$ such that 
${\cal K}_{(1^3)}(\ld,\mu) \vert_{D}=1+\ld_1-\mu_1.$

$(ii)$ Take $n=4.$ In this case we don't have a complete  description of 
the polynomiality
domains of the function $(\ld,\mu) \longrightarrow K_{\ld\mu\eta}(1).$ 
 Instead, we are going to describe ``the dominant chamber'' 
$Z_{\eta}^{++}$  for the latter function, i.e. the maximal domain $D$ in 
the set
 
$Z_{\eta}^{+}:=\{(\ld,\mu)\in Z_{\eta} \mid  \ld-\mu \in Y_{\eta}^{+} \}$ 
such that $K_{\ld\mu\eta}(1) \vert_{D}={\cal K}_{\eta}(\ld,\mu) \vert_{D}= 
K_{\Phi(\eta)}(\ld -\mu).$
\begin{pr} Assume that $\eta=(1^4),$ and consider the sets 
$$W_{4}^{(1)}= \{(\ld,\mu) \in Z_{(1^4)}^{+} \mid \mu_i \ge m_{i+1}, i=1,2,3;
 ~~2\mu_2 \ge \ld_2+\ld_3  \}, ~~ and $$
$$W_{4}^{(2)}= \{(\ld,\mu) \in Z_{(1^4)}^{+} \mid \mu_i \ge m_{i+1}, i=1,2,3;
 ~~2\mu_2 \le \ld_2+\ld_3, ~~\ld_1+\ld_3 \le \mu_1+\mu_2  \}.$$
Then 
$$ {\cal K}_{(1^4)}(\ld,\mu) \vert_{W_{4}^{(1)}} = K_{\Phi(1^4)}(\ld-\mu) -
\left(\begin{array}{c} \max(\ld_1+\ld_3-\mu_1-\mu_2,0)+2\\ \ 3\end{array}\right),$$
$${\cal K}_{(1^4)}(\ld,\mu) \vert_{W_{4}^{(2)}} = K_{\Phi(1^4)}(\ld-\mu).$$
\end{pr}
\begin{pr} We have 
$$Z_{(1^4)}^{++}=\{(\ld,\mu) \in Z_{(1^4)} \mid \mu_i \ge m_{i+1}, i=1,2,3;
 ~~ \ld_1+\ld_3 \le \mu_1+\mu_2  \},$$ 
and furthermore, 
$K_{\ld\mu}(q) \vert_{Z_{(1^4)}^{++}}=K_{\Phi(1^4)}(\ld-\mu |~q).$
\end{pr}
\begin{prb} Describe explicitly ``the dominant chamber'' $Z_{\eta}^{++}$ in
general case.
\end{prb}
\end{exs}
At the end of this Remark we would like to say a few words about the 
Littlewood--Richardson numbers $c_{\ld,\mu}^{\nu}$ considered as a function of
$\ld,\mu$ and $\nu.$  To start with, let us consider the following set:
$$ {\cal Z}_{n}:=\{(\ld,\mu,\nu) \in \Z_{\ge 0}^{3n} \mid \ld_1 \ge
\cdots \ge \ld_n, \mu_1 \ge \cdots \ge \mu_n, \nu_1 \ge 
\cdots \ge \nu_n, |\ld|+|\mu|=|\nu| \}.$$
The next Proposition is an easy corollary of Theorem~3.25.
\begin{pr} The Littlewood--Richardson number $c_{\ld,\mu}^{\nu}$ is a 
continuous piecewise polynomial function in $\ld_1,\cdots,\ld_n,\mu_1,\cdots,
\mu_n,\nu_1,\cdots,\nu_n$ on the set ${\cal Z}_{n}.$
\end{pr}
\begin{prb} Describe ``the dominant chamber'' for the function $(\ld,\mu,\nu) 
\rightarrow c_{\ld,\mu}^{\nu}$,  i.e. the maximal domain ${\cal
D}_{n} \subset  {\cal Z}_{n}$ such that the restriction 
$c_{\ld,\mu}^{\nu} \vert_{{\cal D}_{n}}$ is a polynomial with 
non--negative rational coefficients.
\end{prb}
\begin{prb} Generalize the results obtained by B.V.~Lidskii \cite{Lid2} for
the function $(\ld,\mu) \rightarrow K_{\ld\mu}(1),$ to the case of the
function $(\ld,\mu,\nu) \rightarrow c_{\ld,\mu}^{\nu}.$
\end{prb}
 }
\end{rem}

\qed

\begin{rem} {\rm It is not difficult to see that Rationality Theorems~4.14 and
4.17, Polynomiality Theorem (Corollary~4.15) and Corollary~4.18, are still 
valid for the {\it level} ~~$l$ restricted parabolic Kostka numbers 
$K_{\ld\mu\eta}^{(l)}(1)$ and the
{\it level} ~~$l$  restricted parabolic Kostka polynomials 
$K_{\ld\mu\eta}^{(l)}(q).$  Remember that the latter can be defined as follows
\begin{equation}
K_{\ld\mu\eta}^{(l)}(q)=\sum_{w \in \Sigma_{n,l}} (-1)^{l(w)}~K_{\ld,w \circ 
\mu,\eta}(q).
\end{equation}
See Section~2.6 for a explanation of notation we have used.}
\end{rem}

\begin{rem} { \rm In Section~4  we have studied a behavior of the parabolic
Kostka polynomials $K_{n\ld,n\mu,\eta}(q)$ as a function of $n.$ We always
have assumed that a composition $\eta$ is fixed. Here we would like to
discuss briefly what happens if a composition $\eta$ is also varied. A naive
way to vary $\eta,$ say to consider $n\eta,$ gives rise to a trivial result. We
suggest the following way. In order to start, we need one definition, namely,
let $\mu=(\mu_1,\mu_2,\cdots)$ be a composition. Define
$$ \mu^{\langle n \rangle}=(\underbrace{\mu_1,\ldots,\mu_1}_{n},
\underbrace{\mu_2,\ldots,\mu_2}_{n},\cdots).
$$
Let us remark that $(n\mu)'=\mu^{\langle n \rangle}.$
\begin{theorem} There exists the limit
$$ \lim_{n \to \infty} q^{c(n\ld,\mu^{\langle n \rangle} \Vert 
\eta^{\langle n \rangle})}~K_{n\ld,\mu^{\langle n \rangle},
\eta^{\langle n \rangle}}(q^{-1}):=X_{\ld\mu\eta}(q),$$
which is a formal power series in q.
\end{theorem}
$(\maltese)$  We {\bf expect} that if $\mu$ is a partition, then the formal 
power series 
$X_{\ld\mu\eta}(q)$ has non--negative integer coefficients. For example,

$X_{(3,2,1),(1^6),(1^6)}(q)= \prod_{n \ge 1}(1-q^n)^{-2}.$  \\

 However, we would like to remark that the limit
$$ \lim_{n \to \infty} q^{-a(n\ld,\mu^{\langle n \rangle} \Vert 
\eta^{\langle n \rangle})}~~K_{n\ld,\mu^{\langle n \rangle},
\eta^{\langle n \rangle}}(q)
$$
does not exist in general.

Finally, it looks as an interesting {\bf problem} to study the generating 
functions 
$$\sum_{n \ge 0} K_{\ld^{\langle n \rangle},\mu^{\langle n \rangle},
\eta^{\langle n \rangle}}(q)~t^n ~~ and ~~\sum_{n \ge 0} K_{\ld^{\langle n \rangle},\mu^{\langle n \rangle},n\eta}(q)~t^n. $$
$(\maltese)$ We {\bf expect} that the latter  generating function is a rational
function in $q$ and $t.$    }
\end{rem}

\begin{rem} ( {\bf Parabolic Hall--Littlewood polynomials 
$Q_{\mu,\eta}(X;q)$} )

{ \rm Let $\mu$ and $\eta$ be compositions such that $|\eta| \ge
ll(\mu),$ and $X=(x_1,\cdots,x_n)$ be the set of variables. Define the 
{\it modified} parabolic Hall--Littlewood polynomials $Q_{\mu,\eta}'(X;q)$ 
as follows:
$$ Q_{\mu,\eta}'(X;q)=\sum_{\ld}K_{\ld\mu\eta}(q)~s_{\ld}(X),$$
and the parabolic Hall--Littlewood polynomial $Q_{\mu,\eta}(X;q)$ using 
the plethystic transformation:
$$Q_{\mu,\eta}(X;q)= Q_{\mu,\eta}'(X(1-q);q).$$

\begin{theorem}( {\bf Rationality theorem for parabolic Hall--Littlewood 
polynomials} )   \\ 
The generating function
$ \sum_{n \ge 0} Q_{n\mu,\eta}(X;q)~t^n$ is a {\it rational} function in 
 $q,t$ and $X.$
\end{theorem} 

In particular, the generating function $\sum_{n \ge 0} s_{n\ld}(X)~t^n$  
for Schur functions  is a {\it rational} function in  $t$ and $X.$   \\

On the other hand, the generating function for the double Kostka polynomials
 $$ Z_{\ld,\mu}(q,t,x):=\sum_{n \ge 0} K_{n\ld,n\mu}(q,t)~x^n$$
 is a formal power series in $q$, $t$ and $x$ which, in general, cannot 
be equal to any rational function.  }
\end{rem}

\section{Parabolic Kostka polynomials:   \\
 Examples}

\subsection{ Parabolic Kostka and Kostka--Foulkes polynomials }

 $1^0$ [{\bf Kostka--Foulkes and parabolic Kostka polynomials }]  \\
 Let $\ld$ be a partition and $R=(R_1,R_2,\cdots,R_r)$ be
a dominant sequence of rectangular shape partitions.

$(i)$ Let $R_a$ be the single row $(\mu_a)$ for all $a$, and
$\mu:=(\mu_1,\mu_2,\ldots)$ is a partition of length at most $n$. Then
\begin{equation}
K_{\ld R}(q)=K_{\ld\mu}(q), \label{3.9}
\end{equation}
i.e. $K_{\ld R}(q)$ coincides with the Kostka--Foulkes polynomial
$K_{\ld\mu}(q)$.

$(ii)$ Let $R_a$ be the single column $(1^{\eta_a})$ for all $a$, and
$\eta =(\eta_1,\eta_2,\ldots )$. Then
\begin{equation}
K_{\ld R}(q)={\overline K}_{\ld'\eta^+}(q), \label{3.11}
\end{equation}
the cocharge Kostka--Foulkes polynomial, where $\ld'$ is the conjugate of
the partition $\ld$, and $\eta^+$ is the partition obtained by sorting
the parts of $\eta$ into weakly decreasing order. Formula (\ref{3.11})
follows from that (\ref{3.9}) and  Duality Theorem for parabolic
Kostka polynomials. \\

$2^0$ [{ \bf Parabolic Kostka polynomials and Kostant partition function }] 

Let $\gamma\in \Z^n,$  $|\gamma|=0,$ $N$ be an integer such that
$N+n(\gamma_i-\gamma_{i+1})\ge 0$ ~for all $1\le i\le n,$ where we put
 $\gamma_{n+1}=0$ .
Consider partitions $\ld_N=N(n,n-1,\ldots
,2,1)+\gamma$, $\mu_N=N(n,n-1,\ldots ,2,1)$  and composition
$\eta$, $|\eta|=n$. Then
\begin{equation}
K_{\Phi(\eta)}(\gamma|~q)=K_{\ld_N,\mu_N,\eta}(q). \label{5.27}
\end{equation}

$3^0$ [{ \bf Skew Kostka--Foulkes  and parabolic Kostka polynomials}]  \\
Let $\lambda\supset\mu$ be partitions, $l(\ld)=n,$ and $\bnu$ be
a sequence of partitions.  \\
Define $\mu_0=(\mu,\underbrace{0,\cdots,0}_{n-l(\mu)}))$. Then
$$ K_{\ld\setminus\mu,~{\sbnu}}(q) \bdoteq K_{\ld,(\mu_0,~{\sbnu})}(q).$$ 
If $\mu$ is a rectangular shape partition and $R$ is a dominant sequence of 
rectangular shape partitions, then
$$ K_{\ld\setminus\mu,R}(q) \bdoteq K_{\ld,(\mu ,R)^{+}}(q) 
\bdoteq K_{\ld,(\mu_0,R)}(q),$$
where $(\mu,R)^{+}$ denotes a dominant rearrangement of the sequence of
rectangular shape partitions $(\mu,R)$.

\begin{ex}  Let $\ld$ and $\mu$ be partitions, $\mu \subset \ld,$ 
$|\ld \setminus \mu|=N,$  and the complement $\ld \setminus \mu =
\coprod \ld^{(i)}$  is a disjoint union of partitions 
$\ld^{(i)},$ $|\ld^{(i)}|=
n_{i}$, $i=1,\cdots,s.$  Then ~~~$K_{\ld \setminus \mu,(1^{N})}(q)$=
$$q^{N}~\prod_{i=1}^{s} K_{\ld^{(i)},(1^{n_i})}(q)
~\left[\matrix{N\cr n_1,\ldots ,n_s}\right]_{q}=q^{N+\sum n(\ld^{(i)}{'})}
[N]! /~\prod_{i=1}^{s} H_{\ld^{(i)}}(q)=K_{\ld,(\mu,1^{N})}(q),$$
where for any partition $\ld,$ $H_{\ld}(q)$ denotes the  {\it hook} polynomial
corresponding to $\ld,$ see e.g. \cite{Ma}, p.45.

In particular, if $n \ge m,$  then $K_{(n,m),(n,1^m)}(q) \bdoteq 
\left[\begin{array}{c}n\\ m \end{array}\right]_q.$
\end{ex}

\begin{ex} Let $\ld$ and $\mu$ be partitions, $\mu \subset \ld,$ 
$|\ld \setminus \mu|=l,$  and the complement $\ld \setminus \mu =
\coprod \ld^{(i)}$  is a disjoint union of partitions 
$\ld^{(i)},$ $|\ld^{(i)}|=n_{i}$, $i=1,\cdots,s.$  Define partitions
${\widetilde \ld}=(Nl+|\mu|,\ld)$ and ${\widetilde \mu}=(l,\mu).$ Then
$$K_{{\widetilde \ld} \setminus {\widetilde \mu},(l^N)} \bdoteq \prod_{i=1}^{s}
 \left[\begin{array}{c}N\\ \ld^{(i)}{'}\end{array}\right]_q.$$
\end{ex}

We would like to emphasize that, in general, the parabolic Kostka polynomial
$K_{\ld,(\mu,~{\sbnu})}(q)$ is {\it different} from the skew Kostka--Foulkes
polynomial $K_{\ld\setminus\mu,~{\sbnu}}(q).$

For example, take $\ld=(2,2),$  $\mu=(1)$ and $R=(3)$. Then 
$K_{\ld\setminus\mu,R}(q)= K_{\ld,(\mu_0,R)}(q)=0$, 
but $K_{\ld,(\mu,R)}(q) = -1+q$.

$4^0$  [{ \bf Principal specialization of skew Schur functions}] 

Let $\ld \supset\mu$ be partitions, $|\ld\setminus\mu|=r,$  and $N \ge 1$ 
be an integer number. Then
$$s_{\ld\setminus\mu}(1,q,\cdots,q^{N-1})\bdoteq K_{(Nr,~\ld)\setminus(r,~\mu),
~(\underbrace{r,\ldots,r}_{N})}(q).
$$
If $\mu=\emptyset$, then
\begin{equation}
s_{\ld}(1,q,\cdots,q^{N-1})\bdoteq K_{(N |\ld|,~\ld),
~(\underbrace{r,\ldots,r}_{N+1})}(q) \bdoteq \left[\begin{array}
{c}N\\ \ld'\end{array}\right]_q.
\end{equation}
The second equality in (5.42) together with the fermionic formula (5.44) for
the Kostka--Foulkes polynomials, is a crucial step in a combinatorial proof
of {\it unimodality} of the generalized $q$-Gaussian coefficients 
$\left[\begin{array}{c}N\\ \ld \end{array}\right]_q$, see \cite{Kir3} for
details.
\begin{ex} ( {\bf A q-analogue of  Merris' conjecture, cf \cite{Me}, 
\cite{Kir6}} )

Let $\lambda$ and $\mu$ be partitions such that $\lambda \ge \lambda'$
with respect to the dominance partial ordering,  see 
Section~2.1. Then

 $\bullet$  $a(\lambda ,\mu) \ge a(\lambda',\mu).$

$\bullet$ ( {\bf q-Analogue of Merris' conjecture} )
$$K_{\lambda,\mu}(q) \ge q^{n(\lambda')-n(\ld)}K_{\lambda',\mu}(q).$$

{\bf Question:} If the above inequality is true, what is the case of equality ?

For example, the equality  holds for any partition $\ld$ if $\mu=(1^{n})$.  
It's not difficult to see that the equality also holds if 
$$\ld=(n,m,1^{n-2}) ~~~and ~~~\mu=(2^{n-1+[m/2]},\varepsilon_m)$$ 
for some positive integers $n \ge m$ and $m \le 4.$  Here $\varepsilon_m 
= 0 ~or ~1$ according to the parity of $m$.

{\bf Question:} Could it be {\bf true} that these two examples are the only 
{\rm infinite families} of partitions $\ld$ and $\mu$ such that 
$\ld \gvertneqq \ld'$ and $K_{\ld,\mu}=K_{\ld',\mu}$ ?

$(\maltese)$ Moreover, we {\bf expect} that the difference
$$K_{\lambda,\mu}(q) - q^{n(\lambda')-n(\ld)}K_{\lambda',\mu}(q)$$
is a unimodal polynomial (with non-negative integer coefficients).
In particular, 

$(\maltese)$ we {\bf expect} ~that if $\ld \ge \ld',$ then for any positive 
integer $N$ the difference
$$\left[\begin{array}{c}N\\ \ld \end{array}\right]_q - q^{n(\ld')-n(\ld)}
\left[\begin{array}{c}N\\ \ld' \end{array}\right]_q $$
is a unimodal polynomial (with non-negative integer coefficients).
\end{ex}  

$5^0$ [ {\bf Fermionic formula for polynomials   $K_{\ld,R}(q)$} ]

Let $\ld$ be a partition and $R=((\mu_a^{\eta_a}))_{a=1}^p$ be a sequence
of rectangular shape partitions such that
$$|\ld|=\sum_a|R_a|=\sum_a\mu_a\eta_a.$$

\begin{de}\label{d4.1} A configuration of type $(\ld ;R)$ is a sequence
of partitions $\nu =(\nu^{(1)},\nu^{(2)},\ldots )$ such that
$$|\nu^{(k)}|=\sum_{j>k}\ld_j-\sum_{a\ge 1}\mu_a\max (\eta_a-k,0)
=-\sum_{j\le k}\ld_j+\sum_{a\ge 1}\mu_a\min (k,\eta_a)
$$
for each $k\ge 1$.
\end{de}

Note that if $k\ge l(\ld)$ and $k\ge\eta_a$ for all
$a$, then $\nu^{(k)}$ is empty. So that each configuration contains only
a finite number of partitions. In the sequel (except Corollary~\ref{c4.4})
we make the convention that $\nu^{(0)}$ is the empty partition.

For a partition $\mu$ define the number $Q_n(\mu)=\mu_1'+\cdots
+\mu_n'$, which is equal to the number of cells in the first $n$
columns of $\mu$. 

The {\it vacancy} numbers
$P_n^{(k)}(\nu):=P_n^{(k)}(\nu ;R)$ of a  configuration $\nu$ of 
type $(\ld;R)$ are defined by
$$P_n^{(k)}(\nu )=Q_n(\nu^{(k-1)})-2Q_n(\nu^{(k)})+Q_n(\nu^{(k+1)})
+\sum_{a\ge 1}\min (\mu_a,n)\delta_{\eta_a,k}
$$
for $k,n\ge 1$, where $\delta_{i,j}$ is the Kronecker delta.

\begin{de}\label{d4.2} A configuration $\nu$ of type $(\ld ;R)$ is called
admissible, if
$$P_n^{(k)}(\nu ;R)\ge 0~~{\rm for~all}~~ k,n\ge 1.$$
\end{de}

We denote by $C(\ld ;R)$ the set of all admissible configurations
of type $(\ld ;R)$, and call a vacancy number $P_n^{(k)}(\nu
;R)$ {\it essential}, if $m_n(\nu^{(k)})>0$.

Finally, for a configuration $\nu$ of type $(\ld ;R)$ let us define its
{\it charge}
$$c(\nu )=\sum_{k,n\ge 1}\left(\begin{array}{c}\alpha_n^{(k-1)}-
\alpha_n^{(k)}+\sum_a\theta (\eta_a-k)\theta (\mu_a-n)\\ 2\end{array}
\right),
$$
and {\it cocharge}
$$\overline c(\nu )=\sum_{k,n\ge 1}\left(\begin{array}{c}\alpha_n^{(k-1)}
-\alpha_n^{(k)}\\ 2\end{array}\right),
$$
where $\alpha_n^{(k)}=(\nu^{(k)})_n'$ denotes the size of the $n$--th
column of the $k$--th partition $\nu^{(k)}$ of the configuration $\nu$; thus,
$\al_n^{(0)}=0, \forall n \ge 1.$  For any real number $x\in \R$ we put 
$\theta (x)=1$, if $x\ge 0$, and $\theta (x)=0$, if $x<0$.

\begin{theorem}\label{t4.3} ( {\bf Fermionic formula for parabolic Kostka
polynomials \cite{Kir6,KSS}} )

 Let $\ld$ be a partition and $R$
be a dominant sequence of rectangular shape partitions. Then
\begin{equation}
K_{\ld R}(q)=\sum_{\nu}q^{c(\nu )}\prod_{k,n\ge 1}\left[\begin{array}{c}
P_n^{(k)}(\nu ;R)+m_n(\nu^{(k)})\\ m_n(\nu^{(k)})\end{array}\right]_q
\label{4.1}
\end{equation}
summed over all admissible configurations $\nu$ of type $(\ld ;R)$;
$m_n(\ld)$ denotes the number of parts of the partition $\ld$ of size
$n$.
\end{theorem}

\begin{cor}\label{c4.4} {\rm (Fermionic formula for Kostka--Foulkes
polynomials \cite{Kir1})}

 Let $\ld$ and $\mu$ be partitions of the same size. Then
\begin{equation}
K_{\ld\mu}(q)=\sum_{\nu}q^{c(\nu)}\prod_{k,n\ge 1}\left[\begin{array}{c}
P_n^{(k)}(\nu ,\mu)+m_n(\nu^{(k)})\\ m_n(\nu^{(k)})\end{array}\right]_q
\label{4.1a}
\end{equation}
summed over all sequences of partitions $\nu
=\{\nu^{(1)},\nu^{(2)},\ldots\}$ such that

$\bullet$ $|\nu^{(k)}|=\sum_{j>k}\ld_j$, $k=1,2,\ldots$;

$\bullet$  $P_n^{(k)}(\nu ,\mu):=Q_n(\nu^{(k-1)})-2Q_n(\nu^{(k)})
+Q_n(\nu^{(k+1)})\ge 0$ for all $k,n\ge 1$, where by definition we
put $\nu^{(0)}=\mu$;

\begin{equation}
\bullet~~c(\nu ):=\sum_{k,n\ge
1}\left(\begin{array}{c}(\nu^{(k-1)})_n' -(\nu^{(k)})_n'\\
2\end{array}\right).~~~~~~~~~~~~~~~~~~~~~~~~~~~~~\label{4.3*}
\end{equation}
\end{cor}
 
\subsection{ Parabolic Kostka polynomials and Littlewood--Richardson numbers} 

$(1^0)$ Let $\ld ,\mu ,\nu$ be partitions,
$|\nu|=|\ld|+|\mu|$, $l(\ld)=p$, $l(\mu)= s$. Consider
partition
$$\wt\ld =(\ld_1+\mu_1,\ldots ,\ld_1+\mu_s,\ld_1,\ld_2,\ldots ,\ld_p)$$
and a dominant rearrangement $\wt R$ of the sequence of
rectangular shape partitions        \\
$R=\{(\ld_1^s),\nu \}$.  Then
\begin{equation}
K_{\wt\ld ,\wt R}(q)=q^{a(\wt\ld,\wt R)}\{ c_{\ld\mu}^{\nu}+
\cdots +q^{n(\nu)-n(\ld)-n(\mu)}\}, \label{3.3}
\end{equation}

where $c_{\ld\mu}^{\nu}$
denotes the {\it Littlewood--Richardson number}, i.e.
$c_{\ld\mu}^{\nu}={\rm Mult}[V_{\nu}:V_{\ld}\otimes V_{\mu}].$

Furthermore, $ a(\wt\ld,\wt R) \ge \ds\sum_{j\le\ld_1}\nu_j'- |\ld|,$ 
and $a(\wt\ld,\wt R)= \ds\sum_{j\le\ld_1}\nu_j'- |\ld|$ if and only if
$c_{\ld,\mu}^{\nu} \ne 0.$      \\
In other words, if $a(\wt\ld,\wt R)= \ds\sum_{j\le\ld_1}\nu_j'- |\ld|,$ 
then $c_{\ld,\mu}^{\nu} \ne 0,$ and  \\

$(\clubsuit)$ ~the coefficient ~$b(\wt\ld ,\wt R)$ is equal to the 
Littlewood-Richardson number $c_{\ld\mu}^{\nu}=
c_{(\ld_{1}^{s}),~\nu}^{{\wt \ld}}.$    \\

$(\clubsuit\clubsuit)$  Moreover, $K_{\wt\ld ,\wt R}(1)$ is equal to the 
number $\# |Tab^{(2)}(\Lambda^{(2)},\nu)|$ of semistandard
{\it domino} tableaux of the shape $\Lambda^{(2)}$ and content $\nu$, where 
$\Lambda^{(2)}$ is a unique partition such that

$\bullet$     ~~~~~$2$-$core(\Lambda^{(2)})= \emptyset$,

$\bullet$     ~~~~~$2$-$quotient(\Lambda^{(2)})=(\ld,\mu).$ 

The partition $\Lambda^{(2)}:=\Lambda^{(2)}(\ld,\mu)$  can be constructed, 
see e.g. \cite{FFLP}, as follows:

Take an integer $r \ge \max(l(\ld),l(\mu)),$ then 
$$ \Lambda^{(2)}(\ld,\mu)+(2r,2r-1,\cdots,2,1)$$
$$=(2\ld_1+2r-1,\cdots,2\ld_k+2(r-k)+1,\cdots,2\ld_r+1) \cup (2\mu_1+2r,\cdots,2\mu_j+2(r-j),2\mu_r+2).$$
Remember, \cite{Ma}, p.6, that if $\ld$ and $\mu$ are partitions, then 
$\ld \cup \mu$ denotes the partition whose parts are those of $\ld$ and 
$\mu,$ arranged in descending order.
\begin{ex} Take $\ld=\mu=(2,1)$ and $\nu=(3,2,1)$. Then ${\wt \ld}=(4,3,2,1),$
${\wt R}=(3,(2,2),2,1)$ and $K_{{\wt \ld},{\wt R}}(q)=q^2(2,3,1).$ More 
generally,
$$ \sum_{n \ge 0}K_{n{\wt \ld},n{\wt R}}(q)~t^n= (1-q^8t^2)/(1-q^3t)^2(1-q^4t)^3(1-q^5t).$$
It is easy to see that $a({\wt \ld},{\wt R})=3=|\nu|-|\mu|$ and 
$b({\wt \ld},{\wt R})=2=c_{\ld,\mu}^{\nu}.$

Furthermore, $\Lambda:=\Lambda^{(2)}=(4,4,2,2),$ and the spin polynomial 
\cite{CB}, and the 
charge-spin polynomial $K_{\ld\mu}^{\nu}(q,t)$  \cite{Kir9} are equal to:
$$\sum_{T \in Tab^{(2)}(\Lambda,\nu)} t^{spin(T)}=t+3t^2+2t^3, 
~K_{\ld\mu}^{\nu}(q,t)=\sum_{T \in Tab^{(2)}(\Lambda,\nu)} q^{charge(T)}~
t^{spin(T)}=q^3t(1+qt)(1+t+qt).$$
Thus, $c_{\ld,\mu}^{\nu}(t)=K_{\ld\mu}^{\nu}(q,t) \vert_{q^3}=t+t^2,$  
where $c_{\ld,\mu}^{\nu}(t)$ denotes 
the $LLT$  $t$-analog of the $LR$-number $c_{\ld,\mu}^{\nu}.$ 
\end{ex} 
$(\spadesuit)$  Finally we want to remark that $K_{\wt\ld, R}(q)=q^{\sum_{j >
 \ld_{1}}\nu_j'}K_{\wt \ld,\wt R}(q),$ 

and $\deg_{q}K_{\Phi(\ld_1,1^{|\nu|})}(\wt \ld -(\ld_1^{s},\nu) ~|~q)
=n(\ld)+\deg_{q}K_{\wt \ld, R}(q).$
\\

 $(2^0)$ More generally, let $\lambda\supset\mu$ be partitions such that
the complement $\lambda\setminus\mu$ is a disjoint union of
partitions $\ld^{(1)},\ldots, \ld^{(p)}$, and $l(\mu )=m$. Let
$\nu$ be a partition, define composition $\wt\nu =(\mu,\nu)$ and
partition $\eta=(m,1^{|\nu|})$. Then
\begin{equation}
K_{\lambda\wt\nu\eta}(q)=q^{a(\lambda,\mu,\nu)}
(c^{\nu}_{\ld^{(1)},\ldots ,\ld^{(p)}}+\cdots
+q^{n(\nu)-n(\ld^{(1)})-\cdots -n(\ld^{(p)})}), \label{3.4}
\end{equation}
where
$$c^{\nu}_{\ld^{(1)},\ldots ,\ld^{(p)}}:={\rm
Mult}[V_{\nu}:V_{\ld^{(1)}}\otimes\cdots\otimes V_{\ld^{(p)}}]
$$
denotes the (multiple) Littlewood--Richardson coefficient, and
$a(\lambda,\mu, \nu)\in\Z_{\ge 0}$.

$(\clubsuit)$ Moreover, $K_{\lambda\wt\nu\eta}(1)=
 \# |Tab^{(p)}(\Lambda^{(p)},\nu)|$ 
is equal to the number of semistandard $p$-rim hook tableaux of shape 
$\Lambda^{(p)}$ and content $\nu,$  where $\Lambda^{(p)}$ is a unique 
partition such that

$\bullet$  ~$p$-$core(\Lambda^{(p)})= \emptyset$.

$\bullet$  ~$p$-$quotient(\Lambda^{(p)})=(\ld^{(1)},\ld^{(2)},\cdots,\ld^{(p)}).$

Similar to the case $p=2,$ the partition $\Lambda^{(p)}$ can be constructed
as follows:

Take an integer $r \ge \max(l(\ld^{(1)}),\cdots,l(\ld^{(p)})),$ then
$$\Lambda^{(p)}+(pr,pr-1,\cdots,2,1)= \cup_{k=1}^{p}(p\ld_1^{(k)}+p(r-1)+k,
\cdots,p\ld_j^{(k)}+p(r-j)+k,\cdots,p\ld_r^{(k)}+k).$$

We refer the reader to \cite{Ma}, Chapter~I, Section~1, Example~8, for 
definitions of the $p$-core and $p$-quotient of a partition $\Lambda$,
and \cite{LLT} for the definition of semistandard $p$-rim hook tableaux
(domino tableaux in the case $p=2$). 

$(\spadesuit)$  Note also, that the order of parts in the definition of
composition $\wt\nu$ is important.        \\
$(3^0)$ Let $A=\Lambda\setminus\ld$ and $B=M\setminus\mu$ be skew diagrams
and $\nu$ be a partition. Define partitions
$$ \al=((M_1^{\Lambda_1'})+\Lambda)*M,~~ \eta=(\Lambda_1'+M_1',1^{|\nu|}),
~~\gamma=((M_1^{\Lambda_1'}+\lambda)*\mu)
$$
and composition
$$\beta=(((M_1^{\Lambda_1'}+\lambda)*\mu,0^{M_1'-\mu_1'})*\nu).
$$ 
Assume that $|A|+|B|=|\nu|,$   Then
$$(\spadesuit) ~~K_{\al\beta\eta}(q)=q^{|\nu|}~K_{\al \setminus \gamma,\nu}(q)=
q^{|\nu|}~\{c_{A,B}^{\nu}+  higher ~degree ~terms ~in ~q  \}.$$
Therefore, $a(\al,\beta \Vert \eta) \ge |\nu|,$ and $a(\al,\beta \Vert \eta)=
 |\nu|$ if and only if $c_{A,B}^{\nu} \ne 0.$  In this case
$$b(\al,\beta \Vert \eta)=c_{A,B}^{\nu}= \langle s_{A} s_{B},s_{\nu} \rangle,
$$
where $b(\al,\beta \Vert \eta)$ denotes the initial coefficient of the 
polynomial
$K_{\al\beta\eta}(q)$, see Definition~6.1, $s_A$ and $s_B$ denote the skew 
Schur functions  corresponding to the skew diagrams $A$ and $B,$
and $\langle \bullet, \bullet \rangle,$ denotes the scalar product
( the so--called
{\it Redfield--Hall scalar product} ) on the ring of symmetric functions, see
e.g. \cite{Ma}, Chapter I, Section 4. 
 
We don't know any ``nice'' combinatorial interpretation of the numbers 
$K_{\al\beta\eta}(1)$ or $K_{\al\beta\eta}(-1).$
 
For a nice combinatorial description of
the numbers $c_{A,B}^{\nu}$ in terms of  ``pictures'', see \cite{Zl}.

See also Section~6.8 for a slightly different exposition of connections 
between the Littlewood--Richardson numbers and the parabolic Kostka 
polynomials.   
\subsection{ MacMahon polytope and rectangular Narayana numbers \cite{Kir9}} 

 Take $\ld=(n+k,n,n-1,\ldots,2)$ and $\mu =\ld'=(n,n,n-1,n-2,\ldots,2,1^k)$. 
If $n \ge k \ge 1$,
then for any positive integer $N$

$\bullet$ $a(N\ld,N\mu)=(2k-1)N$;

$\bullet$ $b(N\ld,N\mu)={\rm dim}V_{((n-k+1)^{k-1})}^{\g
l(N+k-1)}=\ds\prod_{i=1}^{k-1}\prod_{j=1}^{n-k+1}\frac{N+i+j-1}{i+j-1}.$

In other words, $b(N\ld,N\mu)$ is equal to the number of (weak)
plane partitions of rectangular shape $((n-k+1)^{k-1})$ whose
parts do not exceed $N,$ see e.g. \cite{Ma}, \cite{St3}. It is well-known,
 see e.g. \cite{St3}, \cite{Kir9}, that the number
$b(N\ld,N\mu)$ is equal also to the number
$i(\M_{k-1,n-k+1};N)$ of rational points ${\bf x}$ in the
MacMahon polytope $\M_{k-1,n-k+1}$ such that the points $N{\bf x}$
have integer coordinates.
 The generating function for the
numbers $b(n\ld,n\mu)$ has the following form
$$\sum_{n\ge 0}b(n\ld,n\mu)t^n=\left(\sum_{j=0}^{(k-2)(n-k)}
N(k-1,n-k+1;j)t^j\right)/(1-t)^{(k-1)(n-k+1)+1},
$$
where $N(k,n;j)$, $0\le j\le (k-1)(n-1)$, denote the rectangular
Narayana  numbers. For definition of the rectangular Narayana numbers and
the MacMahon polytope, see  \cite{Kir9}, Section~2, Exercise~1. \\
~~~For the reader's convenience, we  display the numbers $b(N\ld,N\mu)$ for 
small values of $k$ 
and $N$. \\
 ~~If $k=1,$ then $b(N\ld,N\mu)=1$ for all integer numbers $N \ge 1.$  \\
 ~~If $k=2,$ then $b(N\ld,N\mu)= \left(\begin{array}{c}N+n-1\\N\end{array}
 \right).$  \\
 ~~If $N=1,$ then $b(\ld,\mu)= \left(\begin{array}{c}n\\k-1\end{array}
 \right).$   \\
 ~~If $N=2,$ then $b(2\ld,2\mu)=\frac{1}{k}\left(\begin{array}{c}n\\ k-1\end{array}\right)\left(\begin{array}{c}n+1\\k-1\end{array}\right).$  

Thus, the number $b(2\ld,2\mu)$ is equal to the Narayana number
 $N_{k-1,n+1}.$ \\
Note also, that
$$ b(N\ld,N\mu)=K_{N(k,1^{n-k+1}),N(1^n)}(1).$$
More generally, see e.g. \cite{Kir3}, 
$$ K_{N(k,1^{n-k+1}),N(1^n)}(q)=
q^{kN}\ds\prod_{i=1}^{k-1}\prod_{j=1}^{n-k+1}\frac{1-q^{N+i+j-1}}{1-q^{i+j-1}}
= q^{kN}~\left[\begin{array}{c}N\\ \al\end{array}\right]_q,$$
where $\al$ is a rectangular shape partition $((k-1)^{n-k+1}).$  \\
In particular,
$K_{N(k,1^{n-k+1}),N(1^n)}(q)$ is a symmetric and {\it unimodal} polynomial
in $q.$

This example and many others, suggests the following
\begin{prb} Define a $q$-analog of the numbers $d(\ld,\mu \Vert \eta),$
 in particular the numbers $b(\ld,R),$  which generalizes the $q$-analog
 of the $LR$--numbers
introduced by A.~Lascoux, B.~Leclerc and J.-Y.~Thibon, see e.g.\cite{LT}.
\end{prb} 

\subsection{ Gelfand--Tsetlin's polytope $GT((2^k,1^n),(1^{2k+n}))$ }

Let $\ld=(2^k,1^n),$ $k > 0,$ be a two--column  partition, and 
$\mu=\eta=(1^{2k+n}).$ In this Section we are going to study in more details 
the polynomials $P_{k,n}(q,t):=P_{\ld\mu\eta}(q,t), P_{k,n}(t):=
P_{\ld\mu\eta}(1,t),$  $J_{k,n}(q):=J_{\ld\mu\eta}(q),$
 as well as the Gelfand--Tsetlin polytope $GT_{k,n}:=GT(\ld,\mu).$ 

We refer the reader to \cite{KB}, \cite{Kir9}, \cite{LM}, \cite{St3}, vol.2,
  for the definition and basic properties of the Gelfand--Tsetlin polytope 
$GT(\ld,\mu)$ corresponding to a partition $\ld$ and composition $\mu.$  

First of all, let us remember \cite{Kir9} the formula for the dimension of
Gelfand--Tsetlin's polytope $GT(\ld,\mu),$  namely, if $\ld$ and $\mu$ are
partitions, $l(\ld)=r, ~l(\mu)=s,$ then
$$\dim GT(\ld ,\mu)=(r-1)(s-1)-\left(\begin{array}{c}r\\
2\end{array}\right)-\sum_{i=1}^{r}\left(\begin{array}{c}\ld_i'
-\ld_{i+1}'\\ 2\end{array}\right),
$$
where $\ld_i':= {\#}\{j \mid \ld_j \ge i \}.$ 

In particular, $\dim GT_{k,n}=n(2k-1)+(k-1)^2.$
\begin{pr}
$(1)$ $\deg_{t} P_{k,n}(t)= \dim GT_{k,n} +1 -k-n=(k-1)(2n+k-2);$

$(2)$ $P_{k,n}(q,t)= (-1)^{a_{k,n}}q^{b_{k,n}}t^{c_{k,n}}
P_{k,n}(q^{-1},t^{-1}),$

where $a_{k,n}= (\sum_{j \in J_{k,n}} j ) - \dim GT_{k,n} -1,$  and $b_{k,n},$
and $c_{k,n}$ are certain non--negative integers.

In particular, $P_{k,n}(t)$ is a symmetric polynomial (with non--negative 
coefficients).
\end{pr}
$(\spadesuit)$ We will say that a polynomial $P(q,t)$ is a {\it reciprocal} 
one  if it satisfies the following condition:
$$ P(q,t)=(-1)^{a}~q^{b}~t^{c}~P(q^{-1},t^{-1})$$
for some non--negative integers $a,$ $b$ and $c$.
\begin{exs}
$(i)$  Take $k=4, n=0,$ then $\dim GT_{4,0}=9, K_{\ld\mu}(q)=J_{4,0}(q)=$

$q^4(1,0,1,1,2,1,2,1,2,1,1,0,1),$ and
$a_{4,0}=4, b_{4,0}=112, c_{4,0}=10.$  Moreover,
$$P_{4,0}(t)= 1+4t+31t^2+40t^3+31t^4+4t^5+t^6.$$
In particular, the normalized volume of Gelfand--Tsetlin's polytope $GT_{4,0}$
is equal to 

$112=2^4 \cdot 7.$

It seems interesting to compare the above formulae with the corresponding 
formulae for the the Gelfand--Tsetlin polytope corresponding to the conjugate
partition $\ld'=(4,4)$ and the same $\mu$ and $\eta.$ It's not difficult to 
see that   $\dim GT((4,4),(1^8))=5,$  $J_{(4,4),(1^8),(1^8)}(q)=
\{12,14,15,16,18,20,24 \},$ 
$P_{(4,4),(1^8),(1^8)}(t)=(1,8,22,8,1)$ and $P_{(4,4),(1^8),(1^8)}(q,t)$
is a reciprocal polynomial.
In particular, the normalized volume of the polytope $GT((4,4),(1^8))$ is
equal to $40.$

$(ii)$ Take $k=3,~n=2,$ then $\dim GT_{3,2}=14,$  $a_{3,2}=10, b_{3,2}=130, 
c_{3,2}=14,$  $J_{3,2}(q)= q^3(1,1,2,{\bf 3},2,2,2,1,1,1,1,1,1),$  but
$K_{(2^3,1^2),(1^8)}(q)=q^3(1,1,2,{\bf 2},3,3,4,3,3,2,2,1,1).$ Therefore, the 
difference $K_{(2^3,1^2),(1^8)}(q)-J_{3,2}(q)$ is a polynomial with one 
{\bf negative} coefficient. Moreover,
$$P_{3,2}(t)=1+13t+225t^2+1350t^3+4088t^4+5768t^5+4088t^6+1350t^7+225t^8
+13t^9+t^{10}.$$
Therefore, the normalized volume of Gelfand--Tsetlin's polytope $GT_{3,2}$
is equal to 

$17112=2^3\cdot 3 \cdot 23 \cdot 31.$

On the other hand, for the conjugate partition $\ld'=(5,3)$ we have
$\dim GT((5,3),(1^8))=6,$

$J_{(5,3),(1^8),(1^8)}=\{13,14,15,16,17,18,19,22,23,25 \}$
 and $P_{(5,3),(1^8),(1^8)}(t)=(1,21,105,98,20),$ 

and therefore, the
polynomial $P_{(5,3),(1^8),(1^8)}(q,t)$ does not satisfy the condition (2) of
Proposition~5.7.

$(iii)$ Take $k=5,~n=0$, then $\dim GT_{5,0}=16,$ and
$$P_{5,0}(t)= (1,25,718,8059,43679,116840,161912,116840,43679,8059,718,25,1).
$$
In particular, the normalized volume of Gelfand--Tsetlin's polytope $GT_{5,0}$
is equal to 

$500556=2^2\cdot 3 \cdot 7 \cdot 59 \cdot 101.$

Note that $\dim GT((5,5),(1^{10}))=7$,  
$P_{(5,5),(1^{10}),(1^{10})}(t)=(1,34,295,565,295,34,1),$

and $J_{(5,5),(1^{10}),(1^{10})}(q)=\{20,22,23,24,25,26,28,30,32,35,40 \}.$ 

In particular, the normalized volume of the polytope $GT((5,5),(1^{10}))$ is
equal to 

$1225=35^2.$ One can check that $P_{(5,5),(1^{10}),(1^{10})}(q,t)$ is a
reciprocal polynomial.
\end{exs}
$(\clubsuit)$ It is interesting to note that the polytopes 
$GT((n^{k}),(1^{kn}))$
 and $GT((n^{k-1},n-1),(1^{kn-1}))$  have the same (normalized) volumes and the
same $h$-polynomials, i.e.

$P_{(n^{k}),(1^{kn}),(1^{kn})}(t)=P_{(n^{k-1},n-1),(1^{kn-1}),(1^{kn-1})}(t).$ 

However,
the polynomials $P_{(n^k),(1^{kn}),(1^{kn})}(q,t)$ and 
$P_{(n^{k-1},n-1),(1^{kn-1}),(1^{kn-1})}(q,t)$ are different.

For example, $P_{(3,3),(1^6),(1^6)}(q,t)=1+q^{10}t+q^{20}~t^2,$ 
but $P_{(3,2),(1^5),(1^5)}(q,t)=1-q^{20}~t^3.$ 

Moreover, $J_{(3,3),(1^6),(1^6)}=\{6,8,9,12 \}, ~but ~~J_{(3,2),(1^5),(1^5)}=
\{4,5,6,7,8 \}.$   \\

$(\spadesuit)$ It seems an interesting {\bf problem} to find under what
assumptions  on $\ld$, $\mu$ and $\eta$ the polynomial $P_{\ld\mu\eta}(q,t)$ is
a reciprocal one, i.e. satisfies the condition (2) of Proposition~5.7. One
necessary condition is clear: $P_{\ld\mu\eta}(t)$ have to be a symmetric
polynomial. 

$(\maltese)$  We {\bf expect} that the latter condition is also sufficient. \\

For example, the polynomials $P_{(2^k,1^n),(1^{2k+n}),(1^{2k+n})}(q,t)$ are
reciprocal; we {\bf expect} that polynomials 
$P_{(n^k),(1^{nk}),(1^{nk})}(q,t)$
are also reciprocal. However, there are plenty of other cases. For example,

$P_{(4,3,2),(2,1,2,1,2,1),(2^3)}(q,t)=1+q^5(-3,2)~t-3q^{11}(1,-1,1)~t^2+
3q^{16}(1,1,1,-2)t^3$

$+3q^{23}(-2,1,1,1)~t^4-3q^{29}(1,-1,1)~t^5-q^{36}(-2,3)~t^6
+q^{42}~t^7.$

We have also $J_{(4,3,2),(2,1,2,1,2,1),(2^3)}(q)=3q^5(1,1,1).$ 

$(\maltese)$ On the other hand, we {\bf expect} that the polynomials 
$P_{(n,k),(1^{n+k}),(1^{n+k})}(q,t)$ are reciprocal if and only if 
$k=0,1,n-1,n.$                            \\

In the case $k=2$ we can say more:
\begin{pr}
$(1)$  $\deg_{t} P_{2,n}(q,t)=2n,~ \deg_{q} P_{2,n}(q,t)=2n(n+4);$

$(2)$  $q^{2n(n+4)}P_{2,n}(q^{-1},t)=P_{2,n}(q,t);$

$(3)$ $P_{k,n}(q,t)$ is a polynomial with {\bf non--negative} integer 
coefficients;

$(4)$ $P_{2,n}(1,1)=C_{n}~C_{n+1}.$ In other words, the (normalized)
volume of the Gelfand--Tsetlin polytope $GT_{2,n}$ is equal to the product of 
two consecutive Catalan numbers $C_n$ and $C_{n+1};$

$(5)$ $J_{2,n}(q)=q^2(1,1,\underbrace{2,\ldots,2}_{n-1},
\underbrace{1,\ldots,1}_{n+2}),$ and
$$K_{(2^2,1^n),(1^{n+4})}(q)-J_{2,n}(q)=q^6~
\left[\begin{array}{c}n\\ 2 \end{array}\right]_q.$$
\end{pr}
We end this Section by discussion of some properties of the Littlewood--
Richardson coefficients $c_{\delta_{n},\delta_{n}}^{\ld},$ where $\delta_n=
(n-1,n-2,\cdots,1,0)$  denotes the staircase partition of height $n-1.$

Denote by $\kappa(n,m)$ the maximal value of the $LR$-number 
$c_{\delta_{n},\delta_{n}}^{\ld},$ where $\ld$ runs over all partitions such 
that $l(\ld) \le m.$  Let $\upsilon_{n,m}(r)$ denote the number of partitions 
$\ld,$  $l(\ld) \le m,$ such that $c_{\delta_{n},\delta_{n}}^{\ld}=r.$ It is
well--known (theorem by Kostant) that $\upsilon_{n,n}(1)=2^{n-1}.$

$(\maltese)$ We {\bf expect} that if $ n \le m \le 2n-2,$ then 
$\upsilon_{n,m}(1)=3^{m-n}/2^{m-2n+1}.$
\begin{prb} It is not difficult to see that $\upsilon_{n,n}(\kappa(n,n))=1,$
i.e. there exists a unique partition $\ld:=\ld_{max},$ $l(\ld) \le m,$ with 
the maximal value of the Littlewood--Richardson coefficient 
$c_{\delta_{n},\delta_{n}}^{\ld}.$

{\bf Question:}~~~~~{\it How does this unique partition $\ld_{max}$ look 
like ?} \\

$(\maltese)$ We {\bf expect} that if $n=2k+1, k \ge 1,$ then
$$\ld_{max}=(3k+[(k+1)/2]-1,3k-1,3k-2,\cdots,k+1,[k/2]+1).$$
\end{prb}
\subsection{ One dimensional sums and parabolic Kostka polynomials}  

 $(\spadesuit)$ ~{\bf Polynomials ${\cal P}_{\ld\mu}(q)$  and their interpretations
 \cite{Kir5}}
\vskip 0.3cm

In this Example  we summarize different  interpretations and some
properties of an interesting family of polynomials ${\cal P}_{\ld\mu}(q)$ 
which frequently appear in Combinatorics, Algebraic Geometry, Representation
Theory, Statistical Mechanics, ... .
 
\begin{de} The polynomials ${\cal P}_{\ld\mu}(q)$
are defined as the transition coefficients between the modified 
Hall-Littlewood  polynomials and the monomial symmetric functions
\begin{equation}
 Q'_{\ld}(X_n;q)=\sum_{\mu}{\cal P}_{\ld\mu}(q)m_{\mu}(X_n).
\end{equation}
\end{de}

In other words,
\begin{equation}
{\cal P}_{\ld\mu}(q)=\sum_{\eta}K_{\eta\mu}(1)K_{\eta\ld}(q).
\end{equation}
To put this another way, the polynomial  ${\cal P}_{\ld\mu}(q)$ is a
$q$-analog of the multiplicity of weight $\ld$ in the tensor product 
$\otimes_{i}V_{\mu_{i}}.$ \\

The polynomials ${\cal P}_{\ld\mu}(q)$ admit the following interpretations:

${\bf(1^0)}$ [{\bf Inhomogeneous unrestricted one dimensional sum with 
"special boundary conditions"}]
\begin{equation}
{\cal P}_{\ld\mu}(q)=q^{n(\mu')}\sum_{m\in{\cal P}_{\ld\mu}}q^{E(m)},
\end{equation}
summed over the set ${\cal P}_{\ld\mu}$ of all transportation  matrices $m$ of 
type $(\ld ;\mu )$, i.e. the set of all matrices of non--negative 
integers with row sums $\ld_i$ and column sums $\mu_j$; $E(m)$ stands for 
the value of the {\it energy function} $E(p)$ of the path $p$ which corresponds
to the transportation  matrix $m$ under a natural identification, 
see \cite{Kir5}, of the set of paths ${\cal P}_{\mu}(b_{\max},\ld )$  
with that of transportation  matrices ${\cal P}_{\ld\mu}.$  We refer the reader
to  \cite{KMOTU2}, or \cite{Kir5} Subsection~3.1, Example~$1^0,$ for a 
definition of the set of paths ${\cal P}_{\mu}(b_{\max},\ld ).$  

${\bf (2^0)}$ [{\bf Generating function of a generalized mahonian statistics 
$\varphi$  on the set of transportation  matrices ${\cal P}_{\ld\mu}$}]
$${\cal P}_{\ld\mu}(q)=q^{n(\mu')}\sum_{m\in{\cal P}_{\ld\mu}}q^{\varphi 
(m)}.
$$
For the definition and examples of generalized mahonian statistics 
see \cite{Kir5}. For example, the energy function $E(m)$ defines a generalized
mahonian statistics on the set of transportation matrices.

${\bf (3^0)}$ [{\bf The  Poincare polynomial of the partial flag variety 
${\cal F}_{\mu}^{\ld}/{\bf C}$}]
\begin{equation}
{\cal P}_{\ld\mu}(q)=\sum_{i\ge 0}q^{n(\ld )-i}\dim H_{2i}({\cal 
F}^{\ld}_{\mu};{\bf Z}). 
\end{equation}
 This result is due to R.~Hotta  and  N.~Shimomura  \cite{HS}.

${\bf (4^0)}$ [{\bf The number of ${\bf F}_q$--rational points of the 
partial flag variety ${\cal F}_{\mu}^{\ld}/{\bf F}_q$}]
\begin{equation}
 q^{n(\ld )}{\cal P}_{\ld\mu}(q^{-1})={\cal F}_{\mu}^{\ld}({\bf F}_q). 
\end{equation}

${\bf (5^0)}$ [{\bf The number of chains of subgroups
$$\{ e\}\subseteq H^{(1)}\subseteq H^{(2)}\subseteq\cdots\subseteq 
H^{(m)}\subseteq G
$$
in a finite abelian $p$--group $G$ of type $\ld$, such that each subgroup 
$H^{(i)}$ has order $p^{\mu_1+\cdots +\mu_i}$ }]
\begin{equation}
 \al_{\ld}(S;p)=p^{n(\ld )}{\cal P}_{\ld\mu}(p^{-1}), 
\end{equation}
where $S:=S(\mu )=(\mu_1,\mu_1+\mu_2,\ldots ,\mu_1+\mu_2+\cdots +\mu_m)$, 
and $l(\mu )=m+1$.   \\
 For more details, proofs and an interesting history of this result, see e.g
\cite{But}.
 
${\bf (6^0)}$ [{\bf String function of affine Demazure's module 
$V_w(l\Lambda_L)$ corresponding to the element $w=r_{Ln-1}r_{Ln-2}\ldots 
r_{L+2}r_{L+1}r_L$ of the affine Weyl group $W(A_{n-1}^{(1)})$}]
\begin{equation}
 {\cal P}_{(l^L)\mu}(q) \bdoteq \sum_{n\ge 0}\dim V_w(l\Lambda_L)_{\mu 
-n\delta}q^n.
\end{equation}
 This result has been obtained in \cite{KMOTU2}, where one can find 
necessary definitions, proofs  and further details.

${\bf (7^0)}$ [{\bf Generalized $t$--supernomial coefficients $\left[\matrix{\ld\cr\mu}\right]_{t}^{(0)}$ and $t$--multinomial
 coefficients  $T^{(0)}(\ld ;\mu)$ }]
\begin{equation}
\left[\matrix{\ld\cr\mu}\right]_{t}^{(0)}=\sum_{\eta}K_{\eta\mu}\wt 
K_{\eta\ld}(t)=t^{n(\ld )}\sum_{\eta}K_{\eta\mu}K_{\eta\ld}(t^{-1}), 
\end{equation} 
\begin{equation}
T^{(0)}(\ld ;\mu)=t^{-E_{\min}}{\cal P}_{\ld\mu}(t),
\end{equation}
for some known constant $E_{\min}$.

 The coefficients (5.55) and (5.56) are 
natural generalizations of those introduced by A.~Schilling and S.O.~Warnaar 
in the case $l(\mu )=2$, see \cite{Kir4}, \cite{Sc}, \cite{ScW}, \cite{Wa}.

${\bf (8^0)}$ [{\bf Fermionic expression for polynomials 
${\cal P}_{\ld\mu}(q)$}]

 Let $\ld$ be a partition and $\mu$ be a composition, $l(\mu )=n$, then
\begin{equation}
{\cal P}_{\ld\mu}(q)=\sum_{\{\nu\}}q^{c(\{\nu\})}\prod_{k=1}^{n-1}\prod_{i\ge 
1}\left[\matrix{(\nu^{(k+1)})'_i-(\nu^{(k)})'_{i+1}\cr
(\nu^{(k)})'_i-(\nu^{(k)})'_{i+1}}\right]_{q}, 
\end{equation}
summed over all flags of partitions $\nu =\{ 
0=\nu^{(0)}\subset\nu^{(1)}\subset\cdots\subset\nu^{(n)}=\ld\}$, such 
that $|\nu^{(k)}|=\mu_1+\cdots +\mu_k$, $1\le k\le n$, and 
$$c(\{\nu\} )=\sum_{k=0}^{n-1}\sum_{i\ge 
1}\pmatrix{(\nu^{(k+1})'_i-(\nu^{(k)})'_i\cr 2}.
$$
See  \cite{Kir5}, Sections~3 and 4, and  \cite{HKKOTY} , where further
 details and  applications of 
the fermionic formula (5.57) can be found. 

In particular, the fermionic formula (5.57) gives an explicit expression for
the number $ |{\cal F}_{\mu}^{\ld}({\bf F}_q) |$ of rational points of the
partial flag variety ${\cal F}_{\mu}^{\ld}$ over the finite field ${\bf F}_q.$

\begin{prb} Deduce the fermionic formula (5.57) from the Lefschetz
fixed points formula, applied to the Frobenius automorphism of the variety
${\cal F}_{\mu}^{\ld}.$ 
\end{prb}

${\bf (9^0)}$ [{\bf Truncated form  or finitization of the characters and 
branching functions of (some) integrable representations of the affine 
Lie algebra of type $A_{n-1}^{(1)}$ }]

The observation that certain special limits of polynomials
${\cal P}_{\ld\mu}(q)$ and Kostka--Foulkes polynomials may play an
important role in the representation theory of affine Lie algebras
originally was made in \cite{Kir4}.
It was observed in \cite{Kir4}, that the character formula for the 
level 1 vacuum representation $V(\Lambda_0)$ of the affine Lie algebra of 
type $A_{n-1}^{(1)}$ (see, e.g., \cite{Kac}, Chapter~13) can be obtained as an 
appropriate limit $N\to\infty$ of the modified Hall--Littlewood 
polynomials $Q'_{(1^N)}(X_n;q)$. The proof was based on the following 
well--known formula
$${\cal P}_{(1^N)\mu}(q)=q^{n(\mu')}\left[\matrix{N\cr \mu_1,\ldots 
,\mu_n}\right]_q, 
$$
see \cite{Kir4}, (2.28).

The latter observation about a connection between the character 
ch$(V(\Lambda_0))$ and modified Hall-Littlewood polynomials 
$Q'_{(1^N)}(X_n;q),$ immediately implies that the level 1 branching 
functions $b_{\ld}^{\Lambda_0}(q)$ can be obtained as an appropriate 
limit $\ld_N\to\infty$ of the "normalized" Kostka--Foulkes polynomials 
$q^{-A_N}K_{\ld_N,(1^N)}(q)$. We refer the reader to \cite{Kac}, Chapter~12, 
for definitions and basic properties of the branching functions 
$b_{\ld}^{\Lambda}(q)$ corresponding to an integrable representation 
$V(\Lambda )$ of an affine Lie algebra. 

It was conjectured in \cite{Kir4}, 
Conjecture~4, that the similar result should be valid for the branching 
functions $b_{\ld}^{\Lambda}(q)$ corresponding to the integrable highest 
weight $\Lambda$ irreducible representation $V(\Lambda )$ of the affine 
Lie algebra $\wh{sl}(n)$. This conjecture has been proved in \cite{Kir4} 
in the 
following cases: $\wh{sl}(n)$ and $\Lambda =\Lambda_0$, $\wh{sl}(2)$ and 
$\Lambda =l\Lambda_0$, and $\wh{sl}(n)$ and $\Lambda =2\Lambda_0$. It had 
not been long before A.~Nakayashiki and Y.~Yamada \cite{NY} proved this 
conjecture in the case $\wh{sl}(n)$ and $\Lambda =l\Lambda_i$, $0\le i\le 
n-1$. See also \cite{KKN} for another proof of the result by 
A.~Nakayashiki and Y.~Yamada in the case $i=0$. 
The general case has been investigated in \cite{HKKOTY}. 
It happened that in general the so--called thermodynamical Bethe ansatz 
limit of Kostka--Foulkes polynomials gives the branching function of a 
certain {\it reducible} integrable representation of $\wh{sl}(n)$, see details 
in \cite{HKKOTY}. 

 $(\spadesuit \spadesuit)$ [{\bf Parabolic Kostka polynomials and 1D sums}]

Let $\ld ,\mu$ be partitions, $|\ld|=|\mu|$, and $n,N$ be natural
numbers such that $l(\ld)=r \le n$, $l(\mu)=s \le n$, and $N\ge
\ld_1+\mu_1$. Define partitions $\al_N=(N^n)$ and
$$\beta_N=(N - \ld_r,N- \ld_{r-1},\ldots ,N- \ld_1,\mu_1,\mu_2,\ldots ,\mu_s).
$$

 \begin{theorem} ( {\bf Algebraic version of the Robinson-Schensted-Knuth
correspondence} ) 

Let $\ld ,\mu , n, N, \al_N$ and $\beta_N$ be as above. Then

i)~ $K_{\al_N\beta_N}(q)\le K_{\al_{N+1}\beta_{N+1}}(q)$;
\begin{equation}
ii)\hskip 0.27cm
If~~N\ge|\ld|,~~then~~K_{\al_N\beta_N}(q)\bdoteq\sum_{\eta}K_{\eta\ld}(q)
K_{\eta\mu}(q).~~~~~~~~~\label {7.5}
 \end{equation}
 \end{theorem}

 \begin{theorem}\label{t7.9} ( {\bf Algebraic version of the dual 
Robinson-Schensted-Knuth correspondence} )

Let $\ld ,\mu$ be partitions, $|\ld|=|\mu|$,
 $l(\ld)=r \le n$, $N\ge\ld_1$. Define the rectangular shape
 partition $\al_N=(n^N)$  and dominant sequence of rectangular shape partitions
 $R_N=\{\mu ,(1^{N- \ld_r}),\ldots ,(1^{N- \ld_1})\}$. Then

i)~ $K_{\al_NR_N}(q)\le K_{\al_{N+1}R_{N+1}}(q)$;
\begin{equation}
ii)~~~  If~~ N\ge |\ld|,~~ then ~~
K_{\al_NR_N}(q)\bdoteq\sum_{\eta}
 {\overline K}_{\eta\ld}(q)K_{\eta'\mu}(q).~~~~~
 \label{7.6}
 \end{equation}
 \end{theorem}
In particular, the following numbers
$$ K_{(N^n),((n-1)^N,1^N)}(1)=\sum_{\ld \vdash N, ~l(\ld) \le n} 
(K_{\ld,(1^N)}(1))^2
$$
are equal to the number of permutations $w \in \Sigma_N$ such that the all 
increasing subsequences in $w$  have the length at most $n$.
\begin{theorem} ( {\bf 1D sums and parabolic Kostka polynomials} )

$(i)$ Let $\ld$ and $\mu$ be partitions of the same size $n$. Define
partition $\al_N=(N^n)$ and sequence of compositions 
$$\wt\bmu_{N}=
((N- \ld_r,0^{r-1}),(N- \ld_{r-1},0^{r-1}),\cdots,(N- \ld_1,0^{r-1}),\mu).$$
Then 
$$ K_{\al_N,~{\wt\sbmu}_N}(q) \bdoteq \sum_{\eta}K_{\eta,\ld}(1)K_{\eta,\mu}(q)
={\cal P}_{\mu,\ld}(q).
$$
$(ii)$ Keep notation of the previous item, but define
$$ \wt\bmu_N^{(0)}=(N- \ld_r,N- \ld_{r-1},\cdots,N- \ld_1,(\mu_1,0^{(r-1)}),
\cdots,(\mu_s,0^{(r-1)})).$$
Then
$$ K_{\al_N,~{\wt\sbmu}_N^{(0)}}(q) \bdoteq \sum_{\eta}K_{\eta,\ld}(q)
K_{\eta,\mu}(1)={\cal P}_{\ld,\mu}(q). $$
\end{theorem}
\begin{ex} Take $n=6$, $\ld=(2,2,2)$ and $\mu=(2,2,1,1)$. One can take $N=6$.
Then $\al_6=(6,6,6)$, $\bmu_6=((4),(4),(4),(2),(2),(1),(1))$,
$\wt\bmu_6=((4,0,0),(4,0,0),(4,0,0),(2),(2),(1),(1))$, $\wt\bmu_6^{(0)}=
((4),(4),(4),(2,0,0),(2,0,0),(1,0,0),(1,0,0))$, and
$$ \sum_{\eta}K_{\eta,\ld}(q)K_{\eta,\mu}(q)= q^7(1,1,3,3,5,4,6,3,3,2,1,0,1)=
K_{\al_6,~{\sbmu}}(q),$$
$$ \sum_{eta}K_{\eta,\ld}(q)K_{\eta,\mu}(1)=q^{13}(1,4,8,9,7,3,1)=K_{\al_6,
~{\wt\sbmu}_6^{(0)}}(q), $$
$$ \sum_{\eta}K_{\eta,\ld}(1)K_{\eta,\mu}(q)= q^{31}(3,6,9,7,5,2,1)=K_{\al_6,
~{\wt\sbmu}_6}(q). $$
\end{ex}
\begin{con} ( {\bf Summation formulas for parabolic Kostka polynomials} )

$(i)$ Let $\bmu= ({\mu^{(a)}: 
=(\mu_1^{(a)},\cdots,\mu_{\eta_a}^{(a)})})_{a=1}^{r}$ and
$\bnu$ be two sequences of partitions such that $|\bmu|=|\bnu|$. Take
$n:=\sum_{a=1}^{r}\eta_a$ and $N\ge |\bmu|,$ and define the sequence of
partitions  \\
 $\wt\bmu:=(\wt\mu^{(r)},\wt\mu^{(r-1)},\dots,\wt\mu^{(1)})$,
where
$$ \wt\mu^{(a)}:=(N- \mu_{\eta_a}^{(a)},\cdots,N- \mu_{2}^{(a)},
N- \mu_{1}^{(a)}).
$$
Then
$$ K_{(N^n),(~{\wt\sbmu},~{\sbnu})}(q)\bdoteq \sum_{\ld}K_{\ld,~{\sbmu}}(q)
K_{\ld,~{\sbnu}}(q).
$$
(ii) Define the sequence of partitions $\bmu_0 =
 {(\mu_0^{(a)}:=(\mu_1^{(a)},
\cdots,\mu_{\eta_a}^{(a)},0^{(N-\eta_a)}))}_{a=1}^{r}$  and in a similar way
that $\bnu_0.$  Then
$$ K_{(N^n),(~{\wt\sbmu},~{\sbnu}_0 )}(q)\bdoteq\sum_{\ld}K_{\ld,~{\sbmu}}(q)
K_{\ld,~{\sbnu}}(1),$$
$$ K_{(N^n),(~{\wt\sbmu}_0,~{\sbnu} )}(q)\bdoteq\sum_{\ld}K_{\ld,~{\sbmu}}(1)
K_{\ld,~{\sbnu}}(q).$$
\end{con}
\section{Parabolic Kostka polynomials:   \\
 Conjectures}

We keep notation of Section~2. Thus, $\lambda$ is a partition, $\mu$ and
$\eta$ are compositions such that $|\lambda|=|\mu|$, $|\eta|=n,$ and 
$ll(\mu) \le n.$
Let $K_{\lambda\mu\eta}(q)$ denote the parabolic Kostka polynomial as defined
in Section~4.
\begin{de} Let $\lambda$, $\mu$ and $\eta$ be as above, and assume that $K_{\lambda\mu\eta}(q)\ne 0$. Introduce non-zero numbers $b(\lambda,\mu \Vert \eta)$ 
and
$d(\lambda,\mu \Vert \eta)$, and integer numbers $a(\lambda,\mu \Vert \eta)$ 
and $c(\lambda,\mu \Vert \eta)$ via the decomposition
\begin{equation}
 K_{\lambda\mu\eta}(q)= b(\lambda,\mu \Vert \eta)q^{a(\lambda,\mu \Vert \eta)}+
\cdots+d(\lambda,\mu \Vert \eta)q^{c(\lambda,\mu \Vert \eta)}. \label{3.5}
\end{equation}
If $K_{\ld\mu\eta}(q)=0,$ we put by definition, $a(\ld,\mu \Vert \eta)
=b(\ld,\mu \Vert \eta)=c(\ld,\mu \Vert \eta)= d(\ld,\mu \Vert \eta)=0.$
\end{de} 
If a composition $\mu$ is the concatenation of partitions $\mu^{(1)},\mu^{(2)},\cdots,\mu^{(r)}$, we will use notation $a(\lambda,\bmu):=
a(\lambda,\mu \Vert \eta)$, $b(\lambda,\bmu):=b(\lambda,\mu \Vert \eta)$. If
compositions $\mu$ and $\eta$ correspond to a (dominant) sequence of 
rectangular shape partitions $R,$ we will write $a(\ld,R)$ instead of 
$a(\ld,\mu \Vert \eta),$  $b(\ld,R)$ instead of $b(\ld,\mu \Vert \eta),$ 
and so on.

\subsection{ Non-vanishing conjecture} 

\begin{con}

Let $\ld$ be a partition, $\mu$ and $\eta$ be a composition, 
$ll(\mu) \le |\eta|=n.$  
Then    

 $K_{\ld\mu\eta}(q) \ne 0,$ ~~if and only if ~~$\ld -\mu \in Y_{\eta}.$  \\

$(\clubsuit)$  Moreover,~~~~~~~~~~~$K_{\ld\mu\eta}(q) \le K_{\Phi(\eta)}(\ld-\mu |~q),$

and the equality is attained on a certain polyhedral domain ${\cal D}_{\eta}$
in ``the space of parameters'' $Z_{\eta}=\{ (\ld,\mu) \in \Z_{\ge 0}^n \times
\Z_{\ge 0}^n \mid \ld_1 \ge \cdots \ge \ld_n, \ld-\mu \in Y_{\eta} \}.$
\end{con}

\subsection{Positivity conjecture} 
\begin{con}
    
Let $\ld$ be a partition and $\mu,$ and $\eta$ be compositions such that
$|\ld|=|\mu|,$ $ll(\mu) \le |\eta|$ . Then
 $$ d(\lambda,\mu \Vert \eta) \ge 0, ~~and ~~d(\lambda,\mu \Vert \eta) > 0
  \Leftrightarrow \ld - \mu \in Y_{\eta}.$$
\end{con}

\begin{rem}
 { \rm It may happen that the all coefficients of a parabolic 
Kostka
polynomial $K_{\ld,~{\sbmu}}(q)$, except that $d(\ld,\bmu),$ are negative. 
For example, take $\ld = (2,2)$ and        \\
$\bmu = ((0),(1,0),(1,0),(1),(1))$. Then
$$ K_{n\ld,n{\sbmu}}(q) = -q^{7n-1}[n+\sum_{k=1}^{n}(2n-2k+1)q^{k}]
+(n+1)^{2}q^{8n}.
$$
Note, that in our example $b(n\ld,n{\bmu}) = -n, $ $a(n\ld,n{\bmu}) = 7n-1,$
$c(n\ld,n{\sbmu})= 8n,$  $d(n\ld,n{\sbmu})=(n+1)^{2},$ $K_{n\ld,n{\sbmu}}(1)
 = n+1,$  $K_{n\ld,n{\sbmu}}(-1)= (n+1)^{2}$, and
$$ \sum_{n \ge 0} K_{n\ld,n{\sbmu}}(q)~t^n
=(1-q^6(1+3~q-q^2)t+3q^{14}t^2-q^{23}t^3)/(1-q^7t)^2(1-q^8t)^3.
$$
On the other hand,
$$ K_{(2n,2n),(n,n,n,n)}(q)=q^{2n}\left[\begin{array}{c}n+1\\
1\end{array}\right]_{q^2}. 
$$  }
\end{rem} 
\subsection{ Generalized saturation conjecture for parabolic Kostka 
polynomials} 
\begin{con} $(\blacklozenge)$ ( {\bf Generalized Saturation Conjecture} )
 
Let $\ld$ be a partition, and $\mu$ and $\eta$ be 
compositions, then for any integer $N \ge 1$
\begin{equation}
 c(N\lambda,N\mu \Vert \eta)=Nc(\lambda,\mu \Vert \eta).
\end{equation}

$(\blacklozenge \blacklozenge)$ Let $\ld$ and $\mu$ be partitions and $\eta$ 
be a composition, then for any integer $N \ge 1$
$$ a(N\ld,N\mu \Vert \eta)=N~a(\ld,\mu \Vert \eta) 
$$
$(\blacklozenge \blacklozenge \blacklozenge)$  More generally, 
let $\ld^{(1)},\ld^{(2)},\cdots,\ld^{(s)}$ be a sequence of partitions, $\eta$
be a composition and $\mu^{(1)},\mu^{(2)},\cdots,\mu^{(s)}$ be a sequence
of compositions such that $|\ld^{(j)}|=|\mu^{(j)}|$ and 
$ll(\mu^{(j)}) \le |\eta|$ for all $j.$  
Let $N$, $p_1,p_2,\cdots,p_{s}$ be positive integer numbers.

 For each $i, 1 \le i \le N,$ define partitions
\begin{equation} {\widehat \ld}^{(i)}:= \lbrack (\sum_{j \ge 1}^{s} 
p_{j}\ld^{(j)}+N-i)/N 
\rbrack ~~and ~~{\widehat \mu}^{(i)}:= \lbrack (\sum_{j \ge 1}^{s}
p_{j}\mu^{(j)}+N-i)/N \rbrack.
\end{equation}
 Assume that $|{\widehat \ld}^{(j)}|=|{\widehat \mu}^{(j)}|$ for all $j.$ 
Then   
$$ \sum_{j=1}^{s} ~p_{j} ~c(\ld^{(j)},\mu^{(j)} \Vert \eta) = 
\sum_{i=1}^{N}c({\widehat \ld}^{(i)},{\widehat \mu}^{(i)} \Vert \eta).
$$
\end{con}
$(\maltese)$  If $\ld$ and $\mu$ are partitions, then we {\bf expect} the 
similar conjecture for the numbers $a(\ld,\mu \Vert \eta).$

Remember that for any real number $x$ the symbol $\lbrack x \rbrack$ denotes 
the integer part of $x.$
\begin{rem} {\rm  It is not true in general  that the inequality 
\begin{equation}
deg K_{\Phi(\eta)}(w(\ld+\delta)-\mu-\delta |~q) < degK_{\Phi(\eta)}
(\ld-\mu |~q)
\end{equation}
holds for any permutation  $w\in \Sigma_n$, $w\ne id,$ as it happens in 
the case
$\eta=(1^n)$, see Example~4.2.  If it would be so, the Generalized 
Saturation Conjecture 
would follow easily from Saturation Theorem for the parabolic $q$-Kostant 
partition function, see Corollary~3.14. It is also not true in  general  that
\begin{equation}
 c(\ld,\mu \Vert \eta):=deg K_{\ld\mu\eta}(q)=degK_{\Phi(\eta)}(\ld-\mu |~q),
\end{equation}
even if $\mu$ is a dominant sequence of rectangular shape partitions of the
same {\it length} which
is compatible with $\eta$, see Example~4.2. In fact,
it looks a difficult problem to
find an explicit formula for the numbers $c(\ld,\mu \Vert \eta).$

$(\maltese)$ However, we {\bf expect} the validity of the following inequality 
 \begin{equation}K_{\ld\mu\eta}(q) \le K_{\Phi(\eta)}(\ld-\mu |~q),
\end{equation}
and if $\ld$ is a
partition and $R=(R_a:=(\mu_a^k)_{a=1}^{p})$ is a dominant sequence of
rectangular shape partitions of the same {\it length} $k$, then $d(\ld,R)=1$.

$(\maltese)$  By duality, we {\bf expect} that if $R=(R_a:=
(k^{\eta_a})_{a=1}^{p})$ is 
a sequence of rectangular shape partitions of the same {\it width} $k$,
 then $b(\ld,R)=1.$  }
\end{rem}

\subsection{ Rationality conjecture} 
\begin{con}

Let $\ld$ be a partition, and $\mu$ and $\eta$ be compositions such that
$\ld-\mu \in Y_{\eta},$ and (according to Theorem~4.14 )

 $\bullet$  ~~~~~~~~~~~~$\sum_{n \ge 0} K_{n\ld,n\mu,\eta}(q)~t^n =
P_{\ld\mu\eta}(q,t)/Q_{\ld\mu\eta}(q,t),$    

where $P_{\ld\mu\eta}(q,t)$ and $Q_{\ld\mu\eta}(q,t)$ are mutually prime 
polynomials with integer coefficients, $P_{\ld\mu\eta}(0,0)=1,$

$\bullet$ ~~~~~~~~~~~$Q_{\ld\mu\eta}(q,t)= \prod_{j \in J}~(1-q^j~t)^{n_j}$  \\
for some finite set of integers $J:=J_{\ld\mu\eta}=\{j_{min}=j_1 < j_2 < \cdots < j_s=j_{max}\},$ and a set of non--negative integers $\{n_j \}_{j \in J}.$

$(\blacklozenge)$ Let $P_{\ld\mu\eta}(q,t)= \sum_{k \ge 0}~P_{\ld\mu\eta}^{(k)}(q)~t^k,$
$P_{\ld\mu\eta}^{(0)}(q)=1,$ and ( if $P_{\ld\mu\eta}^{(k)}(q) \ne 0$ )
$$P_{\ld\mu\eta}^{(k)}(q)=\beta_{k}(\ld\mu\eta)q^{\al_{k}(\ld\mu\eta)}+\cdots+\delta_{k}(\ld\mu\eta)q^{\gamma_{k}(\ld\mu\eta)}.$$
Then, for all $k > 0$ such that $P_{\ld\mu\eta}^{(k)}(q) \ne 0,$
the following inequalities 
$$\gamma_{k}(\ld\mu\eta) \le k j_{max}$$
have to be valid.  Moreover, if the equality is attained
for some  value of k, then for the corresponding value of $k$ one should have
 ~~~~$\delta_{k}(\ld\mu\eta) \ge 0.$

$(\blacklozenge \blacklozenge)$ If $\ld$ and $\mu$ are partitions, then 
additionally, 
for all $k > 0$  such that $P_{\ld\mu\eta}^{(k)}(q) \ne 0,$ 
the following inequalities
$$ \al_{k}(\ld\mu\eta) \ge k j_{min}$$
have to be valid, and if the equality 
is attained for some $k$, then for the corresponding value of $k$ one should 
have  ~~~$\beta_{k}(\ld\mu\eta) \ge 0.$ \\
It follows from Remark~4.22, $(\spadesuit)$, that the polynomial 
$P_{\ld\mu\eta}^{(1)}(q)$ may have negative integer coefficients.
\end{con}

It is easily seen that Rationality Conjecture, item $(i),$ implies both
Positivity and Generalized Saturation Conjectures. Rationality Conjecture,
item, $(ii),$ implies the item $(ii)$ of Conjecture~6.5.
 
\begin{quest} Does there exist a ``nice'' combinatorial interpretation of 
the set $J:=J_{\ld\mu\eta}$ and the exponents $\{n_j \}_{j \in J}$ 
which have appeared  in Rationality Conjecture ?
\end{quest}
\begin{exs} For the reader's convenience, we list below a few examples of 
the set $J_{\ld\mu\eta}$.

$(i)$  $J_{(3,3,2,1),(2,1,2,1,2,1),(2^3)}= \{3^2,4^3,5^3,6^2 \}$, 
$\deg_{t}P_{\ld\mu\eta}(q,t)=8.$

$(ii)$ $J_{(4,2,2,1),(2,1,2,1,2,1),(2^3)}= \{4^4,5^6,6^3 \},$
$\deg_{t}P_{\ld\mu\eta}(q,t)=10.$

$(iii)$ $J_{(5,4,2,1),(3,2^4,1),(2^3)}= \{4,5^2,6^3,7^2 \},$
$\deg_{t}P_{\ld\mu\eta}(q,t)=5.$

$(iv)$ $J_{(5,4,2,1),(3,2^4,1),(2,1^2,2)}= \{4,5^2,6^3,7^2,8^2,9^2 \},$
$\deg_{t}P_{\ld\mu\eta}(q,t)=9.$

$(v)$ $J_{(2,2),(0^4,1,3),(1^6)}= \{5,6,7,8,9^2,10,11,12,13,15,17 \},$
$\deg_{t}P_{\ld\mu\eta}(q,t)=12.$

$(vi)$ $J_{(4,4,2,2),(2^6),(1^6)}= \{4,6,8^3,10,12,14,16 \},$ 
$\deg_{t}P_{\ld\mu\eta}(q,t)=7$ 

and $P_{(4,4,2,2),(2^6),(1^6)}(q,t)$ is a reciprocal polynomial.

$(vii)$ $J_{(4,3,2,1),(2^5),(1^5)}=\{3,4,5,6,7,8,9,10 \},$ 
$\deg_{t}P_{\ld\mu\eta}(q,t)=6$ 

and $J_{(4,3,2,1),(2^5),(1^5)}(q,t)$ is a reciprocal polynomial.

$(viii)$ $J_{(4,3,2,1),(1^{10}),(2^5)}=\{10,11^3,12^3,13^2,14^2,15^2,16,17 \},$
$\deg_{t}P_{\ld\mu\eta}(q,t)=8,$ 

but $P_{(4,3,2,1),(1^{10}),(2^5)}(q,t)$ is not a reciprocal polynomial.

$(ix)$ $J_{(6,5),(1^{11}),(1^{11})}=\{25,26,27,28,29,30,31,33,34,35,37,40,41,43,45,50 \}$, 

$\deg_{t}P_{\ld\mu\eta}(q,t)=16$ and $P_{(6,5),(1^{11}),(1^{11})}(q,t)$ is a 
reciprocal polynomial.

$(x)$ $J_{(2^4,1),(1^9),(1^9)}(q)=q^4(1,1,1,2,3,2,2,2,3,2,2,2,1,1,1,1,1),$
$K_{(2^4,1),(1^9)}(q) - J_{(2^4,1),(1^9),(1^9)}(q)$

$=q^9(1,2,2,1,2,2,1,2,1),$  $\deg_{t}P_{\ld\mu\eta}(q,t)=23,$ and  
$P_{(2^4,1),(1^9),(1^9)}(q,t)$ is a reciprocal polynomial. 
\end{exs}

 $(\maltese)$  We {\bf expect} that if $\eta_1$ and $\eta_2$ are two 
compositions such that $\eta_2$ is a subdivision of $\eta_1$, then 
$J_{\ld,\mu,\eta_1} \subseteq  J_{\ld,\mu,\eta_2}.$ 

\subsection{Polynomiality conjecture}  
\begin{con}

 $(\blacklozenge)$ Let $\ld$ be a partition, $\mu$ and $\eta$ be compositions. Then

$d(N\lambda,N\mu \Vert \eta)$ is a polynomial in N with {\bf  non-negative} 
rational coefficients of the following form:

 there exist a non-negative integer D and a sequence of non-negative 
integers  \\
$h_0=1,h_1,\cdots,h_D(\ne 0)$ such that
$$d(N\lambda,N\mu \Vert \eta)=\ds\sum_{k=0}^Dh_k\left(\begin{array}{c}N+D-k\\
D\end{array}\right).
$$

$(\blacklozenge \blacklozenge)$  Let $\ld$ be a partition, $\mu$ 
and $\eta$ be compositions, and
$$ q^{c(n\ld,n\mu \Vert \eta)}K_{n\ld,n\mu,\eta}(q^{-1})= \sum_{k \ge 0} 
d_{\ld\mu\eta}(k;n)~q^k ,$$
so that $d_{\ld\mu\eta}(0;N)=d(N\ld,N\mu \Vert \eta).$        

Then  for a fixed $k \ge 0$, there exists a polynomial with rational 
coefficients
$D_{\ld\mu\eta}^{(k)}(t)$ of degree depending only on $\ld,\mu$ and $\eta,$
but not $k,$ such that if $N \ge k,$ then $d_{\ld\mu\eta}(k;N)=
D_{\ld\mu\eta}^{(k)}(N).$   

Hence, there exists the limit
$$\lim_{n \to \infty}q^{c(n\ld,n\mu \Vert \eta)}K_{n\ld,n\mu,\eta}(q^{-1})/
d(n\ld,n\mu \Vert \eta).$$
$(\maltese)$  Moreover, we {\bf expect} that $D_{\ld\mu\eta}^{(0)}(t)$  and
$D_{\ld\mu\eta}^{(1)}(t)$ have {\bf non--negative} coefficients. 

$(\blacklozenge \blacklozenge \blacklozenge)$  Let $\ld$ and $\mu$ be partitions, and
$$ K_{n\ld,n\mu,\eta}(q)= q^{a(n\ld,n\mu \Vert \eta)} \{\sum_{k \ge 0} 
b_{\ld\mu\eta}(k;n)~q^{k} \}, 
$$
so that $b_{\ld\mu\eta}(0;N)=b(N\ld,N\mu \Vert \eta).$   

Then for a fixed $k \ge 0$, there exists a polynomial with rational 
coefficients
$B_{\ld\mu\eta}^{(k)}(t)$ of degree depending only on $\ld,$ $\mu$ and $\eta,$
but not $k,$ such that if $N \ge k,$ then $b_{\ld\mu\eta}(k;N)=
B_{\ld\mu\eta}^{(k)}(N).$ 

Hence,  there exists the limit
$$\lim_{n \to \infty} q^{-a(n\ld,n\mu \Vert \eta)}K_{n\ld,n\mu,\eta}(q)/
b(n\ld,n\mu \Vert \eta) .$$
$(\maltese)$  Moreover, we {\bf expect} that $B_{\ld\mu\eta}^{(0)}(t)$  and
$B_{\ld\mu\eta}^{(1)}(t)$ have {\bf non--negative} coefficients.  

$(\blacklozenge \blacklozenge\blacklozenge \blacklozenge)$ Let 
$\{(\ld_a,\mu^{(a)}) \}_{a=1}^{(r)}$ be a collection of pairs 
$(\ld_a,\mu^{(a)}),$ where for all  $a$, $ 1 \le a \le r,$  $\ld_a$ is a 
partition and $\mu^{(a)}$ is a composition
of the fake length at most $n.$  Let $\eta$
be a composition of size $n.$ Then, there exists a piecewise polynomial 
function ${\cal M}(t_1,\cdots,t_r)$ with rational  coefficients such that
for each $r$-tuples of non-negative integers $(n_1,\cdots, n_r)$ one has
$$ {\cal M}(n_1,\cdots, n_r)=d(n_1 \ld_1+\cdots+n_r \ld_r, n_1 \mu_1+
\cdots+n_r \mu_r \Vert \eta).$$
$(\maltese)$ Moreover, we {\bf expect} that if all compositions $\mu^{(a)}~'s$
are in fact partitions, the the restriction of ${\cal M}(t_1,\cdots,t_r)$ on
``the dominant chamber'' $\{(n_1 \ge \cdots \ge n_r) \in \Z_{\ge 0}^{r} \}$ is
a polynomial wit {\bf non--negative} rational coefficients. 
\end{con}

Let us note that Polynomiality Conjecture, items $(\blacklozenge)$--
$(\blacklozenge \blacklozenge\blacklozenge),$ 
follow from Rationality Conjecture, except the statements  about 
non--negativity.
\begin{ex} Take $\ld=(5,3,3,2)$, $\mu=(3,3,3,2,1,1)$ and $\eta=(1^6)$. Based
on formulas from Example~4.17 (i), one can find that
$$\lim_{n \to \infty} q^{-a(n\ld,n\mu \Vert \eta)}K_{n\ld,n\mu,\eta}(q)/
b(n\ld,n\mu \Vert \eta)=(1-q)^{-6}(1+q)^{-2} .$$
\end{ex}
\begin{rem} {\rm  Even in the case when $\bmu$ is a dominant sequence of 
rectangular shape partitions, the sequence $(h_0,h_1,\cdots,h_D)$ does not
necessarily turn out to be {\it unimodal}. For example, take
 $\ld = (5,4,3,2,1),
\mu= (2,2,2)$ and $\nu = (6,5,4,3,2,1).$  It is not difficult to compute
the corresponding Littlewood--Richardson numbers:
$$ c_{n\ld,n\mu}^{n\nu}= K_{(n(2,2,2),n(1^6))}(1)=\left(\begin{array}{c}n+4\\
4\end{array}\right)+\left(\begin{array}{c}n+2\\
4\end{array}\right)=(n+1)(n+2)(n^2+3n+6)/12.
$$
Hence, in this case $D=2$ and $(h_0,h_1,h_2)=(1,0,1).$ It is not difficult to
check that
$$ \sum_{n \ge 0} K_{n(2,2,2),n(1^6)}(q)t^n= (1+q^{15}~t^2)/(1-q^3t)(1-q^5t)(1-q^6t)(1-q^7t)(1-q^9t).
$$
We see that in our example $P_{\ld\mu\eta}^{(1)}(q)=0$. We  can show that 
$P_{(2^4),(1^8),(1^8)}^{(1)}(q)=0$ as well. 

$(\maltese)$ However, we expect that 
 if  $n \ge 5$, then $P_{(2^n),(1^{2n}),(1^{2n})}^{(1)}(q) \ne 0$. 

For example, $P_{(2^{5}),(1^{10}),(1^{10})}^{(1)}(q)=
q^{10}(-1,1,1,2,2,2,2,3,1,2,1).$  }
\end{rem}
Let us observe that $c_{(5,4,3,2,1),(2,2,2)}^{(6,5,4,3,2,1)}= 5$  is
equal to the third Catalan number $C_3$. More generally, one can show that
$$ c_{(2n-1,2n-2,\cdots,2,1),(2^n)}^{(2n,2n-1,\cdots,2,1)}= 
K_{(2^n),(1^{2n})}(1)= 
\ds\frac{1}{n+1}\left(\begin{array}{c}2n\\ n\end{array}\right)=C_n
$$
is equal to the $n$-th Catalan number.

For definition of unimodal sequences/polynomials see e.g.
\cite{St2}, where one can find a big variety of examples of
unimodal sequences which frequently appear in Algebra,
Combinatorics and Geometry.
\begin{rem} {\rm
In the particular case when $b(\Lambda,R)= c_{\ld,\mu}^{\nu}$, see Section 5.2,
the fact that the function $f_{\ld,\mu}^{\nu}(N):=
 c_{N\ld,N\mu}^{N\nu}$ is a polynomial in $N$ with rational coefficients 
follows from Polynomiality Theorem for parabolic Kostka polynomials, see
Corollary~4.15, and  has 
been proved independently by the several authors: A.~Knutson (unpublished), 
H.~Derksen and J.~Weyman \cite{DW2}, E.~Rassart \cite{Ras}, ... .  }
\end{rem}
We would like  to state separately two particular cases of Conjecture 6.10.

\subsection{ The generalized Fulton, $d(\ld,\mu \Vert \eta)=2$ and 
$d(\ld,\mu \Vert \eta)=3$  conjectures}
 
\begin{con} $(\blacklozenge)$ ( {\bf The generalized Fulton conjecture} )  

 If $d(k\lambda,k\mu \Vert \eta)=1$ for some positive integer $k,$ then 
$d(N\lambda,N\mu \Vert \eta)=1$ for all positive integers $N.$  

$(\blacklozenge \blacklozenge)$ If $d(\lambda,\mu \Vert \eta)=2,$  then 
$d(N\lambda,N\mu \Vert \eta)=N+1$ for all positive
 integers $N.$
\end{con}
 If $d(\ld,\mu \Vert \eta)=3,$ we {\bf expect} that there are only two 
possibilities:

 either $d(N\ld,N\mu \Vert \eta)= 2N+1,$  or  $d(N\ld,N\mu \Vert \eta)
= \left(\begin{array}{c}N+2\\2\end{array}\right).$    

$(\maltese)$  Therefore, we {\bf expect} that the cases 
$d(N\ld,N\mu \Vert \eta)=  \left(\begin{array}{c}N+2\\2\end{array}\right)+
k\left(\begin{array}{c}N\\2\end{array}\right),$

 $1 \le k \le 3,$
do not occur. For example, we don't know whether or not there  exist 
a partition  $\ld$ and a dominant
sequence of rectangular shape partitions $R$ such that $d(\ld,R)= 3,$  \\
 but  $d(2\ld,2R)\ge\ 7.$
\begin{rem} {\rm In the case when the numbers $b(\ld,R)$ coincide with the 
Littlewood--Richardson numbers, see Section 5.2, the Fulton conjecture
has been proved by A.~Knutson, T.~Tao and C.~Woodward \cite{KTW}.}
\end{rem}
\begin{rem} {\rm If $\mu$ is a composition, but {\it not} a partition,
then Conjecture~6.14 $(\blacklozenge)$ is  not, in general,  valid for the 
numbers  $b(\ld,\mu \Vert \eta)$.
For example, take $\ld=(3,2,1)$ and $\bmu=((0),(2,0),(2),(2))$, see 
Examples~4.6. Then
$a(\ld,\bmu)=3, b(\ld,\bmu)=1$, but $a(2\ld,2\bmu)=7, b(2\ld,2\bmu)=3$
and $a(3\ld,3\bmu)=8, b(3\ld,3\bmu)=-1$. In fact, if $n \ge 3$, then
$a(n\ld,n\bmu)=3n-1, b(n\ld,n\bmu)=2-n$. On the other hand, $c(n\ld,n\bmu)= 9n$
and $d(n\ld,n\bmu)=n+1, \forall n \ge 1.$ 
In particular, we see that $b(N\ld,N\bmu)$ becomes a polynomial in $N$ only
starting from $N=3.$  }
\end{rem}

\subsection{ $q$-Log concavity and $P$--positivity conjectures} 

\begin{con} ( { \bf $q$-Log concavity and $P$--positivity conjectures} )   

$(\blacklozenge)$  ( { \bf $q$-Log concavity conjecture for parabolic Kostka 
polynomials} )

~~$(a)$ Let $\lambda$ and $\mu$ be partitions and $\eta$ be a composition. 
Consider the function
$g_N(q):=g_N^{\ld\mu\eta}(q)=K_{N\lambda,N\mu,\eta}(q).$ Then
$$(g_N(q))^2 \ge\ g_{N-1}(q)~g_{N+1}(q).$$
$(\maltese)$ Moreover, we {\bf expect} that if a composition $\eta_2$ is 
a subdivision of that $\eta_1,$ then
$$(g_N^{\ld\mu\eta_{2}}(q))^2 - g_{N-1}^{\ld\mu\eta_{2}}(q)~
g_{N+1}^{\ld\mu\eta_{2}}(q) \ge (g_N^{\ld\mu\eta_{1}}(q))^2 -g_{N-1}^{\ld\mu\eta_{1}}(q)~g_{N+1}^{\ld\mu\eta_{1}}(q) \ge 0.$$ 

$(b)$ More generally, let $\ld^{(1)},\ld^{(2)},\cdots,\ld^{(s)}$ and 
$\mu^{(1)},\cdots,\mu^{(s)}$ be two sequences of 
partitions, and $\eta$ be a composition such that $|\ld^{(j)}|=|\mu^{(j)}|$ 
and $ll(\mu^{(j)}) \le |\eta|,1 \le j \le s.$   
Let $N$, $p_1,p_2,\cdots,p_s$ be positive integer numbers. 
  Assume that $|{\widehat \ld}^{(j)}|=|{\widehat \mu}^{(j)}|$ for all $j$. 
Then   
$$\prod_{j=1}^{s}(K_{\ld^{(j)}, \mu^{(j)},\eta}(q))^{p_j} \le 
\prod_{i=1}^{N} K_{{\widehat \ld}^{(i)},{\widehat \mu}^{(i)},\eta}(q) .$$

See Conjecture~6.5, $(\blacklozenge \blacklozenge \blacklozenge)$,(6.62), 
 for the  explanation of notation ${\widehat \ld}^{(i)}$ and 
${\widehat \mu}^{(i)}.$

$(\maltese)$ In particular, we {\bf expect} that if $\ld:=(\ld^{(1)}+\cdots+
\ld^{(s)})/N$ and $\mu:=(\mu^{(1)}+\cdots+\mu^{(s)})/N$ are partitions, then
$$ \prod_{j=1}^{s}K_{\ld^{(j)}, \mu^{(j)},\eta}(q) \le 
(K_{\ld\mu\eta}(q))^{N}.$$

~~$(c)$  ( {\bf Strong $q$-log concavity conjecture for parabolic Kostka 
polynomials} ) 

 Let $l \ge k \ge r \ge 1$ be integers, $\ld,\mu$ and $\eta$ be as in 
Conjecture~6.17 $(a),$ $g_{n}(q)=K_{n\ld,n\mu,\eta}(q).$  Then
$$ g_k(q)g_l(q) \ge g_{k-r}(q)g_{l+r}(q).$$
$(\maltese)$  Moreover, we {\bf expect} that the difference 
 $g_k(q)g_l(q)- g_{k-r}(q)g_{l+r}(q)$ is a unimodal polynomial.
 
$(\blacklozenge \blacklozenge)$  ( {\bf $P$--positivity conjecture for 
parabolic Kostka numbers} )  

Let $\al\supset\beta$ be partitions,
$l(\al)=r.$ 
Consider the following  polynomial:
$$ g_{\al\setminus\beta}(q):= g_{\al\setminus\beta}^{\ld\mu\eta}(q)=
det(g_{\al_i-\beta_j-i+j}(q))_{1\le i,j \le r}.
$$
Then  $g_{\al\setminus\beta}(1) \ge 0.$ Equivalently, 
$\{g_{N}^{\ld\mu\eta}(1)\}_{N \ge 1}$ is a  P\'olya frequency sequence.
\begin{rem}
  {\rm If $r \ge 3,$ then it's  not true, in general, that all the 
coefficients of polynomial  $g_{\al\setminus\beta}(q)$ are non--negative.
For example, take $\ld=(4,3,2,1),$  $R=((2,2),(2),(2),(1,1))$ and
 $\al=(2,2,2).$ Then ~$g_{\al}(q)=4q^{22}+28q^{23}+\cdots+7q^{34}-q^{35}.$ }
\end{rem} 

We want to state some special cases of Conjecture~6.17 in its own right.
\begin{con}     

$(\blacklozenge)$  ( {\bf The generalized Okounkov conjecture, I} ) 

Let $\lambda$ be a partition and $R$ be a dominant sequence of rectangular 
shape partitions.  Then
$$ (b(N\lambda,NR))^2 \ge\ b((N-1)\lambda,(N-1)R)~b((N+1)\lambda,(N+1)R).$$

$(\blacklozenge \blacklozenge)$  More generally, let $\ld$ and $\mu$ be 
partitions, and $\eta$ be a composition, then the power series
$$B(t)=\sum_{n\ge 0}b(n\ld,n\mu \Vert \eta)t^{n}$$
is a $P$-series.
\end{con}
Remind that a power series $B(t)= \sum_{n\ge 0}b_n t^n$ is called a 
$P$-{ \it series}, if $det(b_{\ld_i-i+j})\ge 0$ for any partition $\ld.$

$(\blacklozenge \blacklozenge \blacklozenge)$ 

 Let $\ld,\mu,\nu$ be partitions, then 
\begin{equation} c_{\lbrack(\ld+\mu+1)/2\rbrack,\lbrack(\ld+\mu)/2\rbrack}^
{\nu} \ge c_{\ld,\mu}^{\nu},
\end{equation}
For a more general conjecture, see Section~6.8.
\end{con}
In the case then $(\ld+\mu)/2$ is a partition, Conjecture~6.19,
$(\blacklozenge \blacklozenge \blacklozenge),$ was stated by A.Okounkov
\cite{Ok1}, Section~2.5.
 More generally,

$(\maltese)$ we {\bf expect} that for a sequence of partitions $\ld^{(1)},
\cdots,\ld^{(p)},$  the difference of products of Schur 
functions
\begin{equation} \prod_{k=1}^{p}s_{\lbrack(\sum_{j} \ld^{(j)}+p-k)/p\rbrack}-
\prod_{j=1}^{p}s_{\ld^{(j)}}
\end{equation}
is a {\it Schur} or {\it $s$-positive}, i.e. the latter difference can be 
written as a linear combination of Schur functions with {\bf non--negative}
(integer) coefficients, cf Conjecture~6.23 $(\blacklozenge).$ 

In the case of the Littlewood--Richardson numbers Conjecture~6.18, 
$(\blacklozenge),$ was stated by A.~Okounkov \cite{Ok1}.

\begin{rem} {\rm The log-concavity of  numbers
$$ dim~V_{\ld}^{\g l(n)} = s_{\ld}(\underbrace{1,\ldots,1}_n)=\left(\begin{array}{c}n\\ \ld'\end{array}\right),$$
which  can be in a natural way identified  with  certain numbers 
$b(\ld,R)$ for some partitions $\ld$ and dominant sequences of rectangular
 shape partitions $R$, see e.g \cite{Kir3}, has been proved 
by A~.Okounkov \cite{Ok2}.

The $q$-log-concavity of the generalized $q$-Gaussian coefficients for general
partition $\ld$ has been proved by A.~Okounkov \cite{Ok2}, and earlier for some
special cases, by L.~Butler, C.~Krattenthaller, B.~Sagan and others.
 In fact, A.~Okounkov has proved more fine
result, namely, that not only the dimension  of an 
irreducible representation  (or its $q$-dimension),
but the whole skew Schur function is log-concave. 

$(\maltese)$  We { \bf expect}, that the 
modified parabolic skew Hall--Littlewood function is $q$-log-concave as well.}
\end{rem}

\subsection{ The generalized Fomin-Fulton-Li-Poon conjectures}

Let $A=\Lambda\setminus\ld$ and $B=M\setminus\mu$ be skew diagrams
and $\nu$ be a partition. Let $\theta$ be a composition such that $l(\nu) \le
|\theta|.$  Define partitions
$$ \al=\al(A,B):=((M_1^{\Lambda_1'})+\Lambda,M),\eta=(\Lambda_1'+M_1',\theta)
$$
and the composition $\beta=\beta(A,B):=((M_1^{\Lambda_1'})+\ld,\mu,0^{M_1'-\mu_1'},\nu).$

 One can prove that the ratio
$$ K_{A,B,\theta}^{\nu}(q):=q^{-|\nu|}~K_{\al\beta\eta}(q)$$
is in fact a polynomial in $q$ with non--negative integer coefficients. 

More generally, cf Section~5.2, let $A^{(1)}=\Lambda^{(1)} \setminus \ld^{(1)},
\cdots,A^{(k)}=\Lambda^{(k)} \setminus \ld^{(k)}$ be a $k$-tuples
of skew diagrams, $\nu$ and $\theta$ be compositions such that 
$ll(\nu) \le |\theta|.$  Define new partitions $\al=
\al(A^{(1)},\cdots,A^{(k)})$ and $\beta= \beta(A^{(1)},\cdots,A^{(k)})$ in the
following way:

$ if \sum_{i \le r-1}\Lambda_1^{(i)'} < j \le 
\sum_{i \le r}\Lambda_1^{(i)'},~for ~some~1 \le r \le k,$

$~then~~
\al_j=\sum_{i=r+1}^{k}\Lambda_1^{(i)}+\Lambda_{j}^{(r)},~\beta_j=
\sum_{i=r+1}^{k}\Lambda_1^{(i)}+\ld_j^{(r)}, $

where we put by definition, $\Lambda^{(0)}=\ld^{(0)}=\emptyset.$
In addition, define $\eta=(\sum_{i=1}^{k-1}\Lambda_1^{(i)'}+\ld_1^{(k)'},
\theta).$

One can prove that the ratio
\begin{equation}  K_{A^{(1)},\cdots,A^{(k)},\theta}^{\nu}(q):
=q^{-|\nu|}K_{\al\beta\eta}(q)
\end{equation}
is a polynomial in $q$ with non--negative integer coefficients.

 The main intention of this Section is to state a few results, examples  and 
conjectures about the latter polynomials.

\begin{pr} If $\theta=(1^{|\nu|}),$ then

$(\clubsuit)$  $K_{A^{(1)},\cdots,A^{(k)},\theta}^{\nu}(0)=
c_{A^{(1)},\cdots,A^{(k)}}^{\nu},$

where $c_{A^{(1)},\cdots,A^{(k)}}^{\nu}:=
\langle s_{A^{(1)}} \cdots s_{A^{(k)}},s_{\nu} \rangle.$  Remember, that 
 $s_{A^{(i)}}$ denotes the skew Schur function corresponding to the skew 
diagram $A^{(i)},$ and  $\langle ~ , \rangle$ denotes the scalar product
( the so--called {\it Redfield--Hall scalar product} ) on the ring of symmetric functions, see e.g. \cite{Ma}, Chapter I, Section 4. 

In particular, if $\ld^{(1)}=\ld^{(2)}=\emptyset,$ 
then $K_{A^{(1)},A^{(2)},\theta}^{\nu}(0)$
is equal to the $LR$-number $c_{\Lambda^{(1)},\Lambda^{(2)}}^{\nu}.$ 

$(\clubsuit \clubsuit)$ If  $\ld^{(1)}= \cdots = \ld^{(k)}=\emptyset,$ then 
the number 
$K_{A^{(1)},\cdots,A^{(k)},\theta}^{\nu}(1)$ is equal to the number of 
semistandard {\it $k$-rim hook} tableaux of content $\nu$ and a certain shape, 
see details in Section~5.2.
\end{pr}

\begin{con} ( {\bf Strong $q$-log concavity conjecture for polynomials 
$K_{A^{(1)},\cdots,A^{(k)},\theta}^{\nu}(q)$} )

$\{K_{mA^{(1)},\cdots,mA^{(k)},\theta}^{m\nu}(q) \}_{m \ge 1}$ is a strong 
$q$-log concave sequence.

In particular,

$$(K_{mA^{(1)},\cdots,mA^{(k)},\theta}^{m\nu}(q))^2 \ge K_{(m+1)A^{(1)},
\cdots,(m+1)A^{(k)},\theta}^{(m+1)\nu}(q)~
K_{(m-1)A^{(1)},\cdots,(m-1)A^{(k)},\theta}^{(m-1)\nu}(q).$$
\end{con}
Now we are going to state a generalization of the Fomin-Fulton-Li
-Poon conjectures I and II,  concerning the $LR$-numbers $c_{A,B}^{\nu},$ 
see \cite{Ok1}, \cite{FFLP}. To start with, we need 
a bit more notation from the papers quoted above.

$(\spadesuit)$ For an ordered  $k$-tuples  $(\ld^{(1)},\cdots,\ld^{(k)})$ of 
partitions with the same number of components $p,$  let 
$\gamma=\cup_{j=1}^{k} \ld^{(j)} =(\gamma_1 \ge \gamma_2 \ge \cdots \ge 
\gamma_{kp})$  be the decreasing
rearrangement of the $\ld_i^{(j)}~'s, ~1 \le j \le k, ~1\le i \le p.$  
Define partitions  
$$ {\wt \ld^{(j)}}=
(\gamma_j,\gamma_{j+k},\gamma_{j+2k},\cdots,\gamma_{j+(p-1)k}), 1\le j \le k.
$$
Now suppose that $(A^{(1)}=\Lambda^{(1)}\setminus\ld^{(1)},\cdots, 
 A^{(k)}=\Lambda^{(k)} \setminus \ld^{(k)})$ is an ordered $k$-tuples of skew
diagrams and $\theta$ is a composition.  Construct a new ordered $k$-tuples 
$({\wt \Lambda}^{(1)}, \cdots,
{\wt \Lambda}^{(k)})$  from the $k$-tuples  $(\Lambda^{(1)},\cdots,
\Lambda^{(k)}),$ and  $({\wt \ld}^{(1)},\cdots,{\wt \ld}^{(k)}) $ from 
the $k$-tuples $(\ld^{(1)},\cdots,\ld^{(k)}).$ It is easy
to see that ${\wt \ld}^{(j)} \subset {\wt \Lambda}^{(j)}, \forall j.$ 
Finally, define ${\wt A^{(j)}}={\wt \Lambda}^{(j)} \setminus {\wt \ld}^{(j)},
 ~1 \le j \le k$ 
 and ${\wt \eta}=({\wt \lambda}_1^{(k)}+\sum_{j \ge 2}^{k}
{\wt \Lambda}_1^{(j)},\theta).$

It is useful to consider the following modification of the above
construction. Namely, for any an ordered $k$-tuples $(\ld^{(1)},\cdots,
\ld^{(k)})$ of partitions with the same number of components,  define a
new ordered $k$-tuples of partitions $({ \ld^{\dag}}^{(1)},\cdots,
{{\ld^{\dag}}^{(k)}):=(({\widetilde {\ld^{(1)}}}{'}})',
\cdots,({\widetilde {\ld^{(k)}}{'}})')$.
In a similar way, for an ordered $k$-tuples $(A^{(1)},\cdots,A^{(k)})$ of
skew diagrams one can define a new an ordered $k$-tuples of skew diagrams
$({A^{\dag}}^{(1)}, \cdots,{A^{\dag}}^{(k)}).$

Remember that for any partition $\ld$ the symbol $\ld'$ stands for the 
conjugate of the partition $\ld.$

$(\spadesuit \spadesuit)$ For an ordered pair $(\ld,\mu)$ of partitions with 
the same number of
components, define a new ordered pair $(\ld^{*},\mu^{*})$ as follows:
$$\ld_{k}^{*}=\ld_k-k+\# \{j | \mu_j-j \ge \ld_k-k \},
~~~\mu_{j}^{*}=\mu_j-j+1+\# \{k | \ld_k-k > \mu_j-j \}.$$
One can show, see \cite{FFLP}, that $\ld^{*}$ and $\mu^{*}$ are partitions
and $|\ld^{*}|+|\mu^{*}|=|\ld|+|\mu|.$

Now suppose that $A=\Lambda \setminus \ld, B=M\setminus \mu$ are two skew
shapes  and $\theta$ is a composition. Construct 
\footnote{ As we learned from the referee, a similar construction was also
considered by F.~Bergeron, R.~Biagnoli and M.~Rosas, see e.g. \cite{B},
\cite{BBR}, or  \cite{MacN}.}
$\Lambda^{*}$ and $M^{*}$ from the pair $(\Lambda,M),$ and $\ld^{*}$ and
$\mu^{*}$ from the pair $(\ld,\mu).$ It is not difficult to see that
$\ld^{*} \subset \Lambda^{*}$ and $\mu^{*} \subset M^{*}.$  Finally,
define $A^{*}=\Lambda^{*} \setminus \ld^{*},$ $B^{*}=M^{*} \setminus \mu^{*}$
 and $\eta^{*}=({\lambda_1^{*}}+{M_1^{*}},\theta).$

Similarly to the previous case $(\spadesuit)$, for an ordered pair
$(\ld,\mu)$ of partitions, construct a new ordered pair of partitions 
$(\ld^{\ddag},\mu^{\ddag}):=(((\ld')^{*})',((\mu')^{*})'),$ and for an 
ordered pair $(A,B)$ of skew diagrams define a new pair of skew diagrams 
$(A^{\ddag},B^{\ddag}).$

$(\clubsuit)$ One can show, cf  \cite{FFLP}, Section~5.1, that 
\begin{equation} ~~if ~(A^{*},B^{*})=(A,B),~then
~(A^{*},B^{*})=\sigma ({\widetilde A},{\widetilde B});
~~ (({\tilde A})^{*},({\tilde B})^{*})= \sigma ({\widetilde A},{\widetilde B}),
\end{equation}
where  $\sigma$  denotes the twist $\sigma(X,Y)=(Y,X).$ 
\footnote{As it was pointed by the referee, the equalities  (6.69) was
also proved by F.~Bergeron, R.~Biagnoli and M.Rosas, see
e.g. \cite{B}, \cite{BBR}.}

Let us remark that the
transformation $(\ld,\mu) \rightarrow (\ld^{*},\mu^{*}):=(\ld,\mu)^{*}$ is 
not one-to-one in general,
e.g. $((4,4),(5,3,1))^{*}=((5,4),(4,2,1))^{*}=((4,3),(5,4,1))^{*}.$

$(\spadesuit \spadesuit \spadesuit)$ For an ordered $k$-tuples
 $(\ld^{(1)},\cdots,\ld^{(k)})$ of 
partitions with the same number of components $p$, define a new ordered 
$k$-tuples 
of partitions $(\lceil \ld^{(1)} \rceil,\cdots,\lceil \ld^{(k)} \rceil ),$ 
cf (6.65), as follows:
$$ (\lceil \ld^{(j)}\rceil)_i=
\lbrack (\sum _{s=1}^{k} \ld_i^{(s)}+k-j)/k \rbrack, ~1 \le j \le k, 
~1 \le i \le p.
$$
Now suppose that $(A^{(1)}=\Lambda^{(1)} \setminus \ld^{(1)},\cdots,
A^{(k)}=\Lambda^{(k)} \setminus \ld^{(k)})$ is an ordered $k$-tuples of skew
diagrams  and $\theta$ is a composition. Construct in an
obvious way a new ordered $k$-tuples of skew diagrams 
$(\lceil A^{(1)} \rceil, \cdots,\lceil A^{(k)} \rceil)$ 
from the $k$-tuples $(\Lambda^{(1)},\cdots,\Lambda^{(k)})$ and that 
$(\ld^{(1)},\cdots,\ld^{(k)}),$ and put $ \eta=
(\lambda_1^{(1)}+\sum_{j \ge 2}^{k}\Lambda_1^{(j)},\theta).$

By analogy with the case $(\spadesuit),$ for any an ordered $k$-tuples 
$(\ld^{(1)},\cdots,\ld^{(k)})$ of partitions with the same number of 
components,  define a  new ordered $k$-tuples of partitions 

$({\ld^{\#}}^{(1)},\cdots,
{\ld^{\#}}^{(k)}):=(({\lceil \ld^{(1)}{'} \rceil})',
\cdots,({\lceil \ld^{(k)}{'} \rceil})')$.

In a similar way, for an ordered $k$-tuples  $(A^{(1)},\cdots,A^{(k)})$ of
skew diagrams one can define a new ordered $k$-tuples of skew diagrams
$({A^{\#}}^{(1)}, \cdots,{A^{\#}}^{(k)}).$

\begin{theorem} For an ordered $k$-tuples of
skew diagrams  $(A^{(1)},\cdots,A^{(k)})$ we have the following
equalities:
\begin{equation} (\lceil A^{(1)} \rceil, \cdots,\lceil A^{(k)} \rceil)= 
({A^{\dag}}^{(1)}, \cdots,{A^{\dag}}^{(k)}),
({\widetilde A}^{(1)},\cdots,{\widetilde A}^{(k)})=
({A^{\#}}^{(1)}, \cdots,{A^{\#}}^{(k)}).
\end{equation}
\end{theorem}
\begin{con} 

$(\blacklozenge)$  ( {\bf The generalized Fomin-Fulton-Li-Poon
conjecture I, cf  \cite{Ok1}, \cite{FFLP}, and (6.66)})

Let $A^{(1)},\cdots,A^{(k)}$  be skew diagrams, $\theta$ be a composition 
and $\nu$ be a partition. Then
\begin{equation}
K_{{\wt A}^{(1)},\cdots,{\wt A}^{(k)},\theta}^{\nu}(q) \ge 
K_{A^{(1)},\cdots,A^{(k)},\theta}^{\nu}(q).
\end{equation}
 ~~Equivalently,
$$K_{\lceil A^{(1)} \rceil, \cdots, \lceil A^{(k)} \rceil, \theta}^{\nu}
(q) \ge K_{A^{(1)},\cdots,A^{(k)},\theta}^{\nu}(q).$$
In particular, ~~~~$c_{{\wt A}^{(1)},\cdots,{\wt  A}^{(k)}}^{\nu} \ge 
c_{A^{(1)},\cdots,A^{(k)}}^{\nu}, 
~~~~c_{{A^{\dag}}^{(1)}, \cdots,{A^{\dag}}^{(k)}}^{\nu} \ge 
c_{A^{(1)},\cdots,A^{(k)}}^{\nu}.$

$(\clubsuit)$  We see that the generalized Fomin-Fulton-Li-Poon conjecture I,
 (6.71), is equivalent to our conjecture (6.67), which in turn, is a 
generalization of that  (6.66). As it was mentioned, in the case when 
$(\ld+\mu)/2$ is a partition, the conjecture (6.66) was stated by 
A.~Okounkov, \cite{Ok1}. \\

$(\blacklozenge \blacklozenge)$ {\bf The generalized Fomin-Fulton-Li-Poon
conjecture II, cf \cite{FFLP}} )
\footnote{As it was pointed by the referee, a generalization of the 
original  Fomin--Fulton--Li--Poon conjecture II, \cite{FFLP}, Conjecture~5.1, 
 to the case of skew diagrams has
been stated also by F.~Bergeron, R.~Biagnoli and M.~Rosas, see e.g. 
 \cite{B}, \cite{BBR}; see also \cite{MacN}.}

Let $A,B,\Lambda,M,\ld,\mu$ and $\theta$ be as in $(\spadesuit \spadesuit),$ 
then
$$K_{A^{*},B^{*},\theta}^{\nu}(q) \ge K_{A,B,\theta}^{\nu}(q).$$
In particular,~~~$c_{A^{*},B^{*}}^{\nu} \ge c_{A,B}^{\nu},$
~~~$c_{A^{\ddag},B^{\ddag}}^{\nu} \ge c_{A,B}^{\nu}.$     \\

$(\blacklozenge \blacklozenge \blacklozenge)$  ({\bf The generalized 
Okounkov  conjecture  II} )

Let $(A^{(1)}, \cdots,A^{(k)})$ be an ordered $k$-tuples of skew diagrams,
$\nu$ be a partition and $\theta$ be a composition. Let $p_1,\cdots,p_k$ be
non--negative rational numbers, $p_1+\cdots+p_k =1.$ Define ${\check A}^{(i)}=
\sum_{j=1}^{k}p_{k+1-i+j}~A^{(j)}, 1 \le i \le k.$ Assume that the all
${\check A}^{(1)},\cdots,{\check A}^{(k)}$ are skew diagrams. Then
$$ K_{{\check A}^{(1)},\cdots,{\check A}^{(k)},\theta}^{\nu}(q) \ge
K_{A^{(1)},\cdots,A^{(k)},\theta}^{\nu}(q).$$
In particular, $c_{{\check A}^{(1)},\cdots,{\check A}^{(k)}}^{\nu} \ge 
c_{A^{(1)},\cdots,A^{(k)}}^{\nu}.$
\end{con} 
\begin{exs} We elucidate Conjecture~6.23 in the case $k=2.$ To simplify
notation we will write $A,B,\Lambda,M,\ld$ and $\mu$ instead of $A^{(1)},
A^{(2)},\Lambda^{(1)},\Lambda^{(2)},\ld^{(1)}$ and $\ld^{(2)}$ 
correspondingly. 

$(i)$  Take $\Lambda=(5,1), M=(4,3,1),\nu=(6,5,2,1), \theta=(1^4)$
and $\ld=\mu= \emptyset.$ It is easy to check that 

 $({\wt \Lambda},{\wt M})=
((5,3,1),(4,1))=(\Lambda^{\#},M^{\#}); ~(\Lambda^{*},M^{*})=((4,1),(5,3,1));$ 

$(\lceil \Lambda+M \rceil,\lbrack \Lambda+M \rbrack )=((5,2,1),(4,2))=(\Lambda^{\dag},M^{\dag});$ 

$(\Lambda^{\ddag},M^{\ddag})=((4,2),(5,2,1))$  and $\eta=(4,1^4).$
 
Using the fermionic formula
 (5.44) for Kostka-Foulkes polynomials, one can find that

$K_{\Lambda^{*},M^{*},\theta}^{\nu}(q)=K_{{\wt \Lambda},{\wt M},\theta}^{ \nu}(q)= 
(3,11,18,17,11,4,1),$ 
$K_{\lceil \Lambda+M \rceil,~\lbrack \Lambda+M \rbrack,\theta}^{\nu}(q)=
(3,12,19,18,11,4,1),$ 

$K_{\Lambda,M,\theta}^{\nu}(q)=(1,6,12,14,10,4,1).$ 

Therefore, the difference $K_{{\wt \Lambda},{\wt M},\theta}^{ \nu}(q)-K_{\Lambda,M,\theta}^{\nu}(q)$ is equal to $(2,5,6,3,1).$ 

Similar computations show that if we take  $\theta_1=(1^2,2),$ then  
(with $\eta_1={\wt \eta_1}=(4,1^2,2)$)

$K_{\Lambda^{*},M^{*},\theta_1}^{\nu}(q)=K_{{\wt \Lambda},{\wt M},\theta_1}^{ \nu}(q)=(3,9,13,10,5,1),
K_{\Lambda,M,\theta_1}^{\nu}(q)=(1,5,9,9,5,1)$  and 

$K_{\lceil \Lambda+M \rceil,~\lbrack \Lambda+M \rbrack,\theta_1}^{\nu}(q)= 
(3,10,14,11,5,1).$       \\

$(ii)$ Take $A=(5,5,2,2) \setminus (3,1),$ $B=(1,1) \setminus (1),$ $\nu=
(5,3,2,1)$ and $\theta=(1^4).$ It is easy to check that 

$({\wt A},{\wt B})=((5,2,1) \setminus (3,1),(5,2,1)\setminus (1))=
(A^{\#},B^{\#});$
 
$(A^{*},B^{*})=((4,3,1) \setminus (2),(3,2,2,1) \setminus (2,1));
~(A^{\ddag},B^{\ddag})=((2,2,1)\setminus (1),(5,4,1,1) \setminus (3,1));$ 

$(\lceil A+B \rceil,~\lbrack A+B \rbrack)=((3,3,1,1) \setminus (2,1),(3,3,1,1)
\setminus (2))= (A^{\dag},B^{\dag}).$

Using the fermionic formula (5.44) for Kostka--Foulkes polynomials, one
can find that

$K_{A^{*},B^{*},\theta}^{\nu}(q)=(33,82,86,53,21,6,1),
~K_{\lceil A+B \rceil,~\lbrack A+B \rbrack,\theta}^{\nu}(q)=
(12,20,14,5,1),$

$K_{{\wt A},{\wt B},\theta}^{ \nu}(q)=(20,86,139,131,86,43,17,5,1),
~K_{A^{\ddag},B^{\ddag},\theta}^{\nu}(q)=(22,56,61,40,17,5,1),$

~$K_{A,B,\theta}^{\nu}(q)=(4,9,9,4,1).$

Similar computations show that if we take $\eta_1=(1,2,1),$ then 

$K_{A^{*},B^{*},\theta_1}^{\nu}(q)=(33,64,41,9),~
K_{A,B,\theta_1}^{\nu}(q)=(4,7,3)),
~K_{\lceil A+B \rceil,~\lbrack A+B \rbrack, \theta_1}^{\nu}(q)=
(12,15,5),$

$ K_{{\wt A},{\wt B},\theta_1}^{\nu}(q)=(20,73,87,49,13,1),
~K_{A^{\ddag},B^{\ddag},\theta_1}^{\nu}(q)=(22,45,32,9).$
\end{exs}
These examples show that, probably, there are no simple relationships between
polynomials $K_{A^{*},B^{*},\theta}^{\nu}(q),
 K_{{\widetilde A},{\widetilde B},\theta}^{ \nu}(q), 
 K_{A^{\ddag},B^{\ddag},\theta}^{\nu}(q)$  
and  $K_{\lceil A+B \rceil,~\lbrack A+B \rbrack, \theta}^{\nu}(q).$

$(\maltese)$  However, based on examples, we {\bf expect} that
~~~$ \max \{c_{A^{*},B^{*}}^{\nu},c_{{\wt \Lambda},{\wt M}}^{\nu}\} \ge 
c_{\lceil A+B \rceil,~\lbrack A+B \rbrack}^{\nu}.$

$(\maltese)$   We {\bf expect} that if $\theta_1$ and $\theta_2$ are 
compositions such that
$\theta_2$ is a subdivision of $\theta_1,$ see Section~1, Notation, then
$$K_{{\wt A}^{(1)},\cdots, {\wt A}^{(k)}, \theta_2}^{ \nu}(q)- 
K_{A^{(1)},\cdots,A^{(k)},\theta_2}^{\nu}(q) \ge 
K_{{\wt A}^{(1)},\cdots,{\wt A}^{(k)},\theta_1}^{ \nu}(q)- 
K_{A^{(1)},\cdots,A^{(k)},\theta_1}^{\nu}(q) \ge 0,$$ 
$$K_{A^{*},M^{*},\theta_2}^{ \nu}(q)- 
K_{A,B,\theta_2}^{\nu}(q) \ge 
K_{A^{*}, B^{*}, \theta_1}^{ \nu}(q)- 
K_{A,B,\theta_1}^{\nu}(q) \ge 0,$$ 
$$K_{A^{\ddag},M^{\ddag},\theta_2}^{ \nu}(q)- 
K_{A,B,\theta_2}^{\nu}(q) \ge 
K_{A^{\ddag}, B^{\ddag}, \theta_1}^{ \nu}(q)- 
K_{A,B,\theta_1}^{\nu}(q) \ge 0,$$ 
$$K_{\lceil A^{(1)} \rceil,\cdots,\lbrack A^{k)} \rbrack, \theta_2}^{\nu}(q)-
K_{A^{(1)},\cdots,A^{(k)},\theta_2}^{\nu}(q) \ge
K_{\lceil A^{(1)} \rceil,\cdots,\lbrack A^{(k)} \rbrack,\theta_1}^{\nu}(q) -
K_{A^{(k)},\cdots,A^{(k)},\theta_1}^{\nu}(q) \ge 0,
$$
$$K_{{\check A}^{(1)},\cdots, {\check A}^{(k)}, \theta_2}^{ \nu}(q)- 
K_{A^{(1)},\cdots,A^{(k)},\theta_2}^{\nu}(q) \ge 
K_{{\check A}^{(1)},\cdots,{\check A}^{(k)},\theta_1}^{ \nu}(q)- 
K_{A^{(1)},\cdots,A^{(k)},\theta_1}^{\nu}(q) \ge 0.$$ 

\begin{rem} {\rm  We {\bf expect} that Conjecture~6.3 ( {\bf Positivity}), 
Conjecture~6.5 ({\bf Saturation}), Conjecture~6.7 ({\bf Rationality}),
Conjecture~6.10 ({\bf Polynomiality}), Conjecture~6.17 ({\bf $q$-Log concavity and $P$-positivity}), Conjecture~6.24 ({\bf Generalized Fomin-Fulton-Li-Poon's 
conjectures I and II}) are still valid for 
the level ~$l$-restricted parabolic Kostka 
polynomials $K_{\ld\mu\eta}^{(l)}(q),$ see Remark~4.28 for the definition of
the latter. }
\end{rem}

\subsection{Miscellany}
\begin{con} ( {\bf Rationality conjecture for the $LLT$ $q$-analog of  
$LR$-numbers} ) 

Let $\ld,\mu$ and $\nu$ be partitions, and $c_{\ld,\mu}^{\nu}(q)$ stands for
the $q$-analog of Littlewood--Richardson numbers defined in \cite{CB}, 
\cite{LLT}. Then
$$ \sum_{n \ge 0} c_{n\ld,n\mu}^{n\nu}(q)~t^n =P_{\ld,\mu}^{\nu}(q,t)/
Q_{\ld,\mu}^{\nu}(q,t),
$$
where $P_{\ld,\mu}^{\nu}(q,t)$ and $Q_{\ld,\mu}^{\nu}(q,t)$ are mutually prime
polynomials with integer coefficients.  \\
Moreover,
$$ Q_{\ld,\mu}^{\nu}(q,t)=\prod_{i \in I}~(1-q^i~t)^{n_i}
$$
for some finite set of integers $I=I_{\ld,\mu}^{\nu}$, and a set of positive 
integers $n_i$, $i \in I$.
\end{con}
$(\maltese)$ We expect the similar conjecture for the  parabolic
Kazhdan--Lusztig polynomials, see e.g. \cite{LT} for the definition of the 
latter.
\begin{con} ( {\bf Saturation conjecture for the structural constants of the
multiplication of the Schubert polynomials} )

For each $n \ge 1$, let $\Sigma^{(n)}$ denote the set of all permutations $w$ 
such that the code of $w$ has length at most $n.$ Denote by $\Sigma^{(\infty)}$
the union $\bigcup_{n \ge 1} \Sigma^{(n)}$. 

If $w \in \Sigma^{(n)}$ and
$N \ge 1$ is an integer, define the permutation $N \ast w \in \Sigma^{(Nn)}$
to be a unique permutation with the code $(Nc_1,\cdots,Nc_n),$ where $(c_1,
\cdots,c_n)$ is the code of $w.$ 

For each $w \in \Sigma^{(n)}$ denote
by ${\s}_{w} \in P_n:= \Z[x_1,\cdots,x_n]$ the corresponding 
{\it Schubert polynomial}. It is well--known that the ${\s}_{w}$,
$w \in \Sigma^{(n)},$ form a $\Z$-basis of  $P_n$. 

Finally, if $u,v$
are permutations which belong to the infinite symmetric group 
 $\Sigma^{(\infty)},$  denote by $c_{u,v}^{w}$ the structural constants for
the multiplication of Schubert polynomials:
$$ {\s}_u~{\s}_v = \sum_{w \in \Sigma^{(\infty)}}c_{u,v}^{w}~{\s}_w.
$$
Then

$(\clubsuit)$~~~~~~~$c_{N \ast u,N \ast v}^{N \ast w} \ne 0$ for some integer 
$N \ge 1$ if and only if $c_{u,v}^{w} \ne 0.$
\end{con}
$(\maltese)$  We {\bf expect} that the formal power series
$$ \sum_{N \ge 1} c_{N \ast u,N \ast v}^{N \ast w}~t^N $$
is a rational function in $t$  (with the only possible pole at $t=1$ ??).
In other words, the function 
$N \longrightarrow c_{N \ast u,N \ast v}^{N \ast w}$
is a {\bf polynomial} in $N$ with rational (non--negative ??) coefficients.
\begin{prb} ( {\bf Generalized saturation problem for Kazhdan--Lusztig's 
polynomials} )

Let $u,w \in \Sigma_{n}$  be two permutations, denote by
$$ P_{u,w}(q)=1+\cdots+d(u,w)~q^{c(u,w)},~~d(u,w) \ne 0,$$
the corresponding Kazhdan--Lusztig polynomial \cite{KaL}.

$(\clubsuit)$ Prove (or disprove)  that

$(1)$ ~ $c(N \ast u, N \ast w)=N~c(u,w)$ ~~for any positive integer $N;$ 

$(2)$ ~$d(N \ast u, N \ast w)=1$ for some positive integer $N$ if and only if 
$d(u,w)=1.$
\end{prb}
The similar Problem can be stated for the Kazhdan--Lusztig polynomials 
corresponding to the affine symmetric group. 

However, we didn't extensively 
test Conjecture~6.27 and Problem~6.28 on a computer. \\

We want to end this Section by the following question and problem:
\begin{quest} ( {\bf A $q$-analog of the structural constants $c_{u,v}^{w}$} ) 

Does there exist a natural $q$-analog $c_{u,v}^{w}(q) \in \N~[q]$ of the
structural constants $c_{u,v}^{w},$ so that $c_{u,v}^{w}=c_{u,v}^{w}(1),$ 
which for the grassmannian permutations
$u,$ $v$ and $w$ coincides with the  $q$-analog  
$c_{\ld(u),\ld(v)}^{\ld(w)}(q)$  of the $LR$-numbers  ? 

Here $\ld(w)$ denotes the shape of a permutation $w,$ see \cite{Mac} for
a detailed account to the theory of Schubert polynomials. As for a definition 
of the  $q$-analog $c_{\ld,\mu}^{\nu}(q)$ of the $LR$-numbers, 
see e.g. \cite{CB}, \cite{LLT}.
\end{quest} 
\begin{prb} ( {\bf Define the polynomials $c_{\ld,\mu}^{\nu}(q)$ through the
geometry of Schubert varieties} )

Let $n \ge m$ be fixed positive integers, and $\ld,$ $\mu$ and 
$\nu$ be three partitions such that $\max(l(\ld),l(\mu),l(\nu)) \le m,$
~$\max(\ld_1,\mu_1,\nu_1) \le n,$ and $|\ld|+|\mu|=|\nu|.$

It is well--known that the $LR$-number $c_{\ld,\mu}^{\nu}$ counts the number 
of (isolated) points in the triple intersection ${\s}_{\ld} \cap {\s}_{\mu}
 \cap {\s}_{\nu^{*}}$ of the Schubert varieties ${\s}_{\ld},$ ${\s}_{\mu}$ 
and ${\s}_{\nu^{*}}$ in the Grassmannian variety $G(m,n+m),$ see e.g. 
\cite{Fu1} for the explanations of omitted notation, definitions and details.

$(\clubsuit)$  Find a geometric way to attach to each intersection point
$x \in {\s}_{\ld} \cap {\s}_{\mu} \cap {\s}_{\nu^{*}}$ an integer number 
$c(x)$ such that the generating function
$$ \sum_{x \in {\s}_{\ld} \cap {\s}_{\mu} \cap {\s}_{\nu^{*}}}~q^{c(x)}$$
coincides with the $LLT$ $q$-analog $c_{\ld,\mu}^{\nu}(q)$ of the Littlewood--
Richardson number $c_{\ld,\mu}^{\nu}.$
\end{prb}


Anatol N. Kirillov  \\
Research Institute for Mathematical Sciences \\ 
Kyoto University \\
Sakyo-ku, Kyoto 606-8502, Japan \\
e-mail: kirillov@kurims.kyoto-u.ac.jp \\
{\it URL: ~~~~http://www.kurims.kyoto-u.ac.jp/\kern-.05cm$\tilde{\quad}$\kern-.17cm kirillov }

\end{document}